\documentclass[aos]{imsart}
\usepackage{amsthm,amsmath,amsfonts,amssymb,stackengine,comment}
\usepackage{amsfonts,mathrsfs}
\usepackage{color}
\usepackage{enumerate}
\usepackage[textwidth=8em,textsize=scriptsize]{todonotes}
\usepackage{booktabs,float}
\usepackage{graphicx}
\usepackage{multirow}
\usepackage{float}
\usepackage{caption,subcaption}

\RequirePackage[authoryear]{natbib} %% uncomment this for author-year bibliography
\RequirePackage[colorlinks,citecolor=blue,urlcolor=blue]{hyperref}

\startlocaldefs
\numberwithin{equation}{section}
\theoremstyle{plain}
\newtheorem{theorem}{Theorem}[section]
\newtheorem{lemma}[theorem]{Lemma}
\theoremstyle{remark}

\newtheorem{assumption}{Assumption}
\newtheorem{remark}{Remark}

\newtheorem{corollary}{Corollary}

\def\cN{\mathcal{N}}
\def\cS{\mathcal{S}}

\endlocaldefs

\begin{document}
	
\begin{frontmatter}
		\title{
Robust Nonparametric Regression with \\
Deep Neural Networks}
		%\title{A sample article title with some additional note\thanksref{t1}}
		\runtitle{
Robust Deep Nonparametric Regression}
		%\thankstext{T1}{A sample additional note to the title.}
		
\begin{aug}
			%%%%%%%%%%%%%%%%%%%%%%%%%%%%%%%%%%%%%%%%%%%%%%
			%%Only one address is permitted per author. %%
			%%Only division, organization and e-mail is %%
			%%included in the address.                  %%
			%%Additional information can be included in %%
			%%the Acknowledgments section if necessary. %%
			%%%%%%%%%%%%%%%%%%%%%%%%%%%%%%%%%%%%%%%%%%%%%%
			\author[A]{\fnms{Guohao} \snm{Shen}\ead[label=e1,mark]{ghshen@link.cuhk.edu.hk}}\footnote{Guohao Shen and Yuling Jiao contributed equally to this work.}
			\author[B]{\fnms{Yuling} \snm{Jiao}\ead[label=e2, mark]{yulingjiaomath@whu.edu.cn}}$^*$		
			\author[C]{\fnms{Yuanyuan} \snm{Lin}\ead[label=e3,mark]{ylin@sta.cuhk.edu.hk}}
			\and
			\author[D]{\fnms{Jian} \snm{Huang}\ead[label=e4,mark]{jian-huang@uiowa.edu}}
			%%%%%%%%%%%%%%%%%%%%%%%%%%%%%%%%%%%%%%%%%%%%%%
			%% Addresses                                %%		%%%%%%%%%%%%%%%%%%%%%%%%%%%%%%%%%%%%%%%%%%%%%%
\address[A]{Department of Statistics, The Chinese University of Hong Kong,  Hong Kong, China  \printead{e1}}			
\address[B]{School of Mathematics and Statistics, Wuhan University, China  \printead{e2}}						
\address[C]{Department of Statistics, The Chinese University of Hong Kong,  Hong Kong, China  \printead{e3}}			
\address[D]{Department of Statistics and Actuarial Science, University of Iowa, Iowa, USA
  \printead{e4}}	
\end{aug}
		
\begin{abstract}
In this paper, we study the properties of robust nonparametric estimation using deep neural networks for regression models with heavy tailed error distributions. We establish the non-asymptotic error bounds for a class of robust nonparametric regression estimators using deep neural networks with ReLU activation under suitable smoothness conditions on the regression function and mild conditions on the error term. In particular, we only assume that the error distribution has a finite $p$-th moment with $p$ greater than one. We also show that the deep robust regression estimators are able to circumvent the curse of dimensionality when the distribution of the predictor is supported on an approximate lower-dimensional set. An important feature of our error bound is that, for ReLU neural networks with network width and network size (number of parameters) no more than the order of the square of the dimensionality $d$ of the predictor, our excess risk bounds depend sub-linearly on the $d$. Our assumption relaxes the exact manifold support assumption, which could be restrictive and unrealistic in practice. We also relax several crucial assumptions on the data distribution, the target regression function and the neural networks required in the recent literature. Our simulation studies demonstrate that
the robust methods can significantly outperform the least squares method when the errors have heavy-tailed distributions and illustrate that the choice of loss function is important in the context of deep nonparametric regression.
% advantages of using robust loss functions over the least squares loss function in the %context of deep nonparametric regression .

\end{abstract}

\iffalse		
		
\fi		
	\end{frontmatter}

\section{Introduction}
Consider a nonparametric regression model
\begin{equation}\label{model}
		Y=f_0(X)+\eta,
\end{equation}
where $Y \in \mathbb{R}$ is a response, $X \in \mathcal{X}\subseteq\mathbb{R}^d$ is a $d$-dimensional vector of predictors,
$f_0: %[0,1]^d
\mathcal{X} \to \mathbb{R}$ is an unknown regression function, $\eta$ is an unobservable error  independent of $X$, whose $p$-th moment is assumed finite.
A basic problem in statistics and machine learning is to estimate the unknown target regression function $f_0$ based on observations that are independent and identically distributed  (i.i.d.) copies of $(X, Y)$, $(X_i, Y_i), i=1, \ldots, n,$ where $n$ is the sample size. In this paper, we establish the non-asymptotic error bounds for
a class of robust nonparametric regression estimators using deep neural networks
under
%suitable smoothness conditions on the regression function and
mild conditions on the error term, allowing heavy-tailed error distributions.
Through simulation studies, we demonstrate that it is important to use a robust method for
heavy-tailed or contaminated data in nonparametric regression using deep neural networks.
We also show that the deep robust regression estimators are able to circumvent the curse of dimensionality when the predictor is supported on an approximate lower-dimensional set.

\begin{figure}[H]
	\centering
%	\begin{subfigure}{\textwidth}
%		\includegraphics[width=1\textwidth]{bumps_normal.png}
%	\end{subfigure}
	
%	\begin{subfigure}{1\textwidth}
%		\includegraphics[width=1\textwidth]{bumps_t.png}
%	\end{subfigure}
	
%	\begin{subfigure}{\textwidth}
		\includegraphics[width=\textwidth]{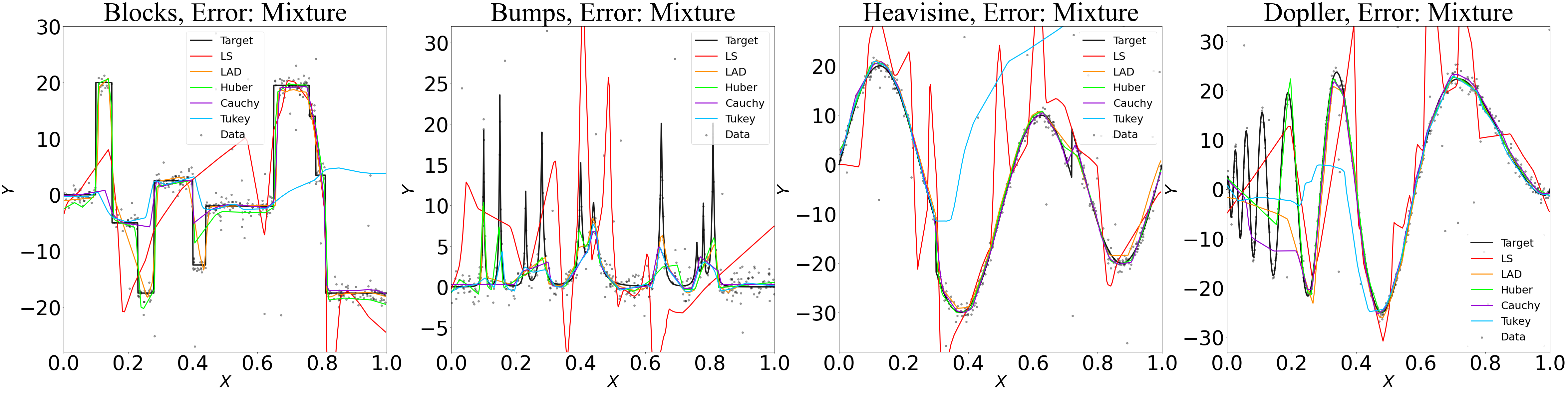}
%	\end{subfigure}
	
\caption{The fitting of four univariate models. The data points are displayed as grey dots, the target functions are depicted as black curves and the colored curves represent the estimated curves
%stand for predictions of the trained estimators
under different loss functions.
%From the top to the bottom, each row corresponds a case with a certain type of error: %$N(0,1),t(2)$ and
Red line: LS, Orange line: LAD, Green line: Huber, Purple line: Cauchy; Blue line: Tukey.
The data are generated based on the  univariate regression models in Section \ref{sim} with contaminated normal
% mixture
error.}
	\label{fig:0}
\end{figure}

To illustrate the importance of using robust methods for nonparametric regression using deep neural networks
in the presence of heavy-tailed errors, we look at the fitting of the regression functions  ``Blocks'', ``Bumps'', ``Heavisine'' and ``Dopller'' \citep{dj1994}, when error follows a contaminated normal distribution \citep{huber1964}, represented by  a mixture of two normals, $\eta \sim  0.8\mathcal{N}(0,1)+0.2 \mathcal{N}(0,10^4)$.
%$\eta\sim \xi \mathcal{N}(0,1)+(1-\xi) \mathcal{N}(0,10^4)$ with $\xi\sim Bernoulli(0.8)$.
The functional form of the models are given in Section \ref{sim}. Fig \ref{fig:0} shows the fitting results using LS, LAD, Huber, Cauchy and Tukey loss functions. We see that the least squares method breaks down in the presence of contaminated normal errors, while the robust methods using LAD, Huber and Cauchy loss functions yield reasonable fitting results. In particular, the Huber method works best among the methods considered in this example.

%\subsection{Literature review}
	
%	There has been a vast theoretical and applied statistical literature on nonparametric %regression over the last few decades. One reason for its growing popularity is that %nonparametric regression models are robust to model misspecification in the sense that %they offers a flexible and general modeling of the relationship between the covariates and %responses. There have been numerous methods to construct estimates,

There is a vast literature on nonparametric regression, many important methods have been developed.  Examples include
%the partitioning estimator \citep{tukey1947non,tukey1961curves},
 kernel regression \citep{nadaraya1964estimating, watson1964smooth},
 tree based regression \citep{cart1984},  local polynomial regression \citep{cleveland1979, fan1992, fan1993}, various spline-based regression methods
\citep{stone1982optimal,wahba1990, frieddman1991multivariate, greens1993}, regression in reproducing kernel Hilbert spaces \citep{wahba1990, scholkopf2018learning},
%the nearest neighbour estimator \cite{fix1985discriminatory},
among others.
%Among many estimation methods for nonparametric regression,
%a basic paradigm is global modeling based on the least squares objective function
Most of the existing methods have been developed using the least squares loss as the objective function over different function spaces.
% \citep{fridedman1991multivariate}.
See also the books \citep{statlearning2001, gyorfi2006distribution, npe2008} %scholkopf2018learning}
for detailed discussions on the general methodologies and theories on nonparametric regression and related topics.

%an extensive coverage of nonparametric regression methods and theories and the %references therein.
% (or least squares estimation).
%Considering the strong power of rapidly developing machine learning technologies, the %adoption of other more flexible nonparametric estimation techniques such as deep %learning offers a wide range of new methods.

Recently, several inspiring works studied the properties of nonparametric regression using deep neural networks
%standard  nonparametric least squares estimation has been studied
% in theory and application
\citep{bauer2019deep,schmidt2019deep,chen2019nonparametric,schmidt2020nonparametric,farrell2021deep}. In these works, convergence rates of the empirical risk minimizers based on the least squares objective function or a convex loss are established under suitable smoothness conditions on the underlying regression function and certain regularity conditions on the model. In particular,
%they assume
the response variable or the error term in the regression model is assumed to be bounded \citep{
gyorfi2006distribution,
farrell2021deep}, or sub-Gaussian \citep{bauer2019deep,chen2019nonparametric,schmidt2019deep,schmidt2020nonparametric}.

The aforementioned results are not applicable to the problems with heavy-tailed errors that commonly arise in practice, due to the sub-Gaussian assumption on the error distribution or
the boundedness assumption on the response variable. Also, it is well known that
nonparametric estimators based on least squares loss function are sensitive to outliers.
%thus they are not robust in analyzing data with heavy-tailed errors that are common in %economics and finance.
To address the non-robustness problem of the methods based on least squares and to deal with heavy- tailed errors, many robust regression methods had been developed in the context of linear regression models.
Important examples include the least absolute deviation (LAD)  \citep{bk1978lad},
%\citep{dodge2008concise},
%quantile regression \citep{koenker1978,
%koenker2001quantile,
%koenker2005},
Huber regression \citep{huber1973robust}, Tukey's biweight regression \citep{beaton1974fitting}, among others.
% etc, are studied in statistical literature.
Indeed, there is a rich literature on robust statistics going back many decades
\citep{anscombe1960RejectionOO,tukey1960, huber1964, tukey1975}.
For systematic discussions on the basic problems in robust statistics, we refer to \cite{huber2004robust} and \cite{hampel2011robust}.
Recently, there is a resurgent interest in developing polynomial time robust estimators achieving near optimal recovery guarantees for many parametric
estimation problems  including mean estimation, linear regression and generalized linear models,
where an analyst is given access to samples, in which a
fraction of the samples may have been adversarially corrupted, see for example, \cite{cherapanamjeri2020optimal}, \cite{jambulapati2021robust} and the references therein.

As every aspects of scientific research and the society are becoming digitized,
statistics and machine learning methods are applied increasingly to problems with more complex and noisier data. It is important that  data modeling methodologies and algorithms are robust to potentially nonnormal noises.
In particular, it has been recognized that robustness is an important and desirable property in applications using deep learning methods
\citep{belagiannis2015robust,wang2016studying, jiang2018mentornet,barron2019general}.  %But the theoretical aspects of robust deep nonparametric regression remain largely %untouched until the recent works \cite{lederer2020risk} and \cite{padilla2020quantile}.
However,  there has been only limited studies on the statistical guarantees of robust nonparametric regression using deep neural networks in the literature.
\cite{lederer2020risk} studied non-asymptotic error bounds for the empirical risk minimizer with robust loss functions under the assumption that the response has a finite second moment. \cite{padilla2020quantile} derived the convergence rates for
quantile regression with ReLU networks.  In both works,
the authors assumed that the true target function belongs to the function class of neural networks %used for estimation
without considering the approximation error. However, it is unrealistic to assume that %Unfortunately, this assumption is unrealistic in practice, since it is usually not known that
the underlying target regression function would be a member of the function class used in estimation. More recently, \cite{farrell2021deep} considered general Lipschitz loss functions under possibly mis-specified case (i.e., the true target function may not belong to the function class used for estimation),  but they assume that the response is bounded, which is restrictive and not compatible with heavy-tailed error models.
%assumption there could be restrictive in the application to  heavy-tailed data.
%		
%Despite the impressive advancements achieved in the existing works, there is still an
%There is another important drawback in the existing results on the
%excess risk bounds and the estimation
%error bounds for the nonparametric regression estimators using deep neural networks %based on the least squares, quantile check loss and other general loss functions. Specifically,
Also, the existing error bound results for the nonparametric regression estimators using deep neural networks contain prefactors that depend exponentially on the ambient dimension $d$ of the predictor  \citep{schmidt2020nonparametric, farrell2021deep, padilla2020quantile}. This adversely affects the quality of the error bounds, especially in the high-dimensional settings when $d$ is large \citep{ghorbani2020discussion}.  In particular, such type of error bounds lead to a sample complexity that depends exponentially on $d$. Such a sample size requirement is difficult to be met even for a moderately large $d$.

	In this paper, we study the properties of robust nonparametric regression using deep neural networks. We establish non-asymptotic upper bounds of the excess risk of the empirical risk minimizer for the nonparametric neural regression using feedforward neural networks with Rectified Linear Unit (ReLU) activation. Specifically, our main contributions are as follows:
	
	\begin{enumerate}[(i)]
		\item  For a heavy-tailed response $Y$ with a finite $p$-th moment for some $p>1$, we establish non-asymptotic bounds on the excess risk of robust nonparametric regression estimators using deep neural networks under continuous Lipschitz loss functions. Moreover, for ReLU neural networks with network width and network size (number of parameters) no more than $O(d^2)$, our excess risk bounds depend sub-linearly on the dimensionality $d$, instead of exponentially in terms of a factor $a^d$ (for some constant $a\geq 2$) as in the existing results.
% \citep{schmidt2020nonparametric,farrell2021deep}.
		
		\item We explicitly describe how the error bounds depend on the neural network parameters, including the width, the depth and the size of the network. The notion of network relative efficiency between two types of neural networks, defined as the ratio of the logarithms of the network sizes needed to achieve the optimal rate of convergence, is used to evaluate the relative merits of network structures.
		
		\item We show that the curse of dimensionality can be alleviated if the distribution of $X$ is assumed to be supported on an approximate low-dimensional manifold. This assumption relaxes the exact manifold support assumption, which could be restrictive and unrealistic in practice.
% \citep{schmidt2019deep,nakada2019adaptive,chen2019efficient}.
 Under such an approximate low-dimensional manifold support assumption, we show that the rate of convergence $n^{-(1-1/p)\alpha/(\alpha+d)}$ can be improved to $n^{-(1-1/p)\alpha/(\alpha+d_0)}$ for $d_0=O(d_{\mathcal{M}}\log(d))$, where $d_{\mathcal{M}}$ is the intrinsic dimension of the low-dimensional manifold and $\alpha\in(0,1]$ is the assumed order of the H\"older-continuity of the target $f^*$.
		
\item We relax several crucial assumptions on the data distribution, the target regression
function and the neural networks required in the recent literature.
%\citep{ bauer2019deep,  schmidt2019deep, schmidt2020nonparametric,farrell2021deep}.
First, we do not assume that the response $Y$ is bounded, instead only its $p$-th moment is finite. Second, we do not require the network to be sparse or have uniformly bounded weights and biases. Third, we largely relax the regularity conditions on the true target function $f^*$. In particular, we only require $f^*$ to be uniformly continuous, as our results are derived in terms of its modulus of continuity, which is  general enough to cover nearly all cases that $f^*$ belongs to different smooth classes of functions.

\item We carry our numerical experiments to evaluate the performance of robust
nonparametric regression estimators using LAD, Huber, Cauchy and Tukey loss functions and compare with the performance of the least squares estimator. The results demonstrate that, with heavy-tailed errors, the robust methods with LAD, Huber and Cauchy loss functions tend to significantly outperform the least squares method.
    %In multiple dimensional models, the robust methods also tend to perform better than %the least squares method even with normal errors.
		
\end{enumerate}
	
The remainder of the paper is organized  as follows.
	%Our main theoretical results are presented in sections \ref{sec2}-\ref{sec4}.
	In Section \ref{sec2} we introduce the setup of the problem,  present some popular robust loss functions, and describe the class of ReLU activated feedforward neural networks used in estimating the regression function.
	In Section \ref{sec3} we present a basic inequality for the excess risk in terms of the stochastic and approximation errors, and describe our approach in analyzing the two types of errors. In Section \ref{sec4} we provide sufficient conditions under which the robust neural regression estimator
	is shown to possess the basic consistency property. We also establish non-asymptotic error bounds for the neural regression estimator using deep feedforward neural networks. Moreover,  we show that the neural regression estimator can circumvent the curse of dimensionality if the data distribution is supported on an approximate low-dimensional manifold. In Section \ref{efficiency} we  present the results on how the error bounds depend on the network structures and consider a notion of network relative efficiency between two types of neural networks, defined as the ratio of the logarithms of the network sizes needed to achieve the optimal convergence rate. This can be used as a quantitative measure for evaluating the relative merits of different network structures. Extensive simulation studies are presented in section \ref{sim}. Concluding remarks are given in section \ref{conclusion}.
	
\section{Preliminaries}
	\label{sec2}
In this section, we introduce the robust nonparametric regressions, including the commonly-used robust loss functions,
% the empirical risk minimization,
and the definition of ReLU feedforward neural networks.
	
\subsection{Robust regression}
	
To accommodate heavy-tailed errors,  we consider a general robust paradigm for estimating $f_0$ by using a robust loss function $L(\cdot,\cdot):\mathbb{R}^2\to\mathbb{R}$.
For example, our results apply to the following robust loss functions.
\begin{itemize}
\item  Least absolute deviation (LAD) loss:
$L(a,y)=\vert a-y\vert, \  (a, y) \in \mathbb{R}^2.$

\item  Quantile loss:
$L(a,y)=\rho_\tau(a-y), \  (a, y) \in \mathbb{R}^2,$
where $\rho_\tau(x)=\tau x$ if $x\geq0$ and $\rho_\tau(x)=(\tau-1)x$ if $x<0$ for some $\tau \in(0,1)$.

\item Huber loss:
 $L(a,y)=h_\zeta(a-y), \ \ (a, y) \in \mathbb{R}^2,  \text{ for some } \zeta > 0,
 $
 where
 %\[h_\zeta(t)  =\left\{ \begin{array}{ll}
%  t^2 & \text{ if } \vert a\vert\leq\zeta  \\
%  \zeta\vert t \vert-\zeta^2/2 & \text{ otherwise}
%  \end{array} \right.
%  \]
where $h_\zeta(x)=x^2/2$ if $\vert x\vert\leq\zeta$ and $h_\zeta(x)=\zeta\vert x\vert-\zeta^2/2$ otherwise.

 \item Cauchy loss:
  $L(a,y)=\log\{1+\kappa^2(a-y)^2\}, \ (a, y) \in \mathbb{R}^2,$ for some $  \kappa > 0.$
  %\in(0,\infty),$$

 \item Tukey's biweight loss:
 $L(a,y)=T_t(a-y), \ (a, y) \in \mathbb{R}^2,  \text{ for some } t>0,
 $
 where $T_t(x)=t^2[1-\{1-(x/t)^2\}^3]/6$ if $\vert a\vert\leq t$ and $T_t(x)=t^2/6$ otherwise.
\end{itemize}
These popular robust loss functions are all continuous and Lipschitz in both of its arguments. Throughout this paper, we consider a general and broad class of loss functions that are continuous and Lipschitz in both of its arguments, as stated in the following assumption.
	
\begin{figure}
		\includegraphics[width=\textwidth]{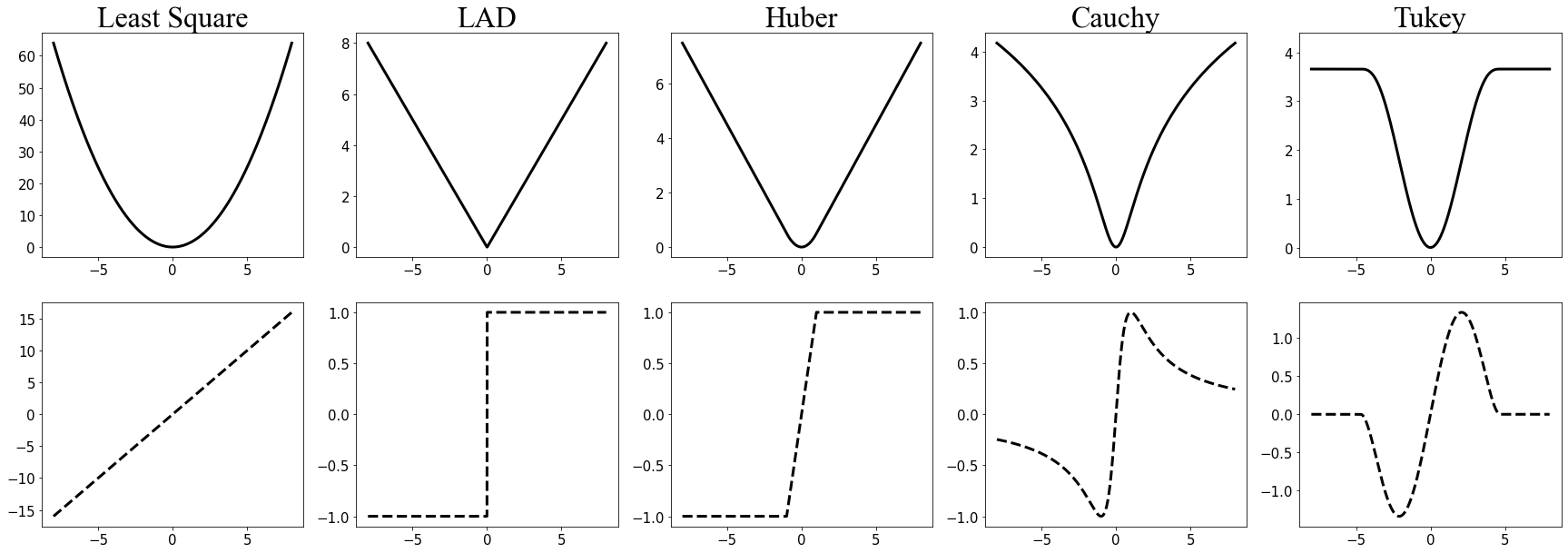}
		\caption{Different loss functions (solid lines in the top row) and their corresponding derivatives (dashed lines in the bottom row). For depicted Huber loss function, $\zeta=1$; Cauchy loss function, $\kappa=1$; Tukey loss function, $t=4.685$.}
		\label{fig1}
			\end{figure}

	\begin{assumption}\label{A1}
		The loss function $L(\cdot,\cdot):\mathbb{R}^2\to\mathbb{R}$ is continuous, and $L(a,y)=0$ if $a=y$ for $(a, y)\in\mathbb{R}^2$. Besides, $L$ is $\lambda_L$Lipschitz in both of its arguments, i.e.
		\begin{align*}
			\vert L(a_1,\cdot)-L(a_2,\cdot)\vert&\leq \lambda_L\vert a_1-a_2\vert,\\
			\vert L(\cdot,y_1)-L(\cdot,y_2)\vert&\leq \lambda_L\vert y_1-y_2\vert,
		\end{align*}
		for any $a_1,a _2,y_1,y_2\in\mathbb{R}$.
	\end{assumption}
	The aforementioned robust loss functions all satisfy Assumption \ref{A1}. To be specific, for the LAD loss, $\lambda_L=1$; for the quantile loss, $\lambda_L=\max\{\tau,1-\tau\}$; for the Huber loss, $\lambda_L=\zeta$; for the Cauchy loss, $\lambda_L=\kappa$; for the Tukey's biweight loss, $\lambda_L=16t/25\sqrt{5}$.
In particular, the LAD, the quantile check loss and the Huber loss functions can be
recast into
a more simplified form, i.e. $L(a,y)=\psi(a-y)$ for for some convex or Lipschitz function  $\psi$ with $\psi(0)=0$. A detailed summary is given in Table \ref{tab:1}, and plots of the above loss functions (except quantile loss function) and their derivatives in comparison with least squares are displayed in Figure \ref{fig1}.
	
\begin{table}[H]
		\centering
		\caption{Summary of different robust loss functions}
		\label{tab:1}
		\resizebox{\textwidth}{!}{\scriptsize %
			\begin{tabular}{cccccc}
				\hline
				& LAD  & Quantile        & Huber  & Cauchy & Tukey \\ \hline
				Hyper parameter & N/A & $\tau \in(0,1)$ & $\zeta\in(0,\infty)$ & $\kappa\in(0,\infty)$  & $t\in(0,\infty)$\\
				$\lambda_L$ & 1    & $\max\{\tau,1-\tau\}$ & $\zeta$ & $\kappa$  & $16t/25\sqrt{5}$  \\
				Continuous  & TRUE & TRUE            & TRUE   & TRUE  & TRUE \\
				Convex      & TRUE & TRUE            & TRUE   & FALSE & FALSE\\
				Differentiable & FALSE & FALSE& TRUE & TRUE & TRUE \\ \bottomrule
			\end{tabular}%
		}
		
		{Note: "N/A" stands for not applicable.}
	\end{table}
	
	For any (random) function $f$, let $Z\equiv (X,Y)$ be a random vector independent of  $f$.
	For any loss function $L$ satisfying Assumption \ref{A1},
	the $L$ risk is defined by $\mathcal{R}(f)=\mathbb{E}_{Z}L(f(X),Y)$. At the population level, an estimation procedure is to find a measurable function $f^*: \mathbb{R}^d\to \mathbb{R}$ satisfying
	\begin{equation}\label{target}
		f^* :=\arg\min_{f} \mathcal{R}(f) =\arg\min_{f}\mathbb{E}_{Z}L(f(X),Y).
	\end{equation}
	%Under the framework of least squares estimation, $f^*$ is also called the {\it Bayes estimator},  	
	%since by conditional on $X$, the minimizer $f^*$ can also be defined point-wisely:
	%$$f^*(x) := \arg\min_f\mathbb{E}_{\eta|X=x}\Vert f_0(x)+\eta-f(x)\Vert^2, \quad \forall x\in\mathcal{X}.$$
	Under proper assumptions on $L$, $X$ and $\eta$,  the true regression function $f_0$ in model (\ref{model}) is the optimal solution $f^*$ on $\mathcal{X}$. For instance, when $L(f(X),Y)=\vert f(X)-Y\vert$ and under the assumption that the conditional median of $\eta$ given $X$ is $0$,  $f^*$
in (\ref{target})	
	is precisely the true regression function $f_0$ on $\mathcal{X}$.
	 For the rest of the paper, under a general robust nonparametric estimation framework
	without specific assumptions on $L,X$ or $\eta$ to ensure $f^*=f_0$ on $\mathcal{X}$, we would focus on estimating $f^*$ instead of $f_0$.

	In real applications, the distribution of $(X,Y)$ is typically unknown and only a random sample $S \equiv \{(X_i,Y_i)\}_{i=1}^n$ with size $n$ is available.
	Define the empirical risk of $f$ on the sample $S$ as
		\begin{equation}
		\label{er1}
		\mathcal{R}_n(f)=\frac{1}{n}\sum_{i=1}^{n}L(f(X_i),Y_i).
	\end{equation}	
	Based on the data, a primary goal is to construct estimators of $f^*$ within a certain class of functions $\mathcal{F}_n$ by minimizing the empirical risk.
	%, such that the empirical risk is small.
	Thus, the empirical risk minimizer (ERM) is defined by
	\begin{equation}\label{erm}
		\hat{f}_n\in\arg\min_{f\in\mathcal{F}_n}\mathcal{R}_n(f).
	\end{equation}
	In this paper, we set $\mathcal{F}_n$  to be a function class consisting of
	feedforward neural networks.
	For any estimator  $\hat{f}_n$,  we evaluate its quality via its  \textit{excess risk},  which is
	defined as the difference between the $L$
	risks of $\hat{f}_n$ and $f^*$:
	%	denoted by $\Vert\hat{f}_n-f_0\Vert^2_{L^2(\nu)}$,	 is of central concern. 	
	%	Specifically,
	\begin{align*}
		%		\Vert\hat{f}_n-f_0\Vert^2_{L^2(\nu)} &=
		\mathcal{R}(\hat{f}_n)-\mathcal{R}(f^*) = \mathbb{E}_{Z}L(\hat{f}_n(X),Y)-\mathbb{E}_{Z}L(f^*(X),Y). \nonumber
		% \\  & = \mathbb{E}_{Z_{0}}\vert \hat{f}_n(X_0)-f_0(X_0)\vert^2,  \nonumber
	\end{align*}
	%Let  $\Vert\hat{f}_n-f_0\Vert^2_n :=\sum_{i=1}^n\vert\hat{f}_n(X_i)-f_0(X_i)\vert^2/n$
	%	be the empirical counterpart of the excess risk $ \Vert\hat{f}_n-f_0\Vert^2_{L^2(\nu)}$.
	%	Intuitively,
	A good estimator  $\hat{f}_n$ shall have a small excess risk
	% keep
	$\mathcal{R}(\hat{f}_n)-\mathcal{R}(f^*) $.
	%or its empirical counterpart $\Vert\hat{f}_n-f_0\Vert^2_n$
	%as small as possible.
	Thereafter,  we shall focus on deriving the non-asymptotic upper bounds of the excess risk $\mathcal{R}(\hat{f}_n)-\mathcal{R}(f^*) $. Note that $\mathcal{R}(\hat{f}_n)$ is  random as $\hat{f}_n$ is produced based on the sample $S$. Thus,  we will derive bounds for the expected excess risk $\mathbb{E}_S\big\{\mathcal{R}(\hat{f}_n)-\mathcal{R}(f^*)\big\}$ as well.
	
	\subsection{ReLU Feedforward neural networks}
	In recent years,  deep neural network modeling has achieved impressive successes in many applications.
	% An important reason is that
	Also, neural network functions have proven to be
	%been recognized as a promising and  a powerful
	an effective tool to approximate high-dimensional functions.
	%The design of deep feedforward neural nets is a crucial  of such a revolution.
	We consider regression function estimators based on the feedforward neural networks with rectified linear unit (ReLU) activation function.
	Specifically, we set the function class $\mathcal{F}_n$ %is set
	to be $\mathcal{F}_{\mathcal{D},\mathcal{W}, \mathcal{U},\mathcal{S},\mathcal{B}}$, a class of  feedforward neural networks $f_\phi: \mathbb{R}^d \to \mathbb{R} $ with parameter $\phi$, depth $\mathcal{D}$, width $\mathcal{W}$, size $\mathcal{S}$, number of  neurons $\mathcal{U}$ and $f_{\phi}$ satisfying $\Vert f_\phi\Vert_\infty\leq\mathcal{B}$ for some $0 <B < \infty$, where
	$\Vert f \Vert_\infty$ is the sup-norm of a function $f$.
	Note that the network parameters may depend on the sample size $n$, but  the dependence is omitted in the notation for simplicity.
	A brief description of the feedforward neural networks are given below.
	
	We begin  with the multi-layer perceptron (MLP), an important subclass of feedforward neural networks. The architecture of a MLP can be expressed as a composition of a series of functions
	\[
	f_\phi(x)=\mathcal{L}_\mathcal{D}\circ\sigma\circ\mathcal{L}_{\mathcal{D}-1}
	\circ\sigma\circ\cdots\circ\sigma\circ\mathcal{L}_{1}\circ\sigma\circ\mathcal{L}_0(x),\  x\in \mathbb{R}^d,
	\]
	where $\sigma(x)=\max(0, x)$ is the rectified linear unit (ReLU) activation function (defined for each component of $x$ if $x$ is a vector) and
	$$\mathcal{L}_{i}(x)=W_ix+b_i,\quad i=0,1,\ldots,\mathcal{D},$$
	where $W_i\in\mathbb{R}^{d_{i+1}\times d_i}$  is a weight matrix,  $d_i$ is the width (the number of neurons or computational units) of the $i$-th layer, and $b_i\in\mathbb{R}^{d_{i+1}}$ is the bias vector in the $i$-th linear transformation $\mathcal{L}_i$.
	The input data consisting of predictor values $X$ is the first layer and the output
	%of $Y$
	is the last layer.
	Such a network $f_\phi$ has $\mathcal{D}$ hidden layers and $(\mathcal{D}+1)$ layers in total.
	We use a $(\mathcal{D}+1)$-vector $(d_0,d_1,\ldots,d_\mathcal{D})^\top$ to describe the width of each layer; particularly,   $d_0=d$ is the dimension of the input $X$ and $d_\mathcal{D}=1$ is the dimension of the response $Y$ in model (\ref{model}). The width $\mathcal{W}$ is defined as the maximum width of hidden layers, i.e.,
	$\mathcal{W}=\max\{d_1,...,d_\mathcal{D}\}$; 	the size $\mathcal{S}$ is defined as the total number of parameters in the network $f_\phi$, i.e $\mathcal{S}=\sum_{i=0}^\mathcal{D}\{d_{i+1}\times(d_i+1)\}$; 	
	the number of neurons $\mathcal{U}$ is defined as the number of computational units in hidden layers, i.e.,
	$\mathcal{U}=\sum_{i=1}^\mathcal{D} d_i$.  Note that the neurons in consecutive layers of a MLP are connected to each other via linear transformation matrices $W_i\in\mathbb{R}^{d_{i+1}\times d_i}$,  $i=0,1,\ldots,\mathcal{D}$.   In other words,
	an MLP is fully connected between consecutive layers and has no other connections.
	For an MLP  $\mathcal{F}_{\mathcal{D},\mathcal{U},\mathcal{W},\mathcal{S},\mathcal{B}}$,
	its parameters satisfy the simple
	relationship
	% \begin{equation}
	%\label{size0}
	\[
	\max\{\mathcal{W},\mathcal{D}\}\leq\mathcal{S}\leq\mathcal{W}(d+1)
	+(\mathcal{W}^2+\mathcal{W})(\mathcal{D}-1)+\mathcal{W}+1
	=O(\mathcal{W}^2\mathcal{D}).
	\]
	%\end{equation}
	%	For general feedforward neural networks, the layers need not be connected in a chain, and skip connections are allowed. For example, some neurons in layer  $i$ can be directly connected to neurons in  layer $(i + 2)$  or higher. In such structured networks,  a neuron  is said to be in layer $i=1,2,\ldots,\mathcal{D}$ if it has a predecessor in layer $i-1$ and no predecessor in any layer $i^\prime\geq i$.  ??Given a citation for reference.??
	
	Different from multilayer perceptron, a general feedforward neural network may not be fully connected. 	
	For such a network, % general feedforward neural network,
	each neuron in layer $i$ may be connected to only a small subset of neurons in layer $i+1$.
	The total number of parameters $\mathcal{S}\leq\sum_{i=0}^\mathcal{D}\{d_{i+1}\times(d_i+1)\}$ is reduced and the computational cost required to evaluate the network will also be reduced.
	
	%Though the multi-layer perceptron are commonly used in practice due to its simplicity,  our theoretical results cover general feedforward neural networks. Moreover, our results for ReLU networks can be extended to networks with piecewise-linear activation functions without further difficulty; see \cite{yarotsky2017error} and \citet{bartlett2019nearly}.

	For notational simplicity,  we write $\mathcal{F}_\phi=\mathcal{F}_{\mathcal{D},\mathcal{W},\mathcal{U},\mathcal{S},\mathcal{B}}$ for short to denote the class of feedforward neural networks with parameters $\mathcal{D},\mathcal{W},\mathcal{U},\mathcal{S}$ and $\mathcal{B}$.
	In the following, we will present our main results: nonasymptotic upper bounds on the excess risk for general feedforward neural networks with piecewise linear activation function.
	
	\section{Basic error analysis}
	\label{sec3}
	
	%\subsection{A basic inequality}
	%\label{sec3.1}
For general robust estimation,
	we start by introducing a generic upper bound on the excess risk of the empirical risk minimizer $\hat f_n$  defined in (\ref{erm})
	with a general function class $\mathcal{F}_n$.
	To this end, for a general loss function $L$ and any estimator $f$ belonging to
	 $\mathcal{F}_n$, the excess risk of $f$
	can be written as
	$$\mathcal{R}(f)-\mathcal{R}(f^*)=\left\{\mathcal{R}(f)-\inf_{f\in\mathcal{F}_n} \mathcal{R}(f)\right\}+ \left\{\inf_{f\in\mathcal{F}_n}\mathcal{R}(f)-\mathcal{R}(f^*) \right\}.$$
The first term on the right hand side is the so-called {\it stochastic error}, which measures the difference between the error of $f$
and the best one in the class  $\mathcal{F}_n$; the second term is the commonly-known {\it approximation error}, which measures
how well the function class $\mathcal{F}_n$ with respect to the loss function $L$ can be used to approximate the function $f^*$. 	
	
It follows  from Lemma 3.1  in
	 \cite{jiao2021deep} that,
	for any random sample $S=\{(X_i, Y_i)_{i=1}^n\}$, %$S=\{Z_i=(X_i^\top,Y_i)^\top\}_{i=1}^n$,
		the excess risk of  the empirical risk minimizer $\hat f_n$   satisfies
		\begin{align}
		\label{uberm}
			\mathcal{R}(\hat{f}_n)-\mathcal{R}(f^*)\leq&  2\sup_{f\in\mathcal{F}_n}\vert \mathcal{R}(f)-\mathcal{R}_{n}(f)\vert+\inf_{f\in\mathcal{F}_n} \mathcal{R}({f})-\mathcal{R}(f^*).	
	\end{align}
	The upper bound of the excess risk of $\hat f_n$ in (\ref{uberm}) does not  depend on $\hat f_n$ any more, but the function class $\mathcal{F}_n$, the loss function $L$ and the random sample $S$.		
	%The excess risk of $\hat f_n$   is upper bounded  by the sum of two terms: $2\sup_{f\in\mathcal{F}_n}\vert \mathcal{R}(f)-\mathcal{R}_n(f)\vert$ and the approximation error $\inf_{f\in\mathcal{F}_n}\mathcal{R}(f)- \mathcal{R}(f^*)$.
	The first term  in the upper bound $2\sup_{f\in\mathcal{F}_n}\vert \mathcal{R}(f)-\mathcal{R}_n(f)\vert$ can be properly handled by the classical empirical process theory \citep{vw1996, anthony1999, bartlett2019nearly}, whose   upper bound is determined by the complexity of the function class $\mathcal{F}_n$.
	% measured by the covering numbers.
	The second term $\inf_{f\in\mathcal{F}_n} \mathcal{R}({f})-\mathcal{R}(f^*)$ is the approximation error of the function class $\mathcal{F}_n$ for $f^*$ under loss $L$.
	The neural networks have been shown to have strong  power to approximate
 high-dimensional functions in both theory and practice.  Important recent results on the approximation theory of neural networks
 can be found in  \citet{yarotsky2017error, yarotsky2018optimal,shen2019nonlinear, shen2019deep}, among others.
	In order to apply those novel approximation theories  to our problem, we need to write $\inf_{f\in\mathcal{F}_n} \mathcal{R}({f})-\mathcal{R}(f^*)$ into
	a form of $\inf_{f\in\mathcal{F}_n} \Vert f-f^*\Vert$ for some functional norm $\Vert\cdot\Vert$.

	\subsection{Stochastic error}
	\label{sec3.2}
	We now proceed to analyze  the stochastic error of ERM $\hat{f}_\phi$ constructed by the feedforward neural networks and provide an upper bound on the expected excess risk.  Here,  the ERM $\hat{f}_\phi$  is defined by
	\begin{equation}\label{erm_nn}
		\hat{f}_\phi\in\arg\min_{f\in\mathcal{F}_\phi} \mathcal{R}_n(f),
	\end{equation}
	where $\mathcal{F}_\phi=\mathcal{F}_{\mathcal{D},\mathcal{W},\mathcal{U},\mathcal{S},\mathcal{B}}$ is the class of feedforward neural networks with parameters $\mathcal{D},\mathcal{W},\mathcal{U},\mathcal{S}$ and $\mathcal{B}$.
	The stochastic error of $\hat{f}_\phi$ can be upper bounded by the uniform results (on the class $\mathcal{F}_\phi$) in empirical process theory. In light of this, we next introduce the conceptions of  metric entropy of a class of functions which are needed  in stochastic error bounds.

	Denote the pseudo dimension for a class  of functions $\mathcal{F}$: $\mathcal{X}\to \mathbb{R}$ by $\text{Pdim}(\mathcal{F}),$, defined as the largest integer $m$ for which there exists $(x_1,\ldots,x_m,y_1,\ldots,y_m)\in\mathcal{X}^m\times\mathbb{R}^m$ such that for any $(b_1,\ldots,b_m)\in\{0,1\}^m$ there exists $f\in\mathcal{F}$ such that $\forall i:f(x_i)>y_i\iff b_i=1$ \citep{anthony1999, bartlett2019nearly}.
	Pseudo dimension is a natural measure of the complexity of a function class, e.g., a class of real-valued functions generated by neural networks.
	Especially, if $\mathcal{F}$ is the class of functions generated by a neural network with a fixed architecture and fixed activation functions,
	it follows from Theorem 14.1 in \cite{anthony1999} that $\text{Pdim}(\mathcal{F})=\text{VCdim}(\mathcal{F})$, where $\text{VCdim}(\mathcal{F})$ is the VC dimension of $\mathcal{F}$.
	In our results, it is required that the sample size $n$ is larger than the pseudo dimension of the class of neural networks under consideration.
	
	For a given sequence $x=(x_1,\ldots,x_n)\in\mathcal{X}^n,$
	%\subseteq\mathbb{R}^n$,
	we define  $\mathcal{F}_\phi|_x=\{(f(x_1),\ldots,f(x_n):f\in\mathcal{F}_\phi\}$,  a  subset of $\mathbb{R}^{n}$. For a positive number $\delta$, let $\mathcal{N}(\delta,\Vert\cdot\Vert_\infty,\mathcal{F}_\phi|_x)$ be the covering number of $\mathcal{F}_\phi|_x$ under the norm $\Vert\cdot\Vert_\infty$ with radius $\delta$.
	And the uniform covering number
	$\mathcal{N}_n(\delta,\Vert\cdot\Vert_\infty,\mathcal{F}_\phi)$ is defined by
	 	\begin{equation}
		\label{ucover}
		\mathcal{N}_n(\delta,\Vert\cdot\Vert_\infty,\mathcal{F}_\phi)=
		\max\{\mathcal{N}(\delta,\Vert\cdot\Vert_\infty,\mathcal{F}_\phi|_x):x\in\mathcal{X}\},
		\end{equation}
			 the maximum of the covering number $\mathcal{N}(\delta,\Vert\cdot\Vert_\infty,\mathcal{F}_\phi|_x)$  over all $x\in\mathcal{X}$.

	 We are now ready to present  a lemma concerning the upper bound of the stochastic error  under a heavy-tailed assumption of the response $Y$.

	\begin{assumption}
		\label{A2}
		The response variable $Y$ has finite $p$-th moment with $p>1$, i.e., there exists a constant $M>0$ such that $\mathbb{E}\vert Y\vert^p\leq M<\infty$.
	\end{assumption}
	
	\begin{lemma}\label{lemma1}
		Consider the $d$-variate nonparametric regression model in (\ref{model}) with an unknown regression function $f_0$ and target function $f^*$ defined in (\ref{target}).  Assume that Assumptions \ref{A1}-\ref{A2} hold, and $ \Vert f^*\Vert_\infty\leq \mathcal{B}$ for $\mathcal{B}\geq 1$. Let $\mathcal{F}_\phi=\mathcal{F}_{\mathcal{D},\mathcal{W},\mathcal{U},\mathcal{S},\mathcal{B}}$ be the class of feedforward neural networks with a continuous piecewise-linear activation function of finite pieces and $\hat{f}_\phi$ be the ERM defined in (\ref{erm_nn}). Then,  for $ \text{Pdim}(\mathcal{F}_\phi)\le 2n$,
		\begin{equation} \label{bound5a}
			%	\mathcal{R}(\hat{f}_\phi)
			\sup_{f\in\mathcal{F}_\phi}\vert \mathcal{R}(f)-\mathcal{R}_{n}(f)\vert \leq c_0\frac{\lambda_L\mathcal{B}}{n^{1-1/p}} \log\mathcal{N}_{2n}(n^{-1},\Vert \cdot\Vert_\infty,\mathcal{F}_\phi),
		\end{equation}
		where $c_0>0$ is a constant independent of $n,d,\lambda_L,\mathcal{B},\mathcal{S},\mathcal{W}$ and $\mathcal{D}$.
		In addition,
		\begin{equation} \label{oracle}
			\mathbb{E}\big\{\mathcal{R}(\hat{f}_\phi)-\mathcal{R}(f^*)\big\}\leq C_0\frac{\lambda_L\mathcal{B}\mathcal{S}\mathcal{D}\log(\mathcal{S})\log(n)}{n^{1-1/p}}+ 2\inf_{f\in\mathcal{F}_\phi}\big\{\mathcal{R}(f)-\mathcal{R}(f^*)\big\},
		\end{equation}
		where $C_0>0$ is a constant independent of $n,d,\lambda_L,\mathcal{B},\mathcal{S},\mathcal{W}$ and $\mathcal{D}$.
	\end{lemma}
	%\todo{take a look at this lemma}
	%	The proof of Lemma \ref{lemma1} is deferred to Appendix.
  Lemma \ref{lemma1} indicates that, the stochastic error is bounded by a term determined by the metric entropy of $\mathcal{F}_\phi$
	in (\ref{bound5a}).
	%$\mathcal{R}(\hat{f}_\phi)=\mathbb{E}\Vert \hat{f}_\phi-f_0\Vert^2_{L^2(\nu)}$
	The metric entropy is measured by the covering number  of $\mathcal{F}_\phi$.
	% with some radius under certain norm.
	To obtain the bound for excess risk in (\ref{oracle}),
	we next bound the covering number of $\mathcal{F}_\phi$ by its pseudo dimension (VC dimension), and apply the result
	in \cite{bartlett2019nearly} that  ${\rm Pdim}(\mathcal{F}_\phi)=O(\mathcal{S}\mathcal{D}\log(\mathcal{S}))$,
	% where ${\rm Pdim}(\mathcal{F}_\phi)$ denotes the pseudo dimension of $\mathcal{F}_\phi$.
	which gives the upper bound for  the expected excess risk as the sum of the stochastic error and  the approximation error of $\mathcal{F}_\phi$ to $f^*$ in (\ref{oracle}). A few important remarks are in order.

	\begin{remark}
		If $Y$ is assumed to be sub-exponentially distributed, i.e., there exists a constant $\sigma_Y>0$ such that $\mathbb{E}\exp(\sigma_Y\vert Y\vert)<\infty$,  the denominator $n^{1-1/p}$ in (\ref{bound5a}) and (\ref{oracle}) can be improved to be $n$, corresponding to the case that $p=+\infty$.
	Notably, the finite $p$-th moment condition and sub-exponential tail condition have intrinsically different impacts on the resulting error bounds in nonparametric regressions.
	\end{remark}

	\begin{remark}
 The non-asymptotic stochastic error bounds in \cite{lederer2020risk} are derived under finite second-moment  assumption of both  covariates and response. Our result covers their result  as a special case by letting $p=2$ in (\ref{bound5a})  in Lemma \ref{lemma1}. Assuming finite  $p$-th moment
 ($p>2$)  would give better rate. Furthermore, our result is more general than that of \cite{lederer2020risk}, since only the well-specified case, i.e. $f^*\in\mathcal{F}_n$,
 is considered and  no approximation error is involved in
 \cite{lederer2020risk}.
	\end{remark}
	
	\begin{remark}
		Upper bounds for least-squares estimator with neural nonparametric regression were studied by \cite{gyorfi2006distribution} and \cite{farrell2021deep} under the assumption that the response $Y$ is bounded.
%Without  the boundedness assumption on $Y$,
\cite{bauer2019deep} and \cite{schmidt2020nonparametric} derived error bounds for the least squares nonparametric regression estimator with deep neural networks under the assumption that %sub-exponentially distributed
the response $Y$ is sub-exponential.
	\end{remark}

	\begin{remark}
	 \cite{han2019convergence} studied general nonparametric least squares regression with heavy-tailed errors whose $p$-th moment ($p\geq1$)
	 is assumed finite and independent of the covariates.  Under the well-specified scenario that $f^*\in\mathcal{F}_n$ they studied the convergence rate for nonparametric least squares estimator.
	They showed that the rate of convergence with respect to the sample size $n$ is $O(\max\{n^{-1/(2+\gamma)},n^{-1/2+1/(2p)}\})$ when the model satisfies a standard entropy condition with exponent $\gamma\in(0,2)$, i.e $\log\mathcal{N}(\delta,\Vert\cdot\Vert_{L^2(\nu)},\mathcal{F}_n)\leq C\delta^{-\gamma}$ for some constant $C>0$, where $\nu$ is the probability measure of $X$.
	Their results shed light on the fact that,
	when the order of moment  $p\geq 1+2/\gamma$, the $L_2$ risk of the nonparametric least squares estimator converges at the same
		rate as the case with Gaussian errors; otherwise for $p< 1+2/\gamma$,
			%	 for any $\gamma$ satisfying $p< 1+2/\gamma$,
				 the convergence rate of the $L_2$ risk of least squares estimator can be strictly slower than other robust estimators.
				 On the other hand,
	the entropy condition with exponent $\gamma\in(0,2)$ on the function class
may not be satisfied in some scenarios,
%could be  demanding under some scenario,
especially when the dimensionality $d$ is high.
	For instance, for the H\"older class $\mathcal{F}=C^{\alpha}_1(\mathcal{X})$ ($f:\mathcal{X}\to\mathbb{R}$ with $\Vert f\Vert_\alpha\leq1$ for all $f\in\mathcal{F}$), its entropy satisfies $\log \mathcal{N}(\delta,\Vert\cdot\Vert_{L^2(\nu)},\mathcal{F})\leq K_{\alpha,d,\mathcal{X}} \delta^{-d/\alpha}$ for some constant $K_{\alpha,d,\mathcal{X}}$; for Lipschitz functions with Lipschitz constant $L$ on $[a,b]^d$, its entropy satisfies $\log \mathcal{N}(\delta,\Vert\cdot\Vert_{L^2(\nu)},\mathcal{F})\leq C \{(b-a)L\}^{-d}\delta^{-d}$. Besides, for VC-class with envelope function $F$ and $r\ge1$, and any probability measure $Q$, its entropy satisfies $\mathcal{N}(\delta\Vert F\Vert_{L^r(Q)},\Vert\cdot\Vert_{L^r(Q)},\mathcal{F})\leq C V(\mathcal{F})(16e)^{V(\mathcal{F})}\delta^{-r\{V(\mathcal{F})-1\}}$ for any $0<\delta<1$.  In some situations above, the entropy condition is not satisfied.
%Some situations above satisfy the ``entropy condition'', but some do not.
\end{remark}
	
	\subsection{Approximation error}
	\label{sec3.3}
	The approximation error depends on the target function $f^*$ and the function class $\mathcal{F}_\phi=\mathcal{F}_{\mathcal{D},\mathcal{W},\mathcal{U},\mathcal{S},\mathcal{B}}$. %It is also known that
	We assume that the distribution of the predictor $X$ is supported on a bounded set, and for simplicity, we assume this bounded set to be $[0,1]^d$. Thereafter, we consider the target functions $f^*$ defined on $\mathcal{X}=[0,1]^d$.
	
	To derive an upper bound of the approximation error $\inf_{f\in\mathcal{F}_\phi}\{R(f) -R^\phi(f^*)\}$ based on existing approximation theories, we need to fill the gap between $\inf_{f\in\mathcal{F}_\phi}\{R(f) -R^\phi(f^*)\}$ and $\inf_{f\in\mathcal{F}_{\phi}}\Vert f-f^*\Vert$.
Indeed, this can be done under Assumption \ref{A1} as follows.

	\begin{lemma}\label{lemma2}
		Assume that Assumption \ref{A1} holds. Let $f^*$ be the target function defined in (\ref{target}) and $\mathcal{R}(f^*)$ be its risk, then we have
		$$\inf_{f\in\mathcal{F}_\phi}\{R(f) -R(f^*)\}\leq \lambda_L\inf_{f\in\mathcal{F}_{\phi}}\mathbb{E}\vert f(X)-f^*(X)\vert=:\lambda_L\inf_{f\in\mathcal{F}_{\phi}}\Vert f-f^*\Vert_{L^1(\nu)},$$
		where $\nu$ denotes the marginal probability measure of $X$ and $\mathcal{F}_\phi=\mathcal{F}_{\mathcal{D},\mathcal{W},\mathcal{U},\mathcal{S},\mathcal{B}}$ denotes the class of feedforward neural networks with parameters $\mathcal{D},\mathcal{W},\mathcal{U},\mathcal{S}$ and $\mathcal{B}$.
	\end{lemma}
	As a consequence of Lemma \ref{lemma2}, we only need to give upper bounds on the approximation error $\inf_{f\in\mathcal{F}_{\phi}}\Vert f-f^*\Vert_{L^1(\nu)}$, so as to bound  the excess risk of the ERM $\hat{f}_\phi$ defined in (\ref{erm_nn}).

	Next, the approximation error bound depends on $f^*$ through its modulus of continuity and is related to	
	the function class $\mathcal{F}_\phi=\mathcal{F}_{\mathcal{D},\mathcal{W},\mathcal{U},\mathcal{S},\mathcal{B}}
	$ via its parameters. Recall that the modulus of continuity $\omega_f$ of a function $f:[0,1]^d\to\mathbb{R}$ is defined as
	\begin{equation}
		\label{mod}
		\omega_f(r) :=\sup\{\vert f(x)-f(y)\vert:x,y\in[0,1]^d,\Vert x-y\Vert_2\leq r\}, {\rm for\ any\ } r\geq0.
	\end{equation}
	A uniformly continuous function $f$ has the property that $\lim_{r\to0}\omega_f(r)=\omega_f(0)=0$.   Based on the modulus of continuity, one can define different equicontinuous families of functions. For instance, the modulus $\omega_f(r)= \theta r$ describes the $\theta$-Lipschitz continuity; the modulus $\omega_f(r)=\theta r^\alpha$ with $\theta,\alpha>0$ defines the H\"older continuity.
	
	In this paper, we only assume that $f^*$ is uniformly continuous (in Assumption \ref{A3} given below), which is a %very loose
mild
assumption on the continuity of the unknown target function $f^*$, as
	the existing works  generally assume stronger smoothness assumptions on $f^*$.

	\begin{assumption}
	\label{A3}
	The target function $f^*:\mathcal{X}\to\mathbb{R}$ is uniformly continuous.
	\end{assumption}
	Assumption \ref{A3} implies that  $\lim_{r\to0}\omega_{f^*}(r)=\omega_{f^*}(0)=0$.  It provides a necessary condition for us to derive an explicit upper bounds of the approximation error in terms of its modulus of continuity $\omega_{f^*}$, and give sufficient conditions for the consistency of some estimators. 	

	The DNN approximation theories are actively studied in recent years. We provide a brief review on the most relevant literature in Appendix A for reader’s convenience, see also \cite{jiao2021deep} for some additional references. We will use the recent approximation results obtained by \cite{shen2019deep}, who established an explicit error bounds for the approximation using deep neural networks. Their results clearly describe how the error bounds depend on the network parameters and are non-asymptotic. Let $\mathbb{N}^+$ be be the set of positive integers. For any $N,M\in\mathbb{N}^+$ and $p\in[1,\infty]$, \cite{shen2019deep} showed that
	\begin{equation}
		\label{shen2019}
		\inf_{f\in\mathcal{F}_\phi}\Vert f-f^*\Vert_{L^p([0,1]^d)}\leq 19\sqrt{d}\, \omega_{f^*}(N^{-2/d}M^{-2/d})
	\end{equation}
	for a  H\"older continuous function $f^*$ with the modulus of continuity
	$\omega_{f^*}(\cdot)$ defined in (\ref{mod})  and the class of functions $\mathcal{F}_\phi$ consisting of
	ReLU activated feedforward networks with
	$$
	\text{width } \mathcal{W}=C_1\max\{d\lfloor N^{1/d} \rfloor, N+1\} \
	\text{ and depth } \ \mathcal{D}=12M+C_2,
	$$
	where $C_1=12,C_2=14$ if $p\in[1,\infty)$,  and $C_1=3^{d+3},C_2=14+2d$ if $p=\infty$. Such an approximation rate in terms of the width and depth is more
	%% generic
	informative and more useful than those characterized by just the size $\mathcal{S}$, as the upper bound represented by width and depth can imply the one in terms of the size $\mathcal{S}$. Moreover, the approximation error bound is explicit in the sense that it does not involve unclearly defined constant, in contrast to other existing approximation error bounds that involve an unknown prefactor or require %sufficiently large
	the network width $\mathcal{W}$ and depth $\mathcal{D}$ greater than some unknown constants.
	
	  Note that the approximation error bound  in Lemma \ref{lemma2} only requires the Lipschitz property of the loss function; if the curvature information of the risk around $f^*$ can be used, the approximation error bound can be further improved.
	
	\begin{assumption}[Local quadratic bound of the excess risk]
		\label{quadratic}
		There exist some constants $\Lambda_{L,f^*}>0$ and $\varepsilon_{L,f^*}>0$ which may depend on the loss function $L(\cdot)$, the risk minimizer $f^*$ and the distribution of $(X,Y)$ such that
		$$\mathcal{R}(f)-\mathcal{R}(f^*)\leq \Lambda_{L,f^*}\Vert f-f_0\Vert^2_{L^2(\nu)},$$
		for any $f$ satisfying $\Vert f-f_0\Vert_{L^\infty(\nu)}\leq \varepsilon_{L,f^*}$, where $\nu$ is the probability measure of $X$.
	\end{assumption}

	\begin{lemma}\label{quadapprox}
		Assume that Assumption \ref{A1} and \ref{quadratic} hold. Let $f^*$ be the target function defined in (\ref{target}) and $\mathcal{R}(f^*)$ be its risk, then we have
		$$\inf_{f\in\mathcal{F}_\phi}\{R(f) -R(f^*)\}\leq \lambda_{L,f^*}\inf_{f\in\mathcal{F}_{\phi}}\Vert f-f^*\Vert^2_{L^2(\nu)},$$
		where $\lambda_{L,f^*}\geq\max\{\Lambda_{L,f^*},\lambda_L/\varepsilon_{L,f^*}\}>0$ is a constant, $\nu$ is the probability measure of $X$ and $\mathcal{F}_\phi=\mathcal{F}_{\mathcal{D},\mathcal{W},\mathcal{U},\mathcal{S},\mathcal{B}}$ denotes the class of feedforward neural networks with parameters $\mathcal{D},\mathcal{W},\mathcal{U},\mathcal{S}$ and $\mathcal{B}$.
	\end{lemma}

	\section{Non-asymptotic error bounds} %Rate of convergence}
	\label{sec4}
	
	Lemma \ref{lemma1} gives the theoretical basis for establishing the consistency and non-asymptotic error bounds. To guarantee consistency, the two terms on the right hand side of (\ref{oracle}) should go to zero as $n\to \infty$. Regarding the non-asymptotic error bound,  the exact rate of convergence will be determined by a trade-off between the stochastic error and the approximation error.
	%depend on the slower rate of the two vanishing items.
	
	\subsection{Consistency and non-asymptotic error bounds}
	In the following, we present the consistency and  non-asymptotic error bound
	of the general robust nonparametric regression estimator using neural networks.
	
	\begin{theorem}[Consistency]\label{thm1}
	For model (\ref{model}), suppose that Assumptions \ref{A1}-\ref{A2} hold. Assume that the target function $f^*$ is continuous function on $[0,1]^d$,  $ \Vert f^*\Vert_\infty\leq \mathcal{B}$ for some $\mathcal{B}\geq1$, and the function class of feedforward neural networks $\mathcal{F}_\phi=\mathcal{F}_{\mathcal{D},\mathcal{W},\mathcal{U},\mathcal{S},\mathcal{B}}$ with continuous piecewise-linear activation function of finite pieces satisfies
		\begin{equation*}
			\mathcal{S}\to\infty \quad {\rm and} \quad \mathcal{B} \frac{\log(n)}{n^{1-1/p}} \mathcal{S}\mathcal{D}\log(\mathcal{S})\to 0,
		\end{equation*}
		as $n\to\infty$.
		Then,  the excess risk of the empirical risk minimizer $\hat{f}_\phi$  defined in (\ref{erm_nn}) is consistent in the sense that
		$$\mathbb{E} \{\mathcal{R}(\hat{f}_\phi)-\mathcal{R}(f^*)\}\to 0 \ \  {\rm as\ } n\to\infty.$$
	\end{theorem}
	Theorem {\ref{thm1}}  is a direct result from Lemma \ref{lemma1} and Theorem 1 in \cite{yarotsky2018optimal}.
	Those conditions in Theorem \ref{thm1} are sufficient  conditions to ensure the consistency of the deep neural regression, and they are relatively
% loose
mild
in terms of the assumptions on the target function $f^*$ and the distribution of $Y$.

	\begin{theorem}[Non-asymptotic excess risk bound] \label{thm2}
		%Let $f^*$ be the target function defined in (\ref{target}),
		Suppose that Assumptions \ref{A1}-\ref{A3} hold, $\nu$ is absolutely continuous with respect to the Lebesgue measure, and $ \Vert f^*\Vert_\infty\leq\mathcal{B}$ for some $\mathcal{B}\geq1$. Then, for any $N, M\in\mathbb{N}^+$, the function class of ReLU multi-layer perceptrons $\mathcal{F}_{\phi}=\mathcal{F}_{\mathcal{D},\mathcal{W},\mathcal{U},\mathcal{S},\mathcal{B}}$ with depth $\mathcal{D}=12M+14$ and width $\mathcal{W}=\max\{4d\lfloor N^{1/d}\rfloor+3d,12N+8\}$,
		%for $n$ large enough,
		for $2n \ge \text{Pdim}(\mathcal{F}_\phi)$,
		the expected excess risk of the ERM $\hat{f}_\phi$ defined in (\ref{erm_nn})
		satisfies
		$$\mathbb{E} \{\mathcal{R}(\hat{f}_\phi)-\mathcal{R}(f^*)\}\leq C\lambda_L\mathcal{B}\frac{\mathcal{S}\mathcal{D}\log(\mathcal{S})\log(n)}{n^{1-1/p}}
		+18\lambda_L\sqrt{d}\omega_{f^*}(N^{-2/d}M^{-2/d}),$$
		where $C>0$ is a constant which does not depend on $n,d,\mathcal{B},\mathcal{S},\mathcal{D},\mathcal{W},\lambda_L,N$ or $M$. Here $\lambda_L$ is the Lipschitz constant of the robust loss function $L$ under Assumption \ref{A1},  and $\omega_{f^*}$ is the modulus of continuity of the uniformly continuous target $f^*$ under Assumption \ref{A3}.
		
		Moreover, aside from Assumptions \ref{A1}-\ref{A3},  if Assumption \ref{quadratic} also holds, then we have
			$$\mathbb{E} \{\mathcal{R}(\hat{f}_\phi)-\mathcal{R}(f^*)\}\leq C\lambda_L\mathcal{B}\frac{\mathcal{S}\mathcal{D}\log(\mathcal{S})\log(n)}{n^{1-1/p}}
		+384\lambda_{L,f^*}d\big\{\omega_{f^*}(N^{-2/d}M^{-2/d})\big\}^2,$$
		where $\lambda_{L,f^*}=\max\{\Lambda_{L,f^*},\lambda_L/\varepsilon_{L,f^*}\}>0$ is a constant defined in Lemma \ref{quadapprox}.
	\end{theorem}

	 Theorem \ref{thm2}  implies that the upper bound of the excess risk %$\mathbb{E} \Vert \hat{f}_\phi-f_0\Vert^2_{L^2(\nu)}$
	is a sum of two terms: the upper bound of the stochastic error $C\lambda_L\mathcal{B}\mathcal{S}\mathcal{D}\log(\mathcal{S})\log(n)/{n^{1-1/p}}$ and that of  the approximation error $18\lambda_L\sqrt{d}\omega_{f^*}(N^{-2/d}M^{-2/d})$.
	A few important comments are in order.  First, our error bound is non-asymptotic and explicit  in the sense that  no unclearly defined %unknown
	constant is involved. Importantly, the prefactor $18\lambda_L\sqrt{d}$ in the upper bound of approximation error depends on $d$ through $\sqrt{d}$ (the root of the dimension $d$), which is quite different from the exponential dependence in many existing results. Second, the approximation rate $\omega_{f^*}(N^{-2/d}M^{-2/d})$
	is in terms of the network width $\mathcal{W}=\max\{4d\lfloor N^{1/d}\rfloor+3d,12N+8\}$ and depth $\mathcal{D}=12M+14$, rather than
	only the network size $\mathcal{S}$.
	This provides insights into the relative merits of different network architecture and provides certain
	qualitative guidance on the network design. Third, our result on the approximation rate accommodates all kinds of uniformly continuous target functions $f^*$ and the result is  expressed in terms of the  modulus of continuity of $f^*$.
	
	To achieve the optimal rate,
	%guarantee a targeted accuracy,
	we need to balance the trade-off between the stochastic error and the approximation error. On one hand, the upper bound on the stochastic error $C\lambda_L\mathcal{B}\mathcal{S}\mathcal{D}\log(\mathcal{S})\log(n)/{n^{1-1/p}}$ increases in the complexity and richness of the function class $\mathcal{F}_{\mathcal{D},\mathcal{W},\mathcal{U},\mathcal{S},\mathcal{B}}$;  larger $\mathcal{D}$, $\mathcal{S}$ and $\mathcal{B}$ lead to a larger upper bound on the stochastic error. On the other hand, the upper bound for the approximation error $18\lambda_L\sqrt{d}\omega_{f^*}(N^{-2/d}M^{-2/d})$  decreases as the size of  $\mathcal{F}_{\mathcal{D},\mathcal{W},\mathcal{U},\mathcal{S},\mathcal{B}}$ increases,
	larger $\mathcal{D}$ and $\mathcal{W}$ result in smaller upper bound on the approximation error.

	 Apart form the excess risk bounds, if the loss function is self-calibrated, the distance between the estimate and the target can also be upper bounded via the self-calibration inequality.
		\begin{assumption}[Self-calibration]	\label{calib}	
		There exists a constant $C_{L,f^*}>0$ which may depend on the loss function $L(\cdot)$, the risk minimizer $f^*$ and the distribution of $(X,Y)$ such that
			\begin{equation*}
				\Delta^2(f,f^*) \leq C_{L,f^*} \big\{\mathcal{R}(f)-\mathcal{R}(f^*)\big\}
			\end{equation*}
		 for any $f:\mathcal{X}\to\mathbb{R}$, where $\Delta^2(\cdot,\cdot)$ is some distance defined in the functional space $\{f :\mathcal{X}\to\mathbb{R}\}$ such that $\Delta^2(f_1,f_2)=0$ if $\Vert f_1-f_2\Vert_{L^\infty(\nu)}=0$.
		\end{assumption}
		
		\begin{remark}
			For least squares estimation, Assumption \ref{calib} naturally holds with  $C_{L,f^*}=1$ and $\Delta^2(f,f^*)=\Vert f-f^*\Vert^2_{L^2(\nu)}=\mathcal{R}(f)-\mathcal{R}(f^*)$. For least absolute deviation or quantile check loss, under mild conditions,  the self-calibration assumption is satisfied with $\Delta^2(f,f^*)=\mathbb{E}\big[\min\{\vert f(X)-f^*(X)\vert,\vert f(X)-f^*(X)\vert^2\}\big]$ and some $C_{L,f^*}>0$. More details can be found in \cite{christmann2007svms,steinwart2011estimating,lv2018oracle} and \cite{padilla2020quantile}.
		\end{remark}
		With the self-calibration condition satisfied, the distance between the estimate and the target is upper bounded as a direct consequence of Theorem \ref{thm2}.
	
		\begin{corollary}\label{coro1}
			Suppose that Assumptions \ref{A1}-\ref{A3} and \ref{calib} hold, $\nu$ is absolutely continuous with respect to the Lebesgue measure, and $ \Vert f^*\Vert_\infty\leq\mathcal{B}$ for some $\mathcal{B}\geq1$. Then, for any $N, M\in\mathbb{N}^+$, the function class of ReLU multi-layer perceptrons $\mathcal{F}_{\phi}=\mathcal{F}_{\mathcal{D},\mathcal{W},\mathcal{U},\mathcal{S},\mathcal{B}}$ with depth $\mathcal{D}=12M+14$ and width $\mathcal{W}=\max\{4d\lfloor N^{1/d}\rfloor+3d,12N+8\}$,
			for $2n \ge \text{Pdim}(\mathcal{F}_\phi)$,
			the expected distance between $f^*$ and the ERM $\hat{f}_\phi$ defined in (\ref{erm_nn}) satisfies
			$$\mathbb{E} \big\{\Delta^2(\hat{f}_\phi,f^*)\big\}\leq C_{L,f^*}\Big\{ C\lambda_L\mathcal{B}\frac{\mathcal{S}\mathcal{D}\log(\mathcal{S})\log(n)}{n^{1-1/p}}
			+18\lambda_L\sqrt{d}\omega_{f^*}(N^{-2/d}M^{-2/d})\Big\},$$
			where $C_{L,f^*}$ is defined in Assumption \ref{calib} and $C>0$ is a constant which does not depend on $n,d,\mathcal{B},\mathcal{S},\mathcal{D},\mathcal{W},\lambda_L,N$ or $M$.
			
			Additionally, other than Assumptions \ref{A1}-\ref{A3} and  \ref{calib}, 	if Assumption \ref{quadratic} also holds, we have
			$$\mathbb{E} \big\{\Delta^2(\hat{f}_\phi,f^*)\big\}\leq C_{L,f^*}\Big[ C\lambda_L\mathcal{B}\frac{\mathcal{S}\mathcal{D}\log(\mathcal{S})\log(n)}{n^{1-1/p}}
			+384\lambda_{L,f^*} d\big\{\omega_{f^*}(N^{-2/d}M^{-2/d})\big\}^2\Big],$$
			where $\lambda_{L,f^*}=\max\{\Lambda_{L,f^*},\lambda_L/\varepsilon_{L,f^*}\}>0$ is a constant defined in Lemma \ref{quadapprox}.
		\end{corollary}

	\subsection{On the curse of dimensionality}
	\label{sec6}
	In many contemporary  statistical and machine learning applications,
	%tasks with deep neural network,
	the dimensionality  $d$ of the input $X$ can be large,  which leads to an extremely slow rate of convergence, even with huge sample size.  This problem is well known as the curse of dimensionality.
%{\color{red}	Note that the approximation error $\inf_{f\in\mathcal{F}_\phi}\Vert f-f^*\Vert^2_{L^2(\nu)}$ in Lemma \ref{lemma1} and \ref{lemma2}} is defined with respect to the probability measure $\nu$.
	%In what follows,  we focus on the support of $X$ to derive a theoretical bound.
	A promising  and effective way to mitigate the curse of dimensionality is to impose additional conditions on $X$ and the target function $f^*$.
	Though the domain of  $f^*$ is high dimensional, when the support of $X$ is concentrated on some neighborhood of a low-dimensional manifold, the upper bound of the approximation error  can be substantially improved in terms of the exponent of the convergence rate  \citep{shen2019deep}.
	
	\begin{assumption}
		\label{A4}
		The input $X$ is supported on $\mathcal{M}_\rho$, a $\rho$-neighborhood of $\mathcal{M}\subset[0,1]^d$, where $\mathcal{M}$ is a compact $d_\mathcal{M}$-dimensional Riemannian submanifold \citep{lee2006riemannian} and
		$$\mathcal{M}_\rho=\{x\in[0,1]^d: \inf\{\Vert x-y\Vert_2: y\in\mathcal{M}\}\leq \rho\}$$
		for $\rho\in(0,1)$.
	\end{assumption}
	
	%For real-world data, they are hardly observed to locate on an exact manifold, instead they could be
	Assumption \ref{A4}  is more realistic than
	the exact manifold support assumption assumed in \citet{schmidt2019deep},\citet{nakada2019adaptive} and \citet{chen2019nonparametric}.	
	In practical application, one can rarely observe data that are located on an exact manifold. It is  more reasonable  to assume that the data are concentrated on a neighborhood of a low-dimensional manifold.
	%viewed as contaminated data  originally supported on a low-dimensional manifold $\mathcal{M}$.

		\begin{theorem}[Non-asymptotic error bound] \label{thm3}
			Under model (\ref{model}), suppose that Assumptions \ref{A1}-\ref{A3} and \ref{A4} hold, the probability measure $\nu$ of $X$ is absolutely continuous with respect to the Lebesgue measure and $ \Vert f_0\Vert_\infty\leq\mathcal{B}$ for some $\mathcal{B}\geq1$. Then for any $N, M\in\mathbb{N}^+$, the function class of ReLU multi-layer perceptrons $\mathcal{F}_{\phi}=\mathcal{F}_{\mathcal{D},\mathcal{W},\mathcal{U},\mathcal{S},\mathcal{B}}$ with depth $\mathcal{D}=12M+14$ and width $\mathcal{W}=\max\{4d_\delta\lfloor N^{1/d_\delta}\rfloor+3d_\delta,12N+8\}$,
			the expected excess risk of the ERM $\hat{f}_\phi$ defined in (\ref{erm_nn}) satisfies
			$$\mathbb{E} \{\mathcal{R}(\hat{f}_\phi)-\mathcal{R}(f^*)\}\leq C_1\lambda_L\mathcal{B}\frac{\mathcal{S}\mathcal{D}\log(\mathcal{S})\log(n)}{n^{1-1/p}}
			+\lambda_L (2+18\sqrt{d_{\delta}}\,) \, \omega_{f^*}\big\{(C_2+1)(NM)^{-2/d_\delta}\big\},$$
			for $2n \ge \text{Pdim}(\mathcal{F}_\phi)$ and $\rho\leq C_2(NM)^{-2/d_{\delta}}(1-\delta)/\{2(\sqrt{d/d_\delta}+1-\delta)\}$, where $d_\delta=O(d_\mathcal{M}{\log(d/\delta)}/{\delta^2})$ is an integer  such that $d_\mathcal{M}$$\leq d_\delta<d$ for any $\delta\in(0,1)$ and $C_1,C_2>0$ are constants which do not depend on $n,d,\mathcal{B},\mathcal{S},\mathcal{D},\mathcal{W},\lambda_L,N$ or $M$. Here $\lambda_L$ is the Lipschitz constant of the general robust loss function $L$ under  Assumption \ref{A1},  and $\omega_{f^*}$ is the modulus of continuity of the uniformly continuous target $f^*$ under Assumption \ref{A3}.
			
			Additionally, aside from Assumptions \ref{A1}-\ref{A3} and \ref{A4},
			if Assumption \ref{quadratic} also holds, then
			\begin{align*}
				\mathbb{E} \{\mathcal{R}(\hat{f}_\phi)-\mathcal{R}(f^*)\}\leq& C_1\lambda_L\mathcal{B}\frac{\mathcal{S}\mathcal{D}\log(\mathcal{S})\log(n)}{n^{1-1/p}}\\
				&\qquad+\lambda_{L,f^*}\Big[(2+18\sqrt{d_{\delta}}\,)\, \omega_{f^*}\big\{(C_2+1)(NM)^{-2/d_\delta}\big\}\Big]^2,
			\end{align*}
			where $\lambda_{L,f^*}=\max\{\Lambda_{L,f^*},\lambda_L/\varepsilon_{L,f^*}\}>0$ is a constant defined in Lemma \ref{quadapprox}.
			
			Furthermore, if Assumption \ref{calib} also holds, then
			\begin{align*}
				\mathbb{E} \{\Delta^2(\hat{f}_\phi,f^*)\}\leq C_{L,f^*}\Bigg( C_1\lambda_L\mathcal{B}&\frac{\mathcal{S}\mathcal{D}\log(\mathcal{S})\log(n)}{n^{1-1/p}}\\
				&+\lambda_{L,f^*}\Big[\big(2+18\sqrt{d_{\delta}}\,\big)\, \omega_{f^*}\big\{(C_2+1)(NM)^{-2/d_\delta}\big\}\Big]^2\Bigg),
			\end{align*}
			where $C_{L,f^*}>0$ is a constant defined in Assumption \ref{calib}.
			\end{theorem}
	
			Theorem \ref{thm3} is established under the assumption that the distribution of $X$ is supported on an approximate Riemannian manifold with an intrinsic dimension lower than the  dimension $d$ of the ambient space $\mathbb{R}^d$.
			This is different from the hierarchical structure assumption on
			the target function  in \citet{bauer2019deep} and  \citet{schmidt2020nonparametric}.
			These  two assumptions are two different lines of assumptions on the underlying target function, and either one  can mitigate the curse of dimensionality.

\section{Error bounds for different network structures}
%Comparing networks with different structures}
\label{efficiency}Theorem \ref{thm2} provides a general
%n explicit
description of how the non-asymptotic error bounds depend on the network parameters.
In this section, we consider the error bounds for some specific network structures.
%which can be used to quantify
We also consider the relative efficiency of  networks with different shapes in terms of the size of the network needed to achieve the optimal error bound.
The calculations given below demonstrate the advantages of deep networks over shallow ones in the sense that deep networks can achieve the same error bound as the shallow networks with a fewer total number of parameters in the network.
%We will make this statement quantitatively clear in terms of the notion of relative %efficiency between
%		networks defined below.

 To give concrete examples on the convergence rate of the excess risk, we assume the target function $f^*$ to be H\"older continuous, i.e. there exist $\theta \ge 0$ and $\alpha>0$ such that
 \begin{equation}
 	\vert f^*(x)-f^*(y)\vert\leq \theta\Vert x-y\Vert_2^\alpha,\quad x,y\in\mathcal{X}=[0,1]^d.
 \end{equation}

%		\subsection{Relative efficiency of network structures} \label{relefficiency}
		Let $\cS_1$ and $\cS_2$ be the sizes of two neural networks
		$\cN_1$ and $\cN_2$ needed to achieve the same non-asymptotic error bound as given in Theorem \ref{thm2}. We define the relative efficiency between two networks $\cN_1$ and $\cN_2$ as
		\begin{equation}
			\label{re}
			\text{REN}(\cN_1, \cN_2) = \frac{\log \cS_2}{\log \cS_1}.
		\end{equation}
		Here we use the logarithm of the size because the size of the network for achieving the optimal error rate  usually has the form $\cS=n^s$ for some $s > 0$ up to a factor only involving the power of $\log n$, as will be seen below.  Let $r=\text{REN}(\cN_1, \cN_2)$.
In terms of sample sizes, this definition of relative efficiency implies that, if it takes a sample of size $n$ for network $\cN_1$ to achieve the optimal error rate, then it will take a sample of size $n^r$ for network $\cN_2$  to achieve the same error rate.
		
For any multilayer neural network in $\mathcal{F}_{\mathcal{D},\mathcal{W},\mathcal{U},\mathcal{S},\mathcal{B}}$, its parameters naturally satisfy
		\begin{equation}
			\label{size}
			\max\{\mathcal{W},\mathcal{D}\}\leq
			\mathcal{S}\leq
			\mathcal{W}(d+1)+(\mathcal{W}^2+\mathcal{W})(\mathcal{D}-1)+\mathcal{W}+1
			=O(\mathcal{W}^2\mathcal{D}).
		\end{equation}
	
		Based on the relationship in (\ref{size})  and Theorem \ref{thm2},  we next present Corollary \ref{c-all} for the networks with some specific structures without assuming
Assumption \ref{quadratic}.
% when Assumption \ref{quadratic} does not hold.
Similar corollaries with Assumption \ref{quadratic} satisfied can be established accordingly, and the relative efficiency of different networks can be compared.
%Network structure design satisfying Assumption \ref{quadratic}  will be illustrated in %section \ref{design}.
Denote
\begin{equation}
\label{npda}
n_*=n_{d,p,\alpha}=n^{\left(1-\frac{1}{p}\right)\frac{d}{d+\alpha}}.
\end{equation}
Recall if the response $Y$ is sub-exponential, we take $p=+\infty.$
For the function class of ReLU multi-layer perceptrons $\mathcal{F}_\phi=\mathcal{F}_{\mathcal{D},\mathcal{W},\mathcal{U},\mathcal{S},\mathcal{B}}$ with depth $\mathcal{D}$,  width $\mathcal{W}$ and size  $\mathcal{S} $ given by
\begin{enumerate}[(a)]
\item Deep with fixed width networks ($\mathcal{N}_{\text{DFW}}$):
% for any $N\in\mathbb{N}^+$,
\begin{align*}
	\mathcal{D}_1&=12\lfloor n_*^{1/2}(\log n)^{-1}\rfloor+14,
%			\mathcal{W}_1&= \max\{4d\lfloor N^{1/d}\rfloor+3d,12N+8\}, \\
\mathcal{W}_1=\max\{7d, 20\},
				\mathcal{S}_1= O(n_*^{1/2}(\log n)^{-1}),
\end{align*}	
\item 	Wide and fixed depth networks ($\mathcal{N}_{\text{WFD}}$):
%for any $M\in\mathbb{N}^+$,
		\begin{align*}
				\mathcal{D}_2&= 26,
%12M+14,
				\mathcal{W}_2= \max\{4d\lfloor n_*^{1/(2d)}\rfloor+3d,12\lfloor n_*^{1/2}\rfloor+8\},
				\mathcal{S}_2=O(n_*(\log n)^{-1}),
			\end{align*}
\item Deep and wide networks ($\mathcal{N}_{\text{DAW}}$):
\begin{align*}
\mathcal{D}_3& =O(n_*^{1/4}),
\mathcal{W}_3=O(n_*^{1/4}),
				\mathcal{S}_3=O(n_*^{3/4}(\log n)^{-2}).
\end{align*}

\end{enumerate}	

\begin{corollary}[Error bounds]
			\label{c-all}
			Let $f^*$ defined in (\ref{target}) be H\"older continuous of order $\alpha>0$ with H\"older constant $\theta\geq0$. Suppose that Assumptions  \ref{A1}-\ref{A3} hold, $\nu$ is absolutely continuous with respect to the Lebesgue measure, and $\Vert f^*\Vert_\infty\leq\mathcal{B}$ for some $\mathcal{B}\geq1$. Then,
%if the depth, width and size of the network satisfies
for any one of the three types of networks $\mathcal{N}_{\text{DFW}}, \mathcal{N}_{\text{WFD}}$ and $\mathcal{N}_{\text{DAW}}$ specified  above,
			the ERM $\hat{f}_\phi$ defined in (\ref{erm_nn}) satisfies
			$$\mathbb{E} \{\mathcal{R}(\hat{f}_\phi)-\mathcal{R}(f^*)\}\leq \lambda_L\left\{C\mathcal{B}
			+18\theta\sqrt{d}\right\}n^{-(1-\frac{1}{p})\frac{\alpha}{d+\alpha}},$$
			%	for large enough $n$,	
			for $2n \ge \text{Pdim}(\mathcal{F}_\phi)$, where $C>0$ is a constant which does not depend on $n,\lambda,\alpha$ or $\mathcal{B}$, and $C$ does not depend on $d$ and $N$ if $d\lfloor N^{1/d}\rfloor+3d/4\leq3N+2$, otherwise $C=O(d^2\lfloor N^{2/d}\rfloor)$.
			Additionally, if Assumption \ref{calib} holds, then
				$$\mathbb{E} \{\Delta^2(\hat{f}_\phi,f^*)\}\leq \lambda_LC_{L,f^*}\left\{C\mathcal{B}
				+18\theta\sqrt{d}\right\}n^{-(1-\frac{1}{p})\frac{\alpha}{d+\alpha}},$$
			where $C_{L,f^*}>0$ is a constant defined in Assumption \ref{calib}.
		\end{corollary}

Corollary \ref{c-all} is a direct consequence of Theorem \ref{thm2} and Corollary \ref{coro1}. By Corollaries \ref{c-all}, the sizes of the  three types of networks
%\textit{deep with fixed width} network
$\mathcal{N}_{\text{DFW}}, \mathcal{N}_{\text{WDF}}$ and
$\mathcal{N}_{\text{DAW}}$
%		and the size of the \textit{wide with fixed depth} network $S_{\text{WFD}}$
 to achieve the same error rate are
		\begin{equation*}
			\label{size1}
			\mathcal{S}_{\text{DFW}}=O(n_*^{1/2}(\log n)^{-1}), \
			\mathcal{S}_{\text{WFD}}=O(n_*(\log n)^{-1}), \ \text{ and } \
\mathcal{S}_{\text{DAW}}=O(n_*^{3/4}(\log n)^{-1}),
		\end{equation*}
		respectively. So we have  the relationship
\[
\mathcal{S}_{\text{DFW}}\approx \mathcal{S}_{\text{WFD}}^{1/2} \approx \mathcal{S}_{\text{DAW}}^{2/3}.
\]

The relative efficiency of these two networks as defined in (\ref{re}) are
		\[
		\text{REN}(\cN_{\text{DAW}}, \cN_{\text{DFW}})
		=\frac{1/2}{3/4} = \frac{2}{3},
		\text{REN}(\cN_{\text{DAW}}, \cN_{\text{WFD}})
		=\frac{1}{3/4} = \frac{4}{3}.
		\]
Thus the deep networks are  twice as efficient as  wide networks. If the sample size for a \textit{deep with fixed width} network to achieve the optimal error rate is $n$, then it is about $n^2$ for a \textit{wide with fixed depth} network and $n^{4/3}$ for a \textit{deep and wide} network.
An explanation is that deep neural networks are generally of greater approximation power than shallow networks. In \cite{telgarsky2016benefits}, it was shown that for any integer $k\geq1$ and dimension $d\geq1$, there exists a function computed by a ReLU neural network with $2k^3+8$ layers, $3k^2+12$ neurons and $4+d$ different parameters such that it can not be approximated by networks activated by piecewise polynomial functions with no more than $k$ layers and $O(2^k)$ neurons.
	
\begin{comment}
		The choice of the network parameters is not unique to achieve the optimal convergence rate.
		The convergence rate depends on the slower one of the stochastic error and the approximation error, i.e.,  $\max\{O(n^{-1+1/p}\mathcal{D}\mathcal{S}\log(\mathcal{S})),O(\sqrt{d}(\mathcal{W}\mathcal{D})^{-2\alpha/d})\}$.
		For deep and wide networks, we have multiple choices to achieve the optimal rate.
		For example, the following two different specifications of the network parameters achieve the
		same convergence rate.
		\begin{eqnarray*}
			\mathcal{D}&=& 12\lfloor n^{(1-1/p)d/2(d+\alpha)}(\log n)^{-2}\rfloor+14, \\
			\mathcal{W} &=& \max\{4d\lfloor \log(n)^{1/d}\rfloor+3d,12\lfloor \log(n)\rfloor+8\},\\
			\mathcal{S}&=&O(n^{(1-1/p)d/2(d+\alpha)});
		\end{eqnarray*}
		and
		\begin{eqnarray*}
			\mathcal{D}&=&12\lfloor \log(n)\rfloor+14, \\
			\mathcal{W}&=&\max\{4d\lfloor n^{(1-1/p)d/2(d+\alpha)}(\log n)^{-1} \rfloor^{1/d}+3d,12\lfloor n^{(1-1/p)d/2(d+\alpha)}\}(\log n)^{-1}\rfloor+8\}, \\
			\mathcal{S}&=&O(n^{(1-1/p)d/2(d+\alpha)}(\log n)^{-2}).
		\end{eqnarray*}
		%\todo{This set of network parameters was originally in Corollary \ref{c3}}
\end{comment}

		\subsection{Optimal design of rectangle networks}
		\label{design}
		We now discuss the optimal design of \textit{rectangle networks}, i.e.,  networks with equal width for each hidden layer. For such networks with a regular shape, we have an exact relationship between the size of the network, the depth and the width:
		\begin{equation}
			\label{size-equal}
			\mathcal{S}=\mathcal{W}(d+1)+(\mathcal{W}^2+\mathcal{W})(\mathcal{D}-1)+
			\mathcal{W}+1=O(\mathcal{W}^2\mathcal{D}).
		\end{equation}
		Based on this relationship
		%between $\mathcal{S}$ and $(\mathcal{D}, \mathcal{W})$
and Theorem \ref{thm2}, we can determine the depth and the width of the network to achieve the optimal rate with the minimal size. Below, we consider the best selection of
network depth and width under Assumption \ref{quadratic}.

Under  Assumption \ref{quadratic} holds, to achieve the optimal rate with respect to the sample size $n$ with a minimal network size, we can set
		\begin{eqnarray*}
			\mathcal{W}&=&\max(7d,20), \\
			\mathcal{D}&=& 12\lfloor n^{(1-1/p)d/{(2d+4\alpha)}}
			(\log n)^{-1}\rfloor+14,\\
			\mathcal{S}&=&
			\big \{\max(49d^2,400)+\max(7d,20)\big\}\times \big\{12\lfloor n^{(1-1/p)d/{(2d+4\alpha)}}(\log n)^{-1}\rfloor+14\big\}\\
			& & + (d+2)\max(7d,20)+ 1  \\
			&=&O(49d^2 n^{(1-1/p)d/{(2d+4\alpha)}}(\log n)^{-1}).
		\end{eqnarray*}
%		Similarly,
Furthermore, if Assumption \ref{A4} is satisfied, we can take $\mathcal{F}_\phi$  to be consisting of fixed-width networks with
		\begin{eqnarray*}
			\mathcal{W}&=&\max(7d_\delta,20), \\
			 \mathcal{D}&=& 12\lfloor n^{(1-1/p)d_\delta/{(2d_\delta+4\alpha)}}
			(\log n)^{-1}\rfloor (\log n)^{-1}+14,\\
			\mathcal{S}&=&\big\{\max(49d_\delta^2,400)+\max(7d_\delta,20)\big\}
			\times \big\{12\lfloor n^{(1-1/p)d_\delta/{(2d_\delta+4\alpha)}}\rfloor+14\big\}\\
			& & + (d_\delta+2)\max(7d_\delta,20)+ 1\\
			&=& O(49d_\delta^2 n^{(1-1/p)d_\delta/{(2d_\delta+4\alpha)}}(\log n)^{-1}).
		\end{eqnarray*}
Then the excess risk of $\hat{f}_\phi$ in Theorem \ref{thm3} satisfies
		\begin{equation*}
			\mathbb{E} \{\mathcal{R}(\hat{f}_\phi)-\mathcal{R}(f^*)\}\leq \{C_1\lambda_L\mathcal{B}+\lambda_{L,f^*}(2+18\sqrt{d_\delta})(C_2+1)^\alpha\theta\} n^{-\left(1-\frac{1}{p}\right)\frac{2\alpha}{d_\delta+2\alpha}},
		\end{equation*}
		where $C_2>0$ is a constant independent of $n,d_\delta,\mathcal{B},\mathcal{S},\mathcal{D},\lambda$ and $\alpha$, and $C_1>0$ is a constant that does not depend on $n,\mathcal{B},\mathcal{S},\mathcal{D}$ or $\lambda$,  and does not depend on $d_\delta$ if $7d_\delta\leq20$, otherwise $C_1=O(d_\delta^2)$. 	

	 Additionally, if Assumption \ref{calib} holds, then
		$$\mathbb{E} \{\Delta^2(\hat{f}_\phi,f^*)\}\leq C_{L,f^*}\{C_1\lambda_L\mathcal{B}+\lambda_{L,f^*}(2+18\sqrt{d_\delta})(C_2+1)^\alpha\theta\} n^{-\left(1-\frac{1}{p}\right)\frac{2\alpha}{d_\delta+2\alpha}},
%n^{-(1-1/p)2\alpha/(d_\delta+2\alpha)},
$$
where $C_{L,f^*}>0$ is a constant defined in Assumption \ref{calib}.

	%The relative efficiencies of networks with different shapes can be considered in a way completely similar to those in Subsections \ref{relefficiency} and \ref{design}.

\section{Experiments}\label{sim}	
We apply the deep nonparametric regression with different robust loss functions to the simulated data. We compare the results with the
nonparametric least squares estimates under different settings of error distributions.
All the deep network computations are implemented on \textit{Python} and mainly via \textit{PyTorch}.

\subsection{Network configuration}
In all settings, we implement the empirical risk minimization by ReLU activated fixed width multilayer perceptrons.
%For simplicity,
We consider %two classes
a class of functions (implemented by fixed width networks)
%in the estimation: (i)
consisting of
ReLU activated multilayer perceptrons with 5 hidden layers and network width being $(d,256,256,256,256,256,1)$, denoted by \textit{Nets-256}.
%\textit{Nets-1}; (ii) ReLU activated multilayer perceptrons with 5 hidden layers and network width being $(d,512,512,512,512,512,1)$, denoted by \textit{Nets-2}.
All weights and bias in each layer are initialized by uniformly samples on bounded intervals according to the default initialization mechanism in \textit{PyTorch}.

The optimization algorithm is \textit{Adam}  \citep{kingma2014adam} with default learning rate $0.01$ and default $\beta=(0.9,0.99)$ (coefficients used for computing running averages of gradient and its square). We set the batch size to be $n/4$ in all cases where $n$ is the size of the training data, and train the network for at least $400$ epochs until the training loss converges or becomes sufficient steady.

\subsection{Estimations and Evaluations}
We consider  $5$ different loss functions: (i) Least square loss function, denoted by LS; (ii) Least absolute deviation loss function, denoted by LAD; (iii) Huber loss function with parameter $\zeta=1.345$, denoted by Huber; (iv) Cauchy loss function with parameter $\kappa=1$, denoted by Cauchy; (v) Tukey's biweight loss function with parameter $t=4.685$.

For each configuration, we use each loss function above as the training loss $L_{train}$ to obtain the ERM $\hat{f}_{train}$ within \textit{Nets-256},
$$\hat{f}_{train}=\arg\min_{f\in {\rm \textit{Nets-256}}} \frac{1}{n}\sum_{i=1}^{n}L_{train}(f(X^{train}_i),Y^{train}_i),$$
where $(X^{train}_i,Y^{train}_i)_{i=1}^n$ is the training data.  For each training loss function, we calculate the  testing risks under five different choices of loss function $L_{test}$ on the testing data $(X^{test}_t,Y^{test}_t)_{t=1}^T$, that is
$$\mathcal{R}(\hat{f}_{train},L_{test})=\frac{1}{T}\sum_{t=1}^{T}L_{test}(\hat{f}_{train}(X^{test}_{t}),Y^{test}_t).$$
We  report the averaged testing risk and its standard deviation over $R=10$ replications under different scenarios.  For each data generation model,  the testing data $(X^{test}_t,Y^{test}_t)_{t=1}^T$, $T=100,000$,  are generated to have the same distribution as the training data $(X^{train}_i,Y^{train}_i)_{i=1}^n$.

\subsection{Data generation: univariate models}
In our simulation, we first consider four univariate target functions, including ``Blocks'', ``Bumps'', ``Heavisine'' and ``Dopller'' \citep{dj1994}. The formulae for these functions are
given below.
\begin{enumerate}
	\item Blocks: $$f_0(x)=\sum h_iI(x>t_i),$$
	where
	\begin{align*}
		(h_i)&=(4,-5,-2.5,4,-3,2.1,4.3,-1.1,-2.1,-4.2),\\
		(t_i)&=(0.1,0.15,0.23,0.28,0.40,0.44,0.65,0.76,0.78,0.81).
	\end{align*}
	\item Bumps: $$f_0(x)=\sum h_i (1+\vert x-t_i\vert/w_i)^{-4},$$
	where
	\begin{align*}
		(h_i)&=(4,5,2.5,4,3,2.1,4.3,1.1,2.1,4.2),\\
		(t_i)&=(0.1,0.15,0.23,0.28,0.40,0.44,0.65,0.76,0.78,0.81),\\
		(w_i)&=(0.005,0.005,0.006,0.01,0.01,0.03,0.01,0.01,0.005,0.008).
	\end{align*}
	
	\item Heavisine: $$f_0(x)=4\sin(4\pi x)-{\rm sgn}(x-0.3)-{\rm sgn}(0.72-x).$$
	\item Dopller: $$f_0(x)=\{x(1-x)\}^{1/2}\sin\{2.2\pi/(x+0.15)\}.$$
\end{enumerate}

\begin{figure}[H]
	\includegraphics[width=\textwidth]{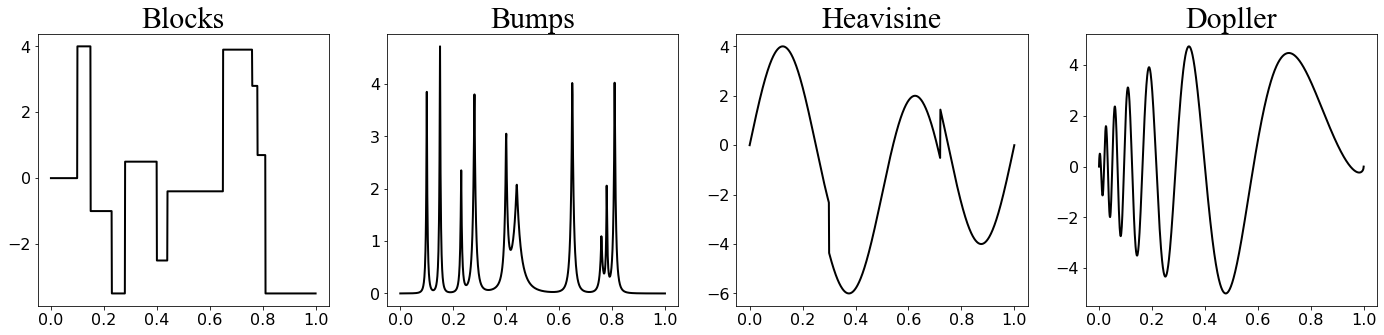}
	\caption{Four univariate target functions  including ``Block'', ``Bump'', ``Heavisine'' and ``Dopller''.}
	\label{fig2}
\end{figure}

\noindent The four target functions are  depicted in Fig \ref{fig2}.

We generate the training data $(X_i^{train},Y_i^{train})_{i=1}^n$,  $n=128$ or $512$ according to model (\ref{model}). The predictor $X^{train}$ is generated from Uniform$[0,1]$. For the error $\eta$, we have tried  several distributions: (i) Standard normal distribution, denoted by $\eta\sim\mathcal{N}(0,1)$; (ii) Student's $t$ distribution with $2$ degrees of freedom, denoted by $\eta\sim t(2)$; (iii) Standard Cauchy distribution(with location parameter 0 and scale parameter 1), denoted by $\eta\sim {\rm Cauchy(0,1)}$; (iv) Normal mixture distribution, denoted by ``Mixture'', where  $\eta\sim \xi \mathcal{N}(0,1)+(1-\xi) \mathcal{N}(0,10^4)$ with $\xi=0.8$
%$\xi \sim Bernoulli(0.8)$.
This is the contaminated normal model often used to depict the situation when $100 \xi \% $ of the observations are corrupt \citep{huber1964}.
Note that Assumption \ref{A2} is satisfied respectively in the way of (i) $\eta\sim\mathcal{N}(0,1)$, $\mathbb{E}\vert\eta\vert^p<\infty$ for $p\in[0,\infty)$; (ii) $\eta\sim t(2)$, $\mathbb{E}\vert\eta\vert^p<\infty$ for $p\in[0,2)$; (iii) $\eta\sim {\rm Cauchy(0,1)}$, $\mathbb{E}\vert\eta\vert^p<\infty$ for $p\in[0,1)$; (iv) $\eta$ follows normal mixture distribution, $\mathbb{E}\vert\eta\vert^p<\infty$ for $p\in[0,\infty)$. In addition, the error  distributions are symmetric; together with the symmetric property of the five loss functions,  the underlying true function $f_0$ in (\ref{model}) coincides with the target function $f^*$ in (\ref{target}) according to the definition of the risk minimizer.

The simulation results are shown in Table \ref{tab:blcok}-\ref{tab:dopller}.
For normal errors, except for ``Heavisine" model, the least squares estimation generally outperforms other methods in the sense that least squares estimator has smaller mean testing excess risk and standard deviation.
However, with heavier-tailed errors or contaminated error,  the performance of the least squares estimator becomes worse or even breaks down, while robust regression methods with LAD, Huber or Cauchy
loss give more accurate and stable results. Therefore, these robust methods
%estimations, which exhibits
possesses the desired robustness property with respect to heavy-tailed or contaminated error terms  in nonparametric estimations. The fitted curves are presented in Figure \ref{fig:fitted} and the training processes are displayed in Figures \ref{fig:blcoks}-\ref{fig:doppler} in Appendix.

\subsection{Data generation: multivariate models}
To simulate a
%nonparametric
multivariate
target function $f^*:[0,1]^d\to\mathbb{R}$,
% in real practice,
we refer to the Kolmogorov-Arnold (KA) representation theorem \citep{hecht1987kolmogorov}, which says that every continuous function can be represented by a specific network with two hidden layers. We formulate the original version of the KA representation theorem according to  \cite{schmidt2021kolmogorov}:
for any continuous function $f^*:[0,1]^d\to\mathbb{R}$, there exist univariate continuous functions $g_k,\psi_{m,k}$ such that
$$f^*(x_1,\ldots,x_d)=\sum_{k=0}^{2d}g_k\big(\sum_{m=1}^d\psi_{m,k}(x_m)\big).$$
In light of this, these $(2d + 1)(d + 1)$ univariate functions suffices for an exact representation of a $d$-variate function. We consider
$f^*$  according to the KA representation by choosing $g_k,\psi_{m,k},k=0,\ldots,2d,\ m=1,\ldots,d$ from a pool of univariate functions $\mathcal{H}=\{h_i(x),i=1,\ldots,7\}$, where $h_1(x)=-2.2x+0.3, h_2(x)=0.7x^3-0.2x^2+0.3x-0.3, h_3(x)=0.3\times{\rm sign}(x)\sqrt{\vert x\vert}, h_4(x)=0.8\times\log(\vert x\vert+0.01), h_5(x)=\exp(\min\{0.2x-0.1,4\}), h_6(x)=\sin(6.28x), h_7(x)={2}/({\vert x\vert+0.1}).$ For a given dimensionality $d$, we first generate $(2d + 1)(d + 1)$ random integers uniformly distributed on $\{1,2,3,4,5,6,7\}$ as indices of the functions  $g_k,\psi_{m,k}$; then we specify the choices of $g_k,\psi_{m,k}$ based on these randomly generated indices. To make our simulation reproducible, the seed is set as $2021$ in the random number generation, so the choices of $g_k,\psi_{m,k}$ can be reproducible.

We generate $X$ from the uniform distribution on $[0,1]^d$, i.e. $X\sim {\rm Unif} [0,1]^d$. We set the dimensionality $d$ of the input $X$ to be $4$
% or $12$
and generate  the training data $(X^{train}_i,Y^{train}_i)_{i=1}^n$ with sample size $n=128$ and $512$
% or $n=512,1024$
for each of $R=10$ replication.
The results are summarized in Tables \ref{tab:2}.
%-\ref{tab:3}.
 It can be seen that except for the case of normally distributed error, estimators based on robust loss functions always outperform least squares estimator in terms of the mean excess risk and standard deviation.
%When the dimensionality gets higher as shown in Table \ref{tab:3},  least squares %estimation is less competitive even with normal errors, but all estimators do not  perform %well  due to the high-dimensionality and the relatively small sample size.

%\begin{comment}
\begin{table}[H]
\setlength{\tabcolsep}{3pt} % Default value: 6pt
\renewcommand{\arraystretch}{1.1} % Default value: 1
\caption{The testing excess risks and standard deviation (in parentheses) across different loss functions of the ERMs trained with respect to different loss functions using \textit{Nets-256}. Data is generated from "Blocks" model with training sample size $n=128$ or $512$ and the number of replications $R=10$.}
\centering
\label{tab:blcok}
\resizebox{\textwidth}{!}{%
\begin{tabular}{@{}cc|ccccc@{}}
\toprule
\multicolumn{2}{c|}{Blocks,   $n=128$}     & \multicolumn{5}{c}{Testing Loss}                                            \\ \midrule
Error                     & Training Loss & LS              & LAD          & Huber        & Cauchy       & Tukey        \\ \midrule
\multirow{5}{*}{$N(0,1)$}
& LS & \textbf{13.534(4.078)}            & \textbf{1.408(0.348)}     & \textbf{1.363(0.340)}     & \textbf{0.766(0.180)} & 0.784(0.182) \\
& LAD           & 18.192(6.686)   & 1.708(0.285) & 1.655(0.281) & 0.890(0.128) & 0.891(0.153) \\
& Huber         & 18.178(6.367)   & 1.732(0.449) & 1.681(0.445) & 0.890(0.196) & 0.871(0.195) \\
& Cauchy        & 42.140(13.466)  & 2.571(0.655) & 2.531(0.653) & 0.967(0.182) & \textbf{0.755(0.112)}\\
& Tukey         & 159.659(15.485) & 8.089(0.683) & 8.020(0.677) & 2.516(0.281) & 1.645(0.225) \\ \midrule
\multirow{5}{*}{$t(2)$}   & LS            & \textbf{19.057(11.877)}  & 1.770(0.487) & 1.713(0.479) & 0.932(0.198) & 0.935(0.202) \\
& LAD           & 28.137(15.167)  & 2.160(0.695) & 2.110(0.694) & 0.974(0.180) & 0.873(0.095) \\
& Huber         & 21.393(7.411)   & \textbf{1.725(0.338)} & \textbf{1.679(0.336)} & 0.807(0.132) & 0.741(0.134) \\
& Cauchy        & 35.695(9.973)   & 2.093(0.450) & 2.057(0.449) & \textbf{0.782(0.126)} & \textbf{0.612(0.081)} \\
& Tukey         & 161.804(18.701) & 7.989(0.585) & 7.925(0.584) & 2.383(0.112) & 1.459(0.092) \\ \midrule
\multirow{5}{*}{$Cauchy(0,1)$} & LS & 1006121.625(2722768.500) & 187.531(486.105) & 187.457(486.096) & 2.914(3.326) & 1.481(0.516) \\
& LAD           & 26.255(16.659)  & 1.778(0.692) & 1.735(0.687) & 0.782(0.213) & 0.687(0.149) \\
& Huber         & \textbf{15.564(17.321)} & \textbf{1.717(0.475)} & \textbf{1.672(0.473)} & \textbf{0.774(0.142)} & 0.692(0.108) \\
& Cauchy        & 31.284(13.132)  & 2.057(0.399) & 2.018(0.400) & 0.785(0.101) & \textbf{0.611(0.066)} \\
& Tukey         & 156.523(18.607) & 7.206(0.624) & 7.145(0.622) & 2.147(0.135) & 1.314(0.108) \\ \midrule
\multirow{5}{*}{Mixture} & LS & 372.860(109.658) &	9.262(1.380) & 9.170(1.376)	 & 2.758(0.276)	& 2.138(0.134) \\
& LAD           & 42.559(16.447) &	2.264(0.458)	& 2.218(0.451)	& 0.943(0.163)	& 0.832(0.164) \\
& Huber         & \textbf{29.233(12.558)} &	1.853(0.529)	&  \textbf{1.806(0.526)}	&0.864(0.173)	&0.814(0.131)\\
& Cauchy        & 34.202(8.009)	&  \textbf{1.842(0.324)}	& 1.807(0.322)	&  \textbf{0.744(0.117)}&  \textbf{0.626(0.111)} \\
& Tukey         & 148.123(19.258) &	6.439(0.683) & 6.386(0.679)	& 1.974(0.227)	& 1.280(0.171) \\ \midrule

\multicolumn{2}{c|}{Blocks, $n=512$}   &   LS              & LAD          & Huber        & Cauchy   & Tukey                                    \\ \midrule
\multirow{5}{*}{$N(0,1)$} & LS            & \textbf{7.547(2.020)}   & \textbf{0.948(0.231)} & \textbf{0.908(0.224)} & \textbf{0.561(0.134)} & 0.596(0.148) \\
& LAD           & 10.319(3.478)   & 1.151(0.281) & 1.109(0.275) & 0.640(0.134) & 0.655(0.128) \\
& Huber         & 9.470(3.630)    & 1.051(0.349) & 1.013(0.344) & 0.581(0.176) & 0.590(0.175) \\
& Cauchy        & 29.639(4.070)   & 1.720(0.229) & 1.693(0.229) & 0.626(0.087) & \textbf{0.474(0.070)} \\
& Tukey         & 166.059(10.031) & 7.957(0.281) & 7.903(0.282) & 2.276(0.016) & 1.322(0.027) \\ \midrule
\multirow{5}{*}{$t(2)$}   & LS            & \textbf{8.003(2.626)}    & 1.047(0.290) & \textbf{1.000(0.277)} & 0.644(0.174) & 0.708(0.191) \\
& LAD           & 8.197(2.587)    & 1.045(0.279) & 1.003(0.271) & 0.624(0.161) & 0.669(0.174) \\
& Huber         & 8.649(3.074)    & \textbf{1.043(0.359)} & 1.002(0.346) & \textbf{0.611(0.221)} & 0.654(0.253) \\
& Cauchy        & 28.553(5.855)   & 1.702(0.358) & 1.674(0.356) & 0.634(0.127) & \textbf{0.497(0.090)} \\
& Tukey         & 171.790(5.422)  & 8.442(0.169) & 8.379(0.166) & 2.499(0.093) & 1.546(0.101) \\ \midrule
\multirow{5}{*}{$Cauchy(0,1)$} & LS & 214.717(529.612)         & 3.321(2.587)     & 3.256(2.580)     & 1.291(0.396) & 1.191(0.280) \\
& LAD           & \textbf{11.284(3.994)}   & 1.004(0.272) & 0.969(0.267) & 0.530(0.121) & 0.525(0.114) \\
& Huber         & 14.093(3.489)   & \textbf{0.981(0.180)} & \textbf{0.947(0.180)} & 0.506(0.068) & 0.491(0.067) \\
& Cauchy        & 26.687(5.734)   & 1.322(0.256) & 1.298(0.255) & \textbf{0.491(0.070)} & \textbf{0.381(0.047)} \\
& Tukey         & 183.687(23.299) & 7.052(0.809) & 7.007(0.804) & 1.897(0.259) & 1.058(0.164) \\ \midrule
\multirow{5}{*}{Mixture} & LS & 102.667(24.698)	& 5.587(0.766) &	5.498(0.763)&	2.313(0.227)	&2.019(0.160) \\
& LAD           & 13.765(4.882)	& 1.096(0.307)	& 1.060(0.301)	& 0.592(0.145)	& 0.606(0.153) \\
& Huber         & \textbf{10.115(2.084)}	&\textbf{0.852(0.182)} & \textbf{0.820(0.175)} & 0.486(0.112) & 0.510(0.125) \\
& Cauchy        & 24.557(4.408)	& 1.254(0.199)	& 1.233(0.198)	& \textbf{0.469(0.068)}	& \textbf{0.375(0.052)} \\
& Tukey         & 164.275(27.099)&	6.745(0.605)&	6.692(0.606)&	2.027(0.107)& 1.297(0.117)\\ \bottomrule
\end{tabular}%
}

\end{table}

\begin{table}[H]
\setlength{\tabcolsep}{5pt} % Default value: 6pt
\renewcommand{\arraystretch}{1} % Default value: 1
\caption{The testing excess risks and standard deviation (in parentheses)  under different loss functions of the ERMs trained with respect to different loss functions using \textit{Nets-256}. Data is generated from "Bumps" model with training sample size $n=128$ or $512$ and the number of replications $R=10$.}
\label{tab:bump}
\resizebox{\textwidth}{!}{%
\begin{tabular}{@{}cc|ccccc@{}}
\toprule
\multicolumn{2}{c|}{Bumps,   $n=128$}    & \multicolumn{5}{c}{Testing Loss}                                                                                      \\ \midrule
Error                   & Training Loss & LS                    & LAD                   & Huber                 & Cauchy                & Tukey                 \\ \midrule
\multirow{5}{*}{$N(0,1)$}      & LS    & \textbf{5.991(1.141)} & \textbf{0.798(0.114)} & \textbf{0.762(0.111)} & 0.493(0.063)          & 0.541(0.065)          \\
& LAD           & 6.859(1.065)          & 0.821(0.091)          & 0.785(0.088)          & \textbf{0.487(0.057)} & \textbf{0.520(0.070)} \\
& Huber         & 6.214(1.289)          & 0.799(0.116)          & 0.763(0.112)          & 0.487(0.065)          & 0.530(0.071)          \\
& Cauchy        & 8.011(0.631)          & 0.898(0.043)          & 0.862(0.039)          & 0.515(0.041)          & 0.538(0.057)          \\
& Tukey         & 8.480(0.519)          & 0.909(0.057)          & 0.874(0.055)          & 0.508(0.036)          & 0.525(0.042)          \\ \midrule
\multirow{5}{*}{$t(2)$} & LS            & 8.697(2.977)          & 1.018(0.223)          & 0.974(0.218)          & 0.593(0.103)          & 0.629(0.096)          \\
& LAD           & \textbf{7.491(0.719)} & 0.869(0.086)          & 0.829(0.081)          & 0.507(0.062)          & 0.535(0.077)          \\
& Huber         & 7.543(0.924)          & \textbf{0.866(0.080)} & \textbf{0.827(0.078)} & 0.498(0.051)          & 0.518(0.064)          \\
& Cauchy        & 8.743(0.582)          & 0.901(0.028)          & 0.865(0.028)          & \textbf{0.492(0.022)} & \textbf{0.497(0.033)} \\
& Tukey         & 9.477(0.951)          & 0.967(0.099)          & 0.930(0.096)          & 0.524(0.059)          & 0.526(0.068)          \\\midrule
\multirow{5}{*}{$Cauchy(0,1)$} & LS    & 774.661(1487.473)     & 3.991(3.920)          & 3.931(3.918)          & 1.031(0.228)          & 0.974(0.187)          \\
& LAD   & 11.500(2.969)         & \textbf{0.810(0.076)} & \textbf{0.773(0.073)} & \textbf{0.451(0.048)} & \textbf{0.458(0.053)} \\
& Huber         & 11.087(2.643)         & 0.852(0.142)          & 0.814(0.136)          & 0.478(0.081)          & 0.488(0.089)          \\
& Cauchy        & \textbf{9.486(2.771)} & 0.855(0.126)          & 0.820(0.122)          & 0.463(0.077)          & 0.463(0.092)          \\
& Tukey         & 11.173(3.318)         & 0.889(0.118)          & 0.853(0.116)          & 0.473(0.064)          & 0.463(0.065)          \\\midrule
\multirow{5}{*}{Mixture} & LS &  432.363(307.985)&	10.093(4.175)&	10.000(4.171)&	2.845(0.623)	&2.106(0.278)        \\
& LAD   &13.203(9.829)	& 1.064(0.378)	&1.027(0.374)	&0.569(0.142)	&0.584(0.123) \\
& Huber         & 21.442(31.043) &	1.232(0.627)&	1.194(0.625)&	0.601(0.135)&	0.598(0.096)        \\
& Cauchy        & 8.017(1.069)&	0.821(0.109)&	0.788(0.106)&	0.467(0.060)&	0.489(0.064)         \\
& Tukey         & \textbf{6.976(0.435)}&	\textbf{0.730(0.032)}	&\textbf{0.699(0.030)}	&\textbf{0.424(0.025)}&	\textbf{0.448(0.034)}        \\\midrule

\multicolumn{2}{c|}{Bumps, $n=512$}      & LS                    & LAD                   & Huber                 & Cauchy                & Tukey                 \\ \midrule
\multirow{5}{*}{$N(0,1)$}      & LS    & 4.248(1.836)          & 0.602(0.183)          & 0.573(0.176)          & 0.382(0.102)          & 0.425(0.110)          \\
& LAD           & 4.377(0.748)          & 0.570(0.061)          & 0.544(0.059)          & 0.350(0.030)          & 0.381(0.029)          \\
& Huber & \textbf{3.810(0.725)} & \textbf{0.530(0.070)} & \textbf{0.505(0.068)} & \textbf{0.332(0.040)} & \textbf{0.366(0.043)} \\
& Cauchy        & 5.605(0.738)          & 0.641(0.050)          & 0.616(0.050)          & 0.368(0.022)          & 0.385(0.023)          \\
& Tukey         & 6.499(0.872)          & 0.687(0.073)          & 0.663(0.071)          & 0.375(0.037)          & 0.378(0.037)          \\\midrule
\multirow{5}{*}{$t(2)$}        & LS    & \textbf{4.192(1.204)} & 0.606(0.134)          & 0.573(0.129)          & 0.386(0.076)          & 0.424(0.081)          \\
& LAD           & 4.704(0.773)          & 0.572(0.062)          & 0.546(0.060)          & 0.334(0.030)          & 0.349(0.031)          \\
& Huber & 4.231(1.020)          & \textbf{0.531(0.097)} & \textbf{0.504(0.094)} & \textbf{0.318(0.053)} & 0.336(0.056)          \\
& Cauchy        & 6.140(0.614)          & 0.650(0.051)          & 0.624(0.049)          & 0.355(0.030)          & 0.356(0.036)          \\
& Tukey         & 6.503(0.750)          & 0.653(0.049)          & 0.628(0.049)          & 0.345(0.020)          & \textbf{0.336(0.017)} \\ \midrule
\multirow{5}{*}{$Cauchy(0,1)$} & LS    & 281.707(421.853)      & 3.520(2.033)          & 3.457(2.027)          & 1.136(0.269)          & 1.053(0.221)          \\
& LAD   & \textbf{5.259(1.686)} & \textbf{0.493(0.054)} & 0.468(0.052)          & \textbf{0.287(0.032)} & 0.296(0.037)          \\
& Huber         & 5.266(1.538)          & 0.493(0.093)          & \textbf{0.467(0.091)} & 0.293(0.049)          & 0.306(0.050)          \\
& Cauchy        & 6.621(0.998)          & 0.575(0.051)          & 0.552(0.050)          & 0.301(0.026)          & \textbf{0.290(0.027)} \\
& Tukey         & 7.514(0.964)          & 0.642(0.041)          & 0.616(0.040)          & 0.328(0.017)          & 0.309(0.017)          \\ \midrule
\multirow{5}{*}{Mixture} & LS    & 142.978(89.599)	& 5.354(1.038)	& 5.270(1.037)	& 2.072(0.207)	& 1.818(0.141)        \\
& LAD   &8.067(3.051) &	0.599(0.135)&	0.574(0.133)&	0.342(0.054)&	0.360(0.048)    \\
& Huber         & \textbf{5.834(0.673)}	& \textbf{0.502(0.049)}	& \textbf{0.478(0.047)}	&0.311(0.033)	& 0.342(0.039)      \\
& Cauchy        & 7.254(0.561)	& 0.541(0.037)	& 0.519(0.036)&	 \textbf{0.310(0.022)}	&0.325(0.026) \\
& Tukey         & 8.339(1.042)	& 0.585(0.054)	& 0.564(0.053)	& 0.317(0.024)	& \textbf{0.321(0.024)}    \\ \bottomrule
\end{tabular}%
}
\end{table}

\begin{table}[H]
\setlength{\tabcolsep}{4pt} % Default value: 6pt
\renewcommand{\arraystretch}{1.1} % Default value: 1
\caption{The testing excess risks and standard deviation (in parentheses)  under different loss functions of the ERMs trained with respect to different loss functions using \textit{Nets-256}. Data is generated from "Heavisine" model with training sample size $n=128$ or $512$ and  replications $R=10$.}
\label{tab:heavisine}
\resizebox{\textwidth}{!}{%
\begin{tabular}{@{}cc|ccccc@{}}
\toprule
\multicolumn{2}{c|}{Heavisine,   $n=128$}       & \multicolumn{5}{c}{Testing Loss}                                                                                      \\ \midrule
Error                          & Training Loss & LS                    & LAD                   & Huber                 & Cauchy                & Tukey                 \\ \midrule
\multirow{5}{*}{$N(0,1)$}      & LS            & 2.623(2.612)          & 0.582(0.455)          & 0.542(0.436)          & 0.432(0.305)          & 0.518(0.365)          \\
& LAD           & 1.742(0.354)          & 0.428(0.070)          & 0.393(0.065)          & 0.330(0.054)          & 0.392(0.065)          \\
& Huber         & \textbf{1.730(0.609)} & 0.441(0.166)          & 0.405(0.153)          & 0.343(0.133)          & 0.414(0.170)          \\
& Cauchy        & 1.823(0.564)          & \textbf{0.394(0.096)} & \textbf{0.366(0.091)} & \textbf{0.291(0.070)} & 0.340(0.085)          \\
& Tukey         & 288.203(541.200)      & 9.012(9.163)          & 8.934(9.144)          & 2.522(1.385)          & 1.747(0.722)          \\\midrule
\multirow{5}{*}{$t(2)$}        & LS            & 3.899(3.483)          & 0.641(0.227)          & 0.597(0.220)          & 0.448(0.127)          & 0.521(0.147)          \\
& LAD           & 2.173(1.005)          & 0.442(0.152)          & 0.408(0.145)          & 0.319(0.099)          & 0.367(0.113)          \\
& Huber  & \textbf{1.839(0.341)} & \textbf{0.412(0.066)} & \textbf{0.378(0.061)} & \textbf{0.307(0.049)} & 0.358(0.059)          \\
& Cauchy        & 2.422(1.596)          & 0.442(0.192)          & 0.410(0.185)          & 0.308(0.111)          & \textbf{0.347(0.117)} \\
& Tukey         & 9102.458(26252.125)   & 38.803(84.576)        & 38.721(84.566)        & 3.504(2.490)          & 1.857(0.550)          \\\midrule
\multirow{5}{*}{$Cauchy(0,1)$} & LS            & 1090.375(3103.706)    & 4.247(5.722)          & 4.184(5.715)          & 1.156(0.562)          & 1.049(0.340)          \\
& LAD           & 3.975(5.956)         & \textbf{0.214(0.051)} & \textbf{0.193(0.047)} & \textbf{0.153(0.036)} & \textbf{0.170(0.040)} \\
& Huber         & 4.150(9.158)         & 0.332(0.102)          & 0.302(0.096)          & 0.244(0.074)          & 0.277(0.089)          \\
& Cauchy        & \textbf{0.650(4.657)} & 0.224(0.038)          & 0.203(0.035)          & 0.160(0.027)          & 0.177(0.031)          \\
& Tukey         & 452.600(1232.173)     & 8.023(19.907)         & 7.989(19.885)         & 1.139(2.016)          & 0.576(0.682)          \\\midrule
\multirow{5}{*}{Mixture} & LS        &  349.663(180.878) &	8.992(2.705)&	8.905(2.701)&	2.652(0.427)&	2.015(0.199)     \\
& LAD         &  3.748(1.782)&	0.609(0.239)&	0.571(0.230)&	0.436(0.153)&	0.514(0.179)  \\
& Huber         &18.776(25.589)	 & 0.888(0.527)	& 0.851(0.525)	& 0.461(0.129)	& 0.496(0.109)     \\
& Cauchy        & \textbf{3.719(3.786)}&	\textbf{0.493(0.283)}&	\textbf{0.464(0.277)}	&\textbf{0.332(0.132)}&	\textbf{0.377(0.123)}        \\
& Tukey         & 1659.819(2865.520) &	19.828(25.013)	& 19.752(25.000)	&3.055(1.834)	&1.749(0.528)      \\\midrule

\multicolumn{2}{c|}{Heavisine, $n=512$}         & LS                    & LAD                   & Huber                 & Cauchy                & Tukey                 \\\midrule
\multirow{5}{*}{$N(0,1)$}      & LS            & 1.097(0.194)          & 0.287(0.060)          & 0.263(0.054)          & 0.225(0.050)          & 0.269(0.060)          \\
& LAD    & \textbf{1.063(0.275)} & \textbf{0.274(0.079)} & \textbf{0.251(0.071)} & \textbf{0.214(0.064)} & \textbf{0.255(0.075)} \\
& Huber         & 1.492(0.729)          & 0.379(0.182)          & 0.346(0.166)          & 0.296(0.144)          & 0.352(0.172)          \\
& Cauchy        & 1.108(0.299)          & 0.280(0.089)          & 0.256(0.080)          & 0.217(0.073)          & 0.256(0.087)          \\
& Tukey         & 24153.281(72343.156)  & 47.157(135.209)       & 47.110(135.182)       & 1.862(3.286)          & 0.909(0.918)          \\\midrule
\multirow{5}{*}{$t(2)$}        & LS            & 4.782(9.469)          & 0.451(0.300)          & 0.419(0.294)          & 0.290(0.126)          & 0.325(0.136)          \\
& LAD           & 1.456(0.464)          & 0.320(0.107)          & 0.291(0.098)          & 0.237(0.082)          & 0.272(0.097)          \\
& Huber         & 1.434(0.458)          & 0.328(0.117)          & 0.298(0.106)          & 0.245(0.091)          & 0.283(0.106)          \\
& Cauchy & \textbf{1.356(0.271)} & \textbf{0.291(0.048)} & \textbf{0.264(0.044)} & \textbf{0.213(0.034)} & \textbf{0.239(0.038)} \\
& Tukey         & 4.210(7.350)          & 0.434(0.470)          & 0.411(0.461)          & 0.248(0.205)          & 0.254(0.193)          \\\midrule
\multirow{5}{*}{$Cauchy(0,1)$} & LS            & 1090.375(3103.706)    & 4.247(5.722)          & 4.184(5.715)          & 1.156(0.562)          & 1.049(0.340)          \\
& LAD           & 3.975(5.956)         & \textbf{0.214(0.051)} & \textbf{0.193(0.047)} & \textbf{0.153(0.036)} & \textbf{0.170(0.040)} \\
& Huber         & 4.150(9.158)         & 0.332(0.102)          & 0.302(0.096)          & 0.244(0.074)          & 0.277(0.089)          \\
& Cauchy        & \textbf{0.650(4.657)} & 0.224(0.038)          & 0.203(0.035)          & 0.160(0.027)          & 0.177(0.031)          \\
& Tukey         & 452.600(1232.173)     & 8.023(19.907)         & 7.989(19.885)         & 1.139(2.016)          & 0.576(0.682)          \\ \midrule
\multirow{5}{*}{Mixture} & LS        &  100.342(40.824)	& 5.224(0.974)&	5.137(0.972)&	2.182(0.229)&	1.928(0.155))     \\
& LAD         &  0.864(0.483)&	0.297(0.090)&	0.271(0.082)&	0.235(0.073)&	0.281(0.087) \\
& Huber         & 0.885(0.371)	& 0.267(0.072)	& 0.244(0.065)	& 0.207(0.059)	& 0.246(0.070)     \\
& Cauchy        & \textbf{0.860(0.291)}&	\textbf{0.235(0.038)}	&\textbf{0.214(0.034)}&	\textbf{0.182(0.032)}&	\textbf{0.213(0.037)}     \\
& Tukey         & 51982.227(155652.234)	& 70.234(200.608)&	70.193(200.588)	&1.984(3.050)	&0.878(0.689) \\\bottomrule

\end{tabular}%
}
\end{table}

\begin{table}[H]
\setlength{\tabcolsep}{4pt} % Default value: 6pt
\renewcommand{\arraystretch}{1.1} % Default value: 1
\caption{The testing excess risks and standard deviation (in parentheses)  under different loss functions of the ERMs trained with respect to different loss functions using \textit{Nets-256}. Data is generated from "Dopller" model with training sample size $n=128$ or $512$ and the number of replications $R=10$.}
\label{tab:dopller}
\resizebox{\textwidth}{!}{%
\begin{tabular}{@{}cc|ccccc@{}}
\toprule
\multicolumn{2}{c|}{Dopller,  $n=128$}         & \multicolumn{5}{c}{Testing Loss}                                                                   \\ \midrule
Error                          & Training Loss & LS                     & LAD            & Huber          & Cauchy                & Tukey                 \\ \midrule
\multirow{5}{*}{$N(0,1)$} & LS    & 14.256(6.635)           & \textbf{1.385(0.359)} & \textbf{1.338(0.358)} & \textbf{0.731(0.130)} & \textbf{0.724(0.118)} \\
& LAD           & 30.378(8.560)          & 2.192(0.583)   & 2.139(0.575)   & 0.983(0.242)          & 0.895(0.202)          \\
& Huber         & \textbf{18.495(6.667)} & 1.630(0.378)   & 1.580(0.373)   & 0.821(0.149)          & 0.795(0.133)          \\
& Cauchy        & 67.333(33.892)         & 3.578(1.406)   & 3.519(1.398)   & 1.328(0.411)          & 1.088(0.270)          \\
& Tukey         & 4923.002(13946.816)    & 30.952(60.211) & 30.857(60.199) & 3.884(2.361)          & 2.253(0.653)          \\ \midrule
\multirow{5}{*}{$t(2)$}   & LS    & \textbf{18.765(9.483)}  & \textbf{1.711(0.563)} & \textbf{1.657(0.557)} & 0.868(0.226)          & 0.843(0.212)          \\
& LAD           & 21.285(4.814)          & 1.762(0.239)   & 1.709(0.237)   & \textbf{0.857(0.109)} & \textbf{0.815(0.117)} \\
& Huber         & 31.808(11.372)         & 2.304(0.588)   & 2.247(0.581)   & 1.031(0.223)          & 0.937(0.204)          \\
& Cauchy        & 78.870(25.555)         & 3.930(0.650)   & 3.875(0.649)   & 1.357(0.140)          & 1.034(0.110)          \\
& Tukey         & 180.530(36.958)        & 8.965(1.507)   & 8.883(1.497)   & 2.837(0.495)          & 1.904(0.339)          \\ \midrule
\multirow{5}{*}{$Cauchy(0,1)$} & LS            & 146.669(177.462)       & 4.229(2.136)   & 4.161(2.129)   & 1.463(0.405)          & 1.238(0.248)          \\
& LAD           & 22.448(8.524)          & 2.279(0.438)   & 2.225(0.431)   & 0.978(0.164)          & 0.859(0.155)          \\
& Huber & \textbf{20.960(15.823)} & \textbf{2.005(0.780)} & \textbf{1.956(0.778)} & \textbf{0.851(0.192)} & \textbf{0.744(0.122)} \\
& Cauchy        & 78.082(54.890)         & 4.079(1.420)   & 4.023(1.417)   & 1.344(0.289)          & 0.981(0.139)          \\
& Tukey         & 197.141(64.730)        & 8.053(1.696)   & 7.982(1.690)   & 2.375(0.419)          & 1.531(0.242)          \\ \midrule
\multirow{5}{*}{Mixture} & LS            & 450.810(193.680)	& 10.148(1.856)	& 10.056(1.856)	& 2.830(0.259)	& 2.113(0.128)        \\
& LAD           & 51.893(70.423) &	2.092(0.928)	&2.042(0.927)	&0.894(0.164)	&0.834(0.115)         \\
& Huber & \textbf{31.130(11.430)} &	\textbf{1.836(0.474)}&	\textbf{1.791(0.471)} &	\textbf{0.837(0.166)} &	\textbf{0.776(0.128)} \\
& Cauchy        & 79.468(26.835)	& 3.426(0.508)	&3.377(0.507)	&1.214(0.123)	&0.956(0.109)        \\
& Tukey         & 7063.935(20676.158)	&25.706(55.641)&25.634(55.634)&	2.770(1.459)&	1.684(0.315)         \\ \midrule

\multicolumn{2}{c}{Dopller, $n=512$}           & LS                     & LAD            & Huber          & Cauchy                & Tukey                 \\ \midrule
\multirow{5}{*}{$N(0,1)$} & LS    & \textbf{13.193(5.793)}  & \textbf{1.267(0.369)} & \textbf{1.226(0.364)} & \textbf{0.654(0.149)} & \textbf{0.633(0.128)} \\
& LAD           & 17.459(9.004)          & 1.458(0.529)   & 1.414(0.524)   & 0.715(0.191)          & 0.682(0.154)          \\
& Huber         & 19.160(5.569)          & 1.558(0.389)   & 1.516(0.383)   & 0.743(0.184)          & 0.699(0.183)          \\
& Cauchy        & 35.604(17.747)         & 2.243(0.649)   & 2.198(0.645)   & 0.922(0.157)          & 0.791(0.097)          \\
& Tukey         & 137.879(102.632)       & 5.814(2.568)   & 5.756(2.558)   & 1.787(0.492)          & 1.267(0.265)          \\ \midrule
\multirow{5}{*}{$t(2)$}        & LS            & \textbf{13.319(5.055)} & 1.374(0.276)   & 1.320(0.274)   & 0.743(0.120)          & 0.738(0.146)          \\
& LAD   & 16.387(7.046)           & \textbf{1.302(0.367)} & \textbf{1.266(0.362)} & \textbf{0.606(0.147)} & \textbf{0.551(0.133)} \\
& Huber         & 20.022(8.077)          & 1.558(0.477)   & 1.513(0.468)   & 0.720(0.195)          & 0.651(0.176)          \\
& Cauchy        & 28.573(11.641)         & 1.847(0.498)   & 1.808(0.496)   & 0.753(0.165)          & 0.630(0.137)          \\
& Tukey         & 537.837(1257.923)      & 11.205(16.875) & 11.143(16.858) & 2.317(1.715)          & 1.366(0.543)          \\  \midrule
\multirow{5}{*}{$Cauchy(0,1)$} & LS            & 463.967(861.745)       & 5.217(4.085)   & 5.150(4.076)   & 1.484(0.552)          & 1.238(0.345)          \\
& LAD   & \textbf{12.482(6.446)}  & \textbf{1.015(0.327)} & \textbf{0.985(0.326)} & \textbf{0.473(0.101)} & \textbf{0.424(0.069)} \\
& Huber         & 16.607(5.793)          & 1.214(0.307)   & 1.180(0.305)   & 0.540(0.094)          & 0.471(0.070)          \\
& Cauchy        & 29.285(16.817)         & 1.800(0.620)   & 1.767(0.620)   & 0.671(0.142)          & 0.525(0.081)          \\
& Tukey         & 133.563(95.067)        & 5.897(3.034)   & 5.843(3.021)   & 1.728(0.705)          & 1.123(0.376)          \\\midrule

\multirow{5}{*}{Mixture} & LS            & 127.257(44.378)&	6.153(1.271)&	6.063(1.269)	&2.394(0.303)	&2.014(0.150)      \\
& LAD           & 20.396(9.476)	& 1.420(0.446)	&1.382(0.443)	&0.668(0.156)	&0.619(0.123)       \\
& Huber & \textbf{14.738(4.872)}& 	\textbf{1.149(0.248)}	&\textbf{1.114(0.245)}	&\textbf{0.582(0.088)}	&\textbf{0.571(0.068)} \\
& Cauchy        & 53.697(28.314)&	2.410(0.814)&	2.373(0.811)&	0.873(0.201)&	0.697(0.125)        \\
& Tukey         & 137.383(83.918)&	5.270(1.839)&	5.215(1.831)&	1.679(0.452)&	1.221(0.338)   \\\bottomrule
\end{tabular}%
}
\end{table}
\clearpage

\begin{figure}[H]
	\centering
	\begin{subfigure}{\textwidth}
		\includegraphics[width=1\textwidth]{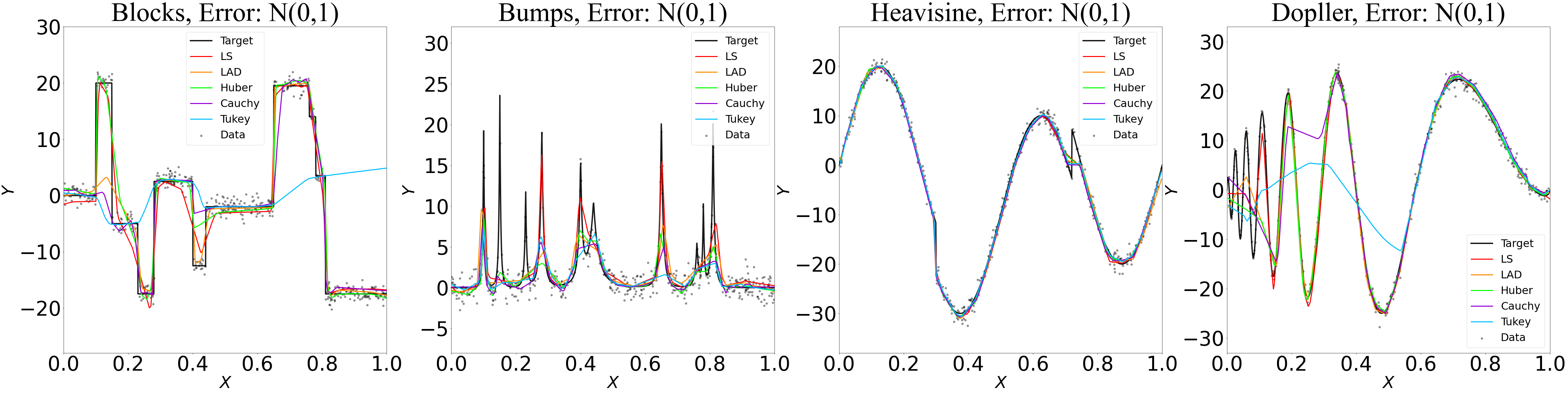}
	\end{subfigure}
	
	\begin{subfigure}{1\textwidth}
		\includegraphics[width=1\textwidth]{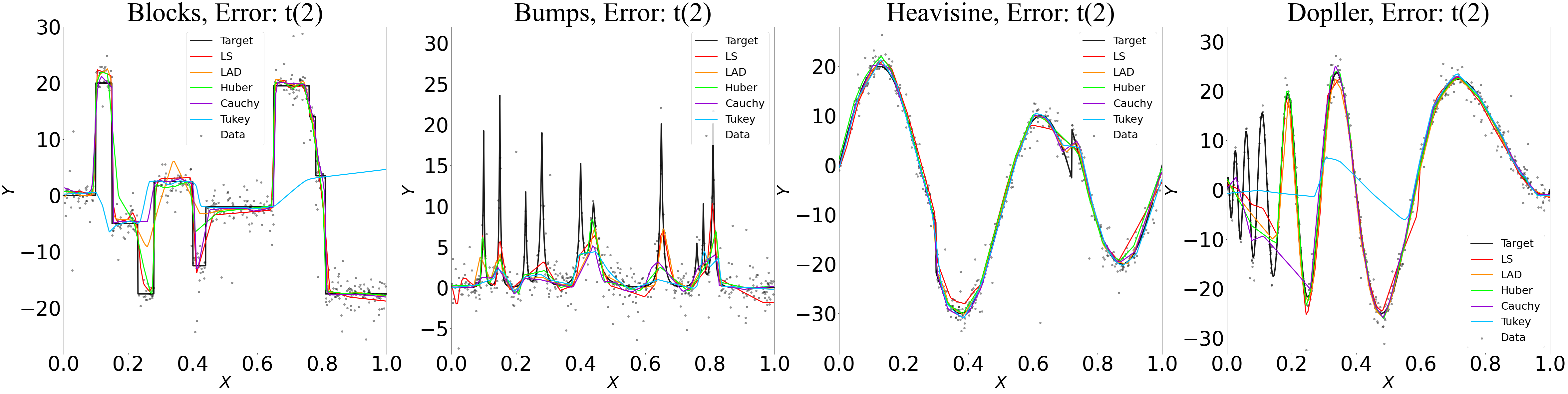}
	\end{subfigure}
	
	\begin{subfigure}{\textwidth}
		\includegraphics[width=\textwidth]{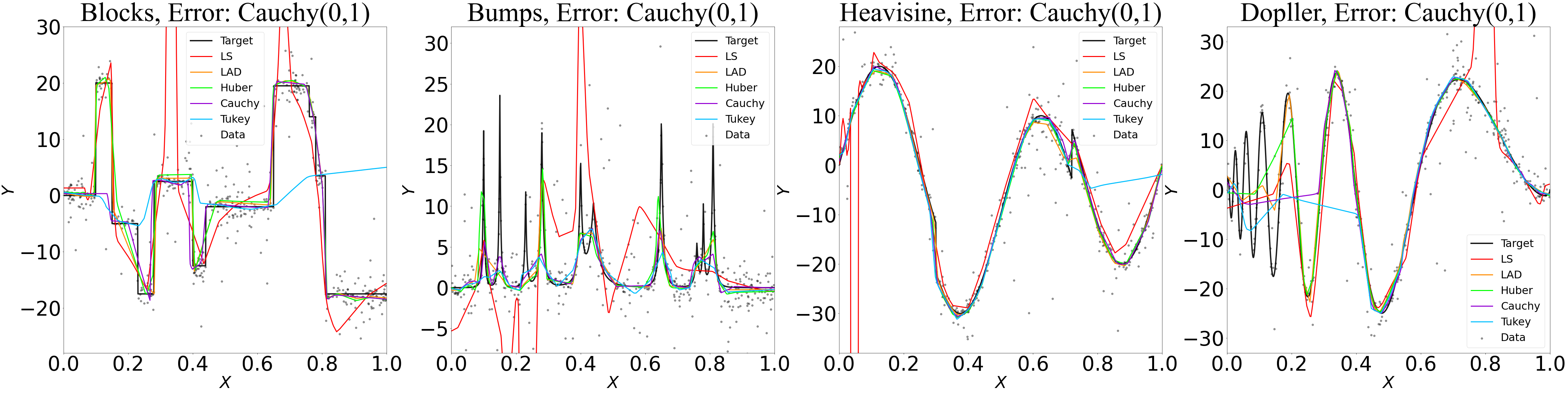}
	\end{subfigure}
	\begin{subfigure}{\textwidth}
		\includegraphics[width=\textwidth]{fitted_mix.png}
	\end{subfigure}
\caption{The fitted curves of univariate models. The training data is depicted as grey dots, the target function is depicted as black curves.
The colored curves represent the estimated curves
%stand for predictions of the trained estimators
under different loss functions.
Red line: LS; Orange line: LAD; Green line: Huber;  Purple line: Cauchy;
Blue line: Tukey.
From the top to the bottom, each row corresponds a  certain type of error: $N(0,1),t(2)$, $Cauchy(0,1)$ and normal mixture. From the left to right, each column corresponds a certain univariate model: ``Blocks'', ``Bumps'', ``Heavisine'' and ``Dopller''.}
	\label{fig:fitted}
\end{figure}

\begin{table}[H]

\setlength{\tabcolsep}{4pt} % Default value: 6pt
\renewcommand{\arraystretch}{1.2} % Default value: 1
\caption{The testing excess risks and standard deviation (in parentheses)  under different loss functions of the ERMs trained with respect to different loss functions using \textit{Nets-256}. The dimension of the input $d=4$, the training sample size $n=128,512$ and the number of replications $R=10$.}
\centering
\label{tab:2}
\resizebox{\textwidth}{!}{%
\begin{tabular}{@{}ccccccc@{}}
\toprule
\multicolumn{2}{c|}{$d=4, n=128$}              & \multicolumn{5}{c}{Testing Loss}                                                                                         \\ \midrule
Error                        & \multicolumn{1}{c|}{Training Loss} & LS                       & LAD                   & Huber                 & Cauchy                & Tukey                 \\ \midrule
\multirow{5}{*}{$N(0,1)$}      & \multicolumn{1}{c|}{LS}            & \textbf{13.678(0.674)}   & 1.336(0.046)          & 1.272(0.045)          & 0.814(0.027)          & 0.923(0.031)          \\
& \multicolumn{1}{c|}{LAD}           & 14.832(1.483)            & 1.368(0.119)          & 1.304(0.114)          & 0.821(0.079)          & 0.926(0.096)          \\
& \multicolumn{1}{c|}{Huber}         & 14.033(0.463)            & \textbf{1.325(0.066)} & \textbf{1.262(0.064)} & \textbf{0.799(0.041)} & \textbf{0.903(0.047)} \\
& \multicolumn{1}{c|}{Cauchy}        & 22.141(2.693)            & 1.505(0.111)          & 1.441(0.108)          & 0.834(0.065)          & 0.923(0.079)          \\
& \multicolumn{1}{c|}{Tukey}         & 59.162(6.536)            & 2.340(0.185)          & 2.275(0.183)          & 0.980(0.062)          & 0.977(0.058)          \\ \midrule
\multirow{5}{*}{$t(2)$}        & \multicolumn{1}{c|}{LS}            & 19.227(3.166)            & 1.725(0.242)          & 1.660(0.238)          & 0.961(0.112)          & 1.029(0.109)          \\
& \multicolumn{1}{c|}{LAD}           & \textbf{18.103(2.433)}   & \textbf{1.483(0.149)} & \textbf{1.422(0.146)} & \textbf{0.819(0.074)} & \textbf{0.879(0.076)} \\
& \multicolumn{1}{c|}{Huber}         & 20.011(1.943)            & 1.684(0.114)          & 1.620(0.112)          & 0.922(0.059)          & 0.981(0.062)          \\
& \multicolumn{1}{c|}{Cauchy}        & 21.666(1.591)            & 1.570(0.136)          & 1.510(0.134)          & 0.842(0.072)          & 0.897(0.074)          \\
& \multicolumn{1}{c|}{Tukey}         & 65.676(14.472)           & 2.509(0.381)          & 2.446(0.379)          & 0.998(0.103)          & 0.955(0.077)          \\ \midrule
\multirow{5}{*}{$Cauchy(0,1)$} & \multicolumn{1}{c|}{LS}            & 2837.600(5740.160)       & 6.181(5.959)          & 6.116(5.955)          & 1.386(0.380)          & 1.223(0.223)          \\
& \multicolumn{1}{c|}{LAD}           & \textbf{18.400(113.676)} & 2.098(0.963)          & 2.041(0.962)          & 0.947(0.150)          & 0.942(0.096)          \\
& \multicolumn{1}{c|}{Huber}         & 130.800(422.469)         & 2.623(1.663)          & 2.563(1.663)          & 1.033(0.160)          & 1.012(0.112)          \\
& \multicolumn{1}{c|}{Cauchy}        & 32.026(8.921)            & \textbf{1.656(0.283)} & \textbf{1.601(0.280)} & \textbf{0.787(0.097)} & \textbf{0.802(0.089)} \\
& \multicolumn{1}{c|}{Tukey}         & 44.800(43.074)           & 2.541(0.236)          & 2.485(0.236)          & 0.971(0.052)          & 0.908(0.033)          \\ \midrule
\multicolumn{2}{c}{ $d=4, n=512$}           & LS                     & LAD            & Huber          & Cauchy                & Tukey                 \\ \midrule
\multirow{5}{*}{$N(0,1)$}      & \multicolumn{1}{c|}{LS}           & \textbf{5.096(0.385)}       & 0.678(0.040) & 0.634(0.037) & 0.459(0.027) & 0.527(0.033) \\
& \multicolumn{1}{c|}{LAD}          & 5.821(0.561)       & 0.729(0.056) & 0.683(0.052) & 0.484(0.042) & 0.552(0.052) \\
& \multicolumn{1}{c|}{Huber}         & 5.172(0.449)       & \textbf{0.678(0.025)} & \textbf{0.634(0.025)} & \textbf{0.457(0.015)} & \textbf{0.525(0.018)} \\
& \multicolumn{1}{c|}{Cauchy}       & 8.733(0.820)       & 0.834(0.037) & 0.787(0.035) & 0.517(0.035) & 0.582(0.045) \\
& \multicolumn{1}{c|}{Tukey}         & 23.338(5.598)      & 1.166(0.134) & 1.120(0.134) & 0.561(0.031) & 0.581(0.031) \\ \midrule
\multirow{5}{*}{$t(2)$}        & \multicolumn{1}{c|}{LS}           & 11.257(8.160)      & 1.028(0.149) & 0.976(0.147) & 0.624(0.050) & 0.696(0.054) \\
& \multicolumn{1}{c|}{LAD}          & 7.113(1.641)       & \textbf{0.753(0.104)} & \textbf{0.710(0.101)} & \textbf{0.464(0.057)} & \textbf{0.511(0.064)} \\
& \multicolumn{1}{c|}{Huber}         & \textbf{5.934(0.493)}       & 0.764(0.077) & 0.718(0.073) & 0.496(0.055) & 0.557(0.068) \\
& \multicolumn{1}{c|}{Cauchy}        & 9.490(0.809)       & 0.820(0.026) & 0.777(0.025) & 0.477(0.018) & 0.515(0.024) \\
& \multicolumn{1}{c|}{Tukey}         & 25.927(5.004)      & 1.273(0.145) & 1.227(0.143) & 0.590(0.051) & 0.594(0.054) \\ \midrule
\multirow{5}{*}{$Cauchy(0,1)$} & \multicolumn{1}{c|}{LS}           & 3705.178(7424.642) & 7.270(5.723) & 7.201(5.717) & 1.567(0.767) & 1.292(0.355) \\
& \multicolumn{1}{c|}{LAD}          & \textbf{9.489(9.208)}  & 0.917(0.128) & 0.872(0.127) & 0.520(0.060) & 0.558(0.069) \\
& \multicolumn{1}{c|}{Huber}        & 10.836(5.023)      & 0.966(0.101) & 0.919(0.100) & 0.552(0.047) & 0.590(0.047) \\
& \multicolumn{1}{c|}{Cauchy}       & 11.392(6.910)      & \textbf{0.873(0.113)}& \textbf{0.830(0.112)} & \textbf{0.476(0.049)} & \textbf{0.501(0.058)} \\
& \multicolumn{1}{c|}{Tukey}        & 34.442(6.141)      & 1.445(0.152) & 1.399(0.149) & 0.615(0.064) & 0.598(0.069) \\ \bottomrule

\end{tabular}%
}
\end{table}

	\section{Concluding remarks}
	\label{conclusion}
	In this work,
%we study the proties of robust nonparametric estimation using deep neural networks %for regression models with heavy tailed error distributions.
we establish the non-asymptotic error bounds for a class of robust nonparametric regression estimators using deep neural networks with ReLU activation under suitable smoothness conditions on the regression function and mild conditions on the error term. In particular, we only assume that the error distribution has a finite $p$-th moment with $p>1$. We also show that the deep robust regression is able to circumvent the curse of dimensionality when the predictor is supported on an approximate lower-dimensional set. An important feature of our error bound is that, for ReLU neural networks with network width and network size (number of parameters) no more than the order of the square of the dimensionality $d$ of the predictor, our excess risk bounds depend sub-linearly on the $d$. Our assumption relaxes the exact manifold support assumption, which could be restrictive and unrealistic in practice. We also relax several crucial assumptions on the data distribution, the target regression function and the neural networks required in the recent literature. Our simulation studies demonstrate the advantages of using robust loss functions over the least squares loss function in the context of deep nonparametric regression in terms of testing excess risks.
	
We only addressed the issue of robustness with respect to the heavy tailed errors or outliers in the response variable. There are some other important sources of non- robustness in nonparametric regression. For example, it is possible that outliers exist in the observed values of the predictors, or that there are contaminations in both the responses and the predictors. How to address these important practical issues in the context of nonparametric regression using deep neural networks deserve further studies in the future.
	
\section*{Acknowledgements}
The work of Y. Jiao is supported in part by the National Science Foundation of China grant 11871474 and by the research fund of KLATASDSMOE of China. The work of Y. Lin is supported by the Hong Kong Research Grants Council (Grant No. 14306219 and 14306620) and Direct Grants for Research, The Chinese University of Hong Kong. The work of J. Huang is partially supported by the U.S. NSF grant DMS-1916199.

%\clearpage
\bigskip
\bibliographystyle{imsart-nameyear}
% Style BST file (imsart-number.bst or imsart-nameyear.bst)
\bibliography{ddr_arXiv.bib}    % Bibliography file (usually '*.bib')
	
\clearpage
	
\begin{appendix}
\section{A brief review of some recent results on neural network approximation}
%BRIEF REVIEW OF SOME RECENT RESULTS ON NEURAL NETWORK %APPROXIMATION}
In this paper, we only assume that $f^*$ is uniformly continuous (Assumption \ref{A3}), which is a mild assumption on the continuity of the unknown target function $f^*$, as
the existing works generally assume %relativele
stronger smoothness assumptions on $f^*$.
For example, in  nonparametric  least squares regression, \citet{stone1982optimal} and \citet{bauer2019deep} require that $f^*$ is $\beta$-H\"older smooth with $\beta \ge 1$, i.e.,  all its partial derivatives  up to order $\lfloor\beta\rfloor$ exist and the partial derivatives of order $\lfloor\beta\rfloor$ are $\lfloor\beta\rfloor-\beta$ H\"older continuous.
	%As discussed in section \ref{review},
	\citet{farrell2021deep} assume that $f^*$ lies in a Sobolev ball with smoothness $\beta\in \mathbb{N}^+$, i.e. $f^*(x)\in\mathcal{W}^{\beta,\infty}([-1,1]^d)$.
	Recent approximation theories on different functional spaces, such as  Korobov spaces, Besev spaces or function space with $f^*\in C^\alpha[0,1]^d$ with $\alpha\geq1$
	can be found in  \cite{liang2016deep}, \cite{lu2017expressive}, \cite{yarotsky2017error},  \cite{mohri2018foundations}, \cite{suzuki2018adaptivity},   among others.

	The function class $\mathcal{F}_\phi$ in this paper is constructed by the feedforward neural networks with the ReLU activation function.
	Some recent  results on the approximation theory of deep neural networks can be found in
	\cite{devore2010approximation,
		hangelbroek2010nonlinear,
		lin2014almost,
		yarotsky2017error,
		yarotsky2018optimal,
		lu2017expressive,
		raghu2017expressive,
		shen2019nonlinear,
		shen2019deep,
		nakada2019adaptive,
		chen2019efficient}.
	
	An important % remarkable
	result proved by \cite{yarotsky2017error} is the following,
	for any $\varepsilon\in(0,1)$, any $d,\beta$, and any $f_0$ in the Sobolev ball $\mathcal{W}^{\beta,\infty}([0,1]^d)$ with $\beta > 0$,
	%  in (\ref{sob}) below,
	there exists a ReLU network $\hat{f}$ with depth $\mathcal{D}$ at most $c\{\log(1/\varepsilon)+1\}$, size $\mathcal{S}$ and number of neurons $\mathcal{U}$ at most $c\varepsilon^{-d/\beta}\{\log(1/\varepsilon)+1\}$ such that $\Vert \hat{f}-f_0\Vert_\infty=\max_{x\in[0,1]^d} \vert \hat{f}(x)-f_0(x)\vert\leq \varepsilon$, where $c$ is some constant depending on $d$ and $\beta$.
	In particular, it is required that the constant $c=O(2^d)$, an exponential rate of $d$,  due to the technicality in the proof.  	
	The main idea of % the proof in
	\cite{yarotsky2017error} is to show that, small neural networks can approximate polynomials locally, and stacked neural networks (by $2^d$ small sub-networks) can further approximate smooth function by approximating its Taylor expansions.
	Later, \citet{yarotsky2018optimal} gave the optimal rate of approximation for general continuous functions by deep ReLU networks, in terms of the network size $\mathcal{S}$ and the modulus of continuity of $f^*$.  It was shown that $\inf_{f\in\mathcal{F}_\phi}\Vert f-f^*\Vert_\infty\leq c_1\omega_{f^*}(c_2\mathcal{S}^{-p/d})$ for some $p\in[1,2]$ and some constants $c_1,c_2$ possibly depending on $d,p$ but not $\mathcal{S}, f^*$.
	The upper bound holds for any $p\in(1,2]$ if the network $\mathcal{F}_\phi=\mathcal{F}_{\mathcal{D},\mathcal{W},\mathcal{U},\mathcal{S},\mathcal{B}}$ satisfies $\mathcal{D}\geq c_3\mathcal{S}^{p-1}/\log(\mathcal{S})$ for some constant $c_3$ possibly depending on $p$ and $d$.
	
	Based on a method similar to that in %proof  strategy as
	\citet{yarotsky2017error}, approximation of
	%rates to approximate
	a target function  on a low-dimensional manifold $\mathcal{M}$ of $[0,1]^d$
	with deep neural nets was studied by \citet{nakada2019adaptive, chen2019efficient}.  The approximation error of ReLU networks to $f^*$ defined
	on an exact manifold  $\mathcal{M}$ will no longer depend on $d$, instead it depends on the intrinsic dimension
	%or the Minkowski dimension
	$d_0< d $ of the manifold $\mathcal{M}$. Specifically, the approximation rate is shown to be  $c\mathcal{S}^{-\beta/d_0}$, where $f^*$ is often assumed to be a $\beta$-H\"older smooth function defined on $[0,1]^d$ and $c$ is some constant depending on $d, d_0$ and $\beta$. This helps alleviate the curse of dimensionality in approximating high-dimensional functions
	% nonparametric regression
	with deep ReLU neural networks.
	
	It is worth noting that	other related works, including \citet{chen2019nonparametric},
	\citet{chen2019efficient},
	\citet{nakada2019adaptive},
	\citet{schmidt2019deep} and
	\citet{schmidt2020nonparametric}, rely on a similar approximation construction of \cite{yarotsky2017error}.
	A common feature
	%characteristic
	of the
	% important
	results from these works is that,  the prefactor of the approximation error is of the order $O(2^d)$
	unless  the size $\mathcal{S}$ of the network  has to grow exponentially with respect to the dimension $d$.
	%to render the approximate error;
	Unfortunately, a prefactor of the order $O(2^d)$ is extremely big for a large $d$ in the high-dimensional settings, which can destroy the approximation error bound even for networks with a large size. For example, each handwritten digit picture in the MNIST dataset \citep{mnist2010}
	is of dimension $d=28\times28$ and the sample size $n$ is about $70,000$, but $2^d$ is approximately $10^{236}$.

	%They also pointed that the condition $n \gtrsim d^{d}$ is not necessary for learning the additive model, see for example \cite{bach2017}.

	%	Different from the aforementioned approximation results

	Lastly, it follows from Proposition 1 of \cite{yarotsky2017error} that, for a neural network, in terms of its computational power and complexity, there is no substantial difference in using the ReLU activation function and  other piece-wise linear activation function with finitely many breakpoints. To elaborate, let $\xi: \mathbb{R}\to\mathbb{R}$ be any continuous piece-wise linear function with $M$ breakpoints ($1\leq M<\infty$). If a network $f_\xi$ is activated by $\xi$, of depth $\mathcal{D}$, size $\mathcal{S}$ and the number of neurons $\mathcal{U}$, then there exists a ReLU activated network with depth $\mathcal{D}$, size not more than $(M+1)^2\mathcal{S}$, the number of neurons not more than $(M+1)\mathcal{U}$, that computes the same function as $f_\xi$. Conversely, let $f_\sigma$ be a ReLU activated network  of depth $\mathcal{D}$, size $\mathcal{S}$ and the number of neurons $\mathcal{U}$, then there exists a network with  activation function $\xi$, of depth $\mathcal{D}$, size $4\mathcal{S}$ and
	the number of neurons $2\mathcal{U}$ that computes the same function $f_\sigma$ on a bounded subset of $\mathbb{R}^d$.

\section{Proofs}
		In this part, we provide the proofs of Lemmas \ref{lemma1}, \ref{lemma2} and \ref{quadapprox}, Theorems \ref{thm2} and \ref{thm3} and {\color{black} Corollary \ref{c-all}}.
			
		\subsection{Proof of Lemma \ref{lemma1}}
	\begin{comment}
		\begin{proof}
	By the definition of the empirical risk minimizer, for any $f\in\mathcal{F}_n$, we have $\mathcal{R}_n(\hat{f}_n)\leq \mathcal{R}_n(f)$. Therefore,
	\begin{align*} \mathcal{R}(\hat{f}_n)-\mathcal{R}(f_0)=&\mathcal{R}(\hat{f}_n)-\mathcal{R}_n(\hat{f}_n)+\mathcal{R}_n(\hat{f}_n)-\mathcal{R}_n(f)+\mathcal{R}_n(f)-\mathcal{R}(f)+\mathcal{R}(f)-\mathcal{R}(f_0)\\
	\leq&\mathcal{R}(\hat{f}_n)-\mathcal{R}_n(\hat{f}_n)+\mathcal{R}_n(f)-\mathcal{R}(f)+\mathcal{R}(f)-\mathcal{R}(f_0)\\
	=&\big\{\mathcal{R}(\hat{f}_n)-\mathcal{R}_n(\hat{f}_n)\big\}+\big\{\mathcal{R}_n(f)-\mathcal{R}(f)\big\}+\big\{\mathcal{R}(f)-\mathcal{R}(f_0)\big\}\\
	\leq& 2\sup_{f\in\mathcal{F}_n}\vert \mathcal{R}(f)-\mathcal{R}_n(f)\vert+\big\{\mathcal{R}(f)-\mathcal{R}(f_0)\big\}.
	\end{align*}
	Since the above inequality holds for any $f\in\mathcal{F}_n$, Lemma \ref{lemma1} is proved by choosing $f$  satisfying $f\in\arg\inf_{f\in\mathcal{F}_n}\mathcal{R}(f)$.
	\end{proof}
	\end{comment}

		\begin{proof}
			Let $S=\{Z_i=(X_i,Y_i)\}_{i=1}^n$ be a sample form the distribution of $Z=(X,Y)$ and $S^\prime=\{Z_i^\prime=(X^\prime_i,Y^\prime_i)\}_{i=1}^n$ be another sample independent of $S$. Define $g(f,Z_i)=L(f(X_i),Y_i)-L(f^*(X_i),Y_i)$ for any $f$ and sample $Z_i$. Note that the empirical risk minimizer $\hat{f}_\phi$ defined in Lemma \ref{lemma1} depends on the sample $S$, and its excess risk is $\mathbb{E}_{S^\prime} \{\sum_{i=1}^ng(\hat{f}_\phi,Z_i^\prime)/n\}$ and its expected excess risk is
			\begin{equation}\label{ine1}
				\mathbb{E}\big\{\mathcal{R}(\hat{f}_\phi)-\mathcal{R}(f^*)\big\}=\mathbb{E}_{S}[\mathbb{E}_{S^\prime} \{\frac{1}{n}\sum_{i=1}^ng(\hat{f}_\phi,Z_i^\prime)\}].
			\end{equation}
			Next we will take 3 steps to complete the proof of Lemma \ref{lemma1}.
			
			\subsubsection*{Step 1: Excess risk decomposition}
			
			Define the `best in class' estimator $f^*_\phi$ as the  estimator in the function class $\mathcal{F}_\phi=\mathcal{F}_{\mathcal{D},\mathcal{W},\mathcal{U},\mathcal{S},\mathcal{B}}$ with minimal $L$ risk:
			\begin{equation*}
				f^*_\phi=\arg\min_{f\in\mathcal{F}_{\phi}} \mathcal{R}(f).
			\end{equation*}
			The approximation error of $f^*$  is $\mathcal{R}(f^*_\phi)-\mathcal{R}(f^*)$. Note that the approximation error only depends on the function class $\mathcal{F}_{\mathcal{D},\mathcal{W},\mathcal{U},\mathcal{S},\mathcal{B}}$ and the distribution of data. By the definition of empirical risk minimizer, we have
			\begin{equation}  \label{ine2}
				\mathbb{E}_{S}\{\frac{1}{n}\sum_{i=1}^ng(\hat{f}_\phi,Z_i)\}\leq \mathbb{E}_S\{\frac{1}{n}\sum_{i=1}^ng(f^*_\phi,Z_i)\}.
			\end{equation}
			Multiply 2 by the both sides of (\ref{ine2}) and add it up with (\ref{ine1}), we have
			\begin{align} \nonumber
				\mathbb{E}\big\{\mathcal{R}(\hat{f}_\phi)-\mathcal{R}(f^*)\big\}&\leq\mathbb{E}_{S}\Big[ \frac{1}{n}\sum_{i=1}^n\big\{-2g(\hat{f}_\phi,Z_i)+\mathbb{E}_{S^\prime}g(\hat{f}_\phi,Z_i^\prime)\big\}\Big]+ 2\mathbb{E}_S\{\frac{1}{n}\sum_{i=1}^ng(f^*_\phi,Z_i)\}\\
				&\leq\mathbb{E}_{S}\Big[ \frac{1}{n}\sum_{i=1}^n\big\{-2g(\hat{f}_\phi,Z_i)+\mathbb{E}_{S^\prime}g(\hat{f}_\phi,Z_i^\prime)\big\}\Big]+ 2\big\{\mathcal{R}(f^*_\phi)-\mathcal{R}(f^*)\big\}. \label{bound0}
			\end{align}
			It is seen that the excess risk is upper bounded by the sum of a expectation of a stochastic term and approximation error.
			\subsubsection*{Step 2: Bounding the stochastic term}
			Next, we will focus on giving an upper bound of the first term on the right-hand side in (\ref{bound0}), and handle it with truncation and classical chaining technique of empirical processes. In the following, for ease of presentation, we write $G(f,Z_i)=\mathbb{E}_{S^\prime}\{g(f,Z_i^\prime)\}-2g(f,Z_i)$ for $f\in\mathcal{F}_\phi$.
			
			%%%%%%%%%%%%%%%%%%%%
			%Given a sequence $x=(x_1,...,x_n)\in\mathcal{X}\subseteq\mathbb{R}^n$, let  %$\mathcal{F}_\phi|_x=\{(f(x_1),...,f(x_n):f\in\mathcal{F}_\phi\}$ be the  subset of %$\mathbb{R}^{n}$. For a positive number $\delta$, let %$\mathcal{N}(\delta,\Vert\cdot\Vert_\infty,\mathcal{F}_\phi|_x)$ be the covering number of %$\mathcal{F}_\phi|_x$ under the norm $\Vert\cdot\Vert_\infty$ with radius $\delta$.
			%Define the uniform covering number %$\mathcal{N}_n=\mathcal{N}(\delta,\Vert\cdot\Vert_\infty,\mathcal{F}_\phi)$ to be the %maximum
			%	over all $x\in\mathcal{X}$ of the covering number %$\mathcal{N}(\delta,\Vert\cdot\Vert_\infty,\mathcal{F}_\phi|_x)$, i.e.,
			%	$$\mathcal{N}(\delta,\Vert\cdot\Vert_\infty,\mathcal{F}_\phi)=
			%\max\{\mathcal{N}(\delta,\Vert\cdot\Vert_\infty,\mathcal{F}_\phi|_x):x\in\mathcal{X}\}.$$

			Given a $\delta$-uniform covering of $\mathcal{F}_\phi$, we denote the centers of the balls by $f_j,j=1,2,\ldots,\mathcal{N}_{2n},$ where $\mathcal{N}_{2n}=\mathcal{N}_{2n}(\delta,\Vert\cdot\Vert_\infty,\mathcal{F}_\phi)$ is the uniform covering number with radius $\delta$ ($\delta<\mathcal{B}$) under the norm $\Vert\cdot\Vert_\infty$, where $\mathcal{N}_{2n}(\delta,\Vert\cdot\Vert_\infty,\mathcal{F}_\phi)$
			is defined in (\ref{ucover}).
			By the definition of covering, there exists a (random) $j^*$ such that $\Vert\hat{f}_\phi(x) -f_{j^*}(x)\Vert_\infty\leq\delta$ on $x=(X_1,\ldots,X_n,X_1^\prime,\ldots,X_n^\prime)\in\mathcal{X}^{2n}$. By the assumptions that $\Vert f^*\Vert_\infty, \Vert f_j\Vert_\infty\leq \mathcal{B}$, $\mathbb{E}\vert Y_i\vert^p<\infty$, and Assumption \ref{A1}, we have
			\begin{align*}
				\vert g(\hat{f}_\phi,Z_i)-g(f_{j^*},Z_i) \vert &= \vert  L(\hat{f}_\phi(X_i),Y_i)-L(f_{j^*}(X_i),Y_i) \vert\\
				&\leq \lambda_L \vert \hat{f}_\phi(X_i)-f_{j^*}(X_i)\vert\\
				&\leq \lambda_L \delta,
			\end{align*}
			and
			\begin{align*}
				\frac{1}{n}\sum_{i=1}^n \mathbb{E}_{S}\big\{g(\hat{f}_\phi,Z_i)\}	\leq\frac{1}{n}\sum_{i=1}^n\mathbb{E}_{S}\big\{&g(f_{j^*},Z_i)\}+\frac{1}{n}\sum_{i=1}^n\mathbb{E}_{S}\big\{g(\hat{f}_\phi,Z_i)-g(f_{j^*},Z_i)\}\\
				\leq\frac{1}{n}\sum_{i=1}^n\mathbb{E}_{S}\big\{&g(f_{j^*},Z_i)\}+\frac{1}{n}\sum_{i=1}^n\mathbb{E}_{S}\vert g(\hat{f}_\phi,Z_i)-g(f_{j^*},Z_i)\vert\\
					\leq\frac{1}{n}\sum_{i=1}^n\mathbb{E}_{S}\big\{&g(f_{j^*},Z_i)\}+\lambda_L\delta.
			\end{align*}
		Then we have
		\begin{align*}
			\vert G(\hat{f}_\phi,Z_i)-G(f_{j^*},Z_i)\vert&=\vert\mathbb{E}_{S^\prime}\{g(\hat{f}_\phi,Z_i^\prime)\}-2g(\hat{f}_\phi,Z_i)-\mathbb{E}_{S^\prime}\{g(f_{j^*},Z_i^\prime)\}+2g(f_{j^*},Z_i)\vert\\
			&\leq\vert\mathbb{E}_{S^\prime}\{g(\hat{f}_\phi,Z_i^\prime)\}-\mathbb{E}_{S^\prime}\{g(f_{j^*},Z_i^\prime)\}\vert+\vert2g(f_{j^*},Z_i)-2g(\hat{f}_\phi,Z_i)\vert\\
			&\leq 3\lambda_L\delta,
		\end{align*}
		and
		\begin{equation} \label{bound1}
				\mathbb{E}_{S}\Big[ \frac{1}{n}\sum_{i=1}^n G(\hat{f}_\phi,Z_i) \Big]
				\leq\mathbb{E}_{S}\Big[ \frac{1}{n}\sum_{i=1}^n G(f_{j^*},Z_i) \Big]+3\lambda_L\delta.
		\end{equation}
			
			Let $\beta_n\geq \mathcal{B}\geq1$ be a positive number who may depend on the sample size $n$. Denote $T_{\beta_n}$ as the truncation operator at level $\beta_n$, i.e., for any $Y\in\mathbb{R}$,
			$T_{\beta_n}Y=Y$ if $\vert Y\vert\leq\beta_n$ and $T_{\beta_n}Y= \beta_n\cdot {\rm sign}(Y)$ otherwise. Define the function $f^*_{\beta_n}:\mathcal{X}\to\mathbb{R}$ pointwisely by
			$$f^*_{\beta_n}(x)=\arg\min_{f(x):\Vert f\Vert_\infty\leq\beta_n} \mathbb{E}\big\{L(f(X),T_{\beta_n}Y)|X=x\big\},$$
			for each $x\in\mathcal{X}$.
			Besides, recall that $\Vert f^*\Vert_\infty\leq\mathcal{B}\leq\beta_n$ and
			$$f^*(x)=\arg\min_{f(x):\Vert f\Vert_\infty\leq\beta_n} \mathbb{E}\big\{L(f(X),Y)|X=x\big\}.$$
			Then for any $f$ satisfying $\Vert f \Vert\infty\leq\beta_n$, the definition above implies that  $\mathbb{E}\{L(f^*_{\beta_n}(X_i),T_{\beta_n}Y_i)\}\leq \mathbb{E}\{L(f(X_i),T_{\beta_n}Y_i)\}$ and $\mathbb{E}\{L(f^*(X_i),Y_i)\}\leq \mathbb{E}\{L(f(X_i),Y_i)\}$.

			 For any $f\in\mathcal{F}_\phi$, we let $g_{\beta_n}(f,Z_i)=L(f(X_i),T_{\beta_n}Y_i)-L(f^*_{\beta_n}(X_i),T_{\beta_n}Y_i)$.
			 Then we have
			\begin{align*}
				\mathbb{E}\{g(f,Z_i)\} &= \mathbb{E}\{g_{\beta_n}(f,Z_i)\}+ \mathbb{E}\{L(f(X_i),Y_i)-L(f(X_i),T_{\beta_n}Y_i)\}\\
				&\qquad\qquad\qquad\quad+\mathbb{E}\{L(f^*_{\beta_n}(X_i),T_{\beta_n}Y_i)-L(f^*(X_i),T_{\beta_n}Y_i)\}\\
				&\qquad\qquad\qquad\quad +\mathbb{E}\{L(f^*(X_i),T_{\beta_n}Y_i)-L(f^*(X_i),Y_i)\}\\
				&\leq \mathbb{E}\{g_{\beta_n}(f,Z_i)\}+ \mathbb{E}\vert L(f(X_i),Y_i)-L(f(X_i),T_{\beta_n}Y_i)\vert\\
				&\qquad\qquad\qquad\quad+\mathbb{E}\vert L(f^*(X_i),T_{\beta_n}Y_i)-L(f^*(X_i),Y_i)\vert\\
				&\leq \mathbb{E}\{g_{\beta_n}(f,Z_i)\} +2\lambda_L\mathbb{E}\{\vert T_{\beta_n}Y_i-Y_i\vert\}\\
				&\leq \mathbb{E}\{ g_{\beta_n}(f,Z_i)\}+2\lambda_L\mathbb{E}\big\{\vert \vert Y_i\vert I(\vert Y_i\vert>\beta_n)\big\}\\
				&\leq  \mathbb{E}\{ g_{\beta_n}(f,Z_i)\} + 2\lambda_L\mathbb{E}\{\vert Y_i\vert \vert Y_i\vert^{p-1}/\beta_n^{p-1}\}\\
				&\leq  \mathbb{E}\{ g_{\beta_n}(f,Z_i)\} + 2\lambda_L\mathbb{E}\vert Y_i\vert^{p}/\beta_n^{p-1}.
			\end{align*}
			By Assumption 2, the response $Y$ has finite $p$-moment and thus $\mathbb{E}\vert Y_i\vert^{p}<\infty$. Similarly,
			\begin{align*}
			\mathbb{E}\{g_{\beta_n}(f,Z_i)\}&= \mathbb{E}\{g(f,Z_i)\} + \mathbb{E}\{L(f^*(X_i),Y_i)-L(f^*_{\beta_n}(X_i),Y_i)\}\\
				&\qquad\qquad\qquad+\mathbb{E}\{L(f(X_i),T_{\beta_n}Y_i)-L(f(X_i),Y_i)\}\\
				&\qquad\qquad\qquad +\mathbb{E}\{L(f^*_{\beta_n}(X_i),Y_i)-L(f^*_{\beta_n}(X_i),T_{\beta_n}Y_i)\}\\
				&\leq \mathbb{E}\{g(f,Z_i)\} + \mathbb{E}\vert L(f(X_i),T_{\beta_n}Y_i)-L(f(X_i),Y_i)\vert\\
				&\qquad\qquad\qquad +\mathbb{E}\vert L(f^*_{\beta_n}(X_i),Y_i)-L(f^*_{\beta_n}(X_i),T_{\beta_n}Y_i)\vert\\
				&\leq \mathbb{E}\{ g(f,Z_i)\} + 2\lambda_L\mathbb{E}\vert Y_i\vert^{p}/\beta_n^{p-1}.
			\end{align*}
			Note that above inequalities also hold for $g(f,Z_i^\prime)$ and $g_{\beta_n}(f,Z_i^\prime)$.  Then for any $f\in\mathcal{F}_\phi$, define $G_{\beta_n}(f,Z_i)=\mathbb{E}_{S^\prime}\{g_{\beta_n}(f,Z_i^\prime)\}-2g_{\beta_n}(f,Z_i)$ and we have
			\begin{equation} \label{bound2}
				\mathbb{E}_{S}\Big[ \frac{1}{n}\sum_{i=1}^n G(f_{j^*},Z_i) \Big]
				\leq\mathbb{E}_{S}\Big[ \frac{1}{n}\sum_{i=1}^n G_{\beta_n}(f_{j^*},Z_i) \Big]+6\lambda_L\mathbb{E}\vert Y_i\vert^{p}/\beta_n^{p-1}.
			\end{equation}

			Besides, by Assumption \ref{A1}, for any $f\in\mathcal{F}_\phi$ we have $\vert g_{\beta_n}(f,Z_i)\vert\leq 4\lambda_L\beta_n$ and $\sigma^2_g(f):={\rm Var}(g_{\beta_n}(f,Z_i))\leq\mathbb{E} \{g_{\beta_n}(f,Z_i)^2\}\leq 4\lambda_L\beta_n\mathbb{E}\{g_{\beta_n}(f,Z_i)\}$. For each $f_j$ and any $t>0$, let $u=t/2+{\sigma_g^2(f_j)}/(8\lambda_L\beta_n)$,
			%similar to the proof in 11.3 of \cite{gyorfi2006distribution},
			by  the Bernstein inequality,
			\begin{align*}
				&P\Big\{\frac{1}{n}\sum_{i=1}^nG_{\beta_n}(f_j,Z_i)>t\Big\}\\
				=&P\Big\{\mathbb{E}_{S^\prime} \{g_{\beta_n}(f_j,Z_i^\prime)\}-\frac{2}{n}\sum_{i=1}^ng_{\beta_n}(f_j,Z_i)>t\Big\}\\
				=&P\Big\{\mathbb{E}_{S^\prime} \{g_{\beta_n}(f_j,Z_i^\prime)\}-\frac{1}{n}\sum_{i=1}^ng_{\beta_n}(f_j,Z_i)>\frac{t}{2}+\frac{1}{2}\mathbb{E}_{S^\prime} \{g_{\beta_n}(f_j,Z_i^\prime)\}\Big\}\\
				\leq& P\Big\{\mathbb{E}_{S^\prime} \{g_{\beta_n}(f_j,Z_i^\prime)\}-\frac{1}{n}\sum_{i=1}^ng_{\beta_n}(f_j,Z_i)>\frac{t}{2}+\frac{1}{2}\frac{\sigma_g^2(f_j)}{4\lambda_L\beta_n}\}\Big\}\\
				\leq& \exp\Big( -\frac{nu^2}{2\sigma_g^2(f_j)+16u\lambda_L\beta_n/3}\Big)\\
				\leq& \exp\Big( -\frac{nu^2}{16u\lambda_L\beta_n+16u\beta_n/3}\Big)\\
				\leq& \exp\Big( -\frac{1}{16+16/3}\cdot\frac{nu}{\lambda_L\beta_n}\Big)\\
				\leq& \exp\Big( -\frac{1}{32+32/3}\cdot\frac{nt}{\lambda_L\beta_n}\Big).
			\end{align*}
			This leads to a tail probability bound of $\sum_{i=1}^n G_{\beta_n}(f_{j^*},Z_i)/n$, that is
			$$P\Big\{\frac{1}{n}\sum_{i=1}^nG_{\beta_n}(f_{j^*},Z_i)>t\Big\}\leq 2\mathcal{N}_{2n}\exp\Big( -\frac{1}{43}\cdot\frac{nt}{\lambda_L\beta_n}\Big).$$
			Then for $a_n>0$,
			\begin{align*}
				\mathbb{E}_S\Big[ \frac{1}{n}\sum_{i=1}^nG_{\beta_n}(f_{j^*},Z_i)\Big]\leq& a_n +\int_{a_n}^\infty P\Big\{\frac{1}{n}\sum_{i=1}^nG_{\beta_n}(f_{j^*},Z_i)>t\Big\} dt\\
				\leq& a_n+ \int_{a_n}^\infty 2\mathcal{N}_{2n}\exp\Big( -\frac{1}{43}\cdot\frac{nt}{\lambda_L\beta_n}\Big) dt\\
				\leq& a_n+ 2\mathcal{N}_{2n}\exp\Big( -a_n\cdot\frac{n}{43\lambda_L\beta_n}\Big)\frac{43\lambda_L\beta_n}{n}.
			\end{align*}
			Choosing $a_n=\log(2\mathcal{N}_{2n})\cdot{43\lambda_L\beta_n}/{n}$, we have
			\begin{equation} \label{bound3}
				\mathbb{E}_S\Big[ \frac{1}{n}\sum_{i=1}^nG_{\beta_n}(f_{j^*},Z_i)\Big]\leq \frac{43\lambda_L\beta_n(\log(2\mathcal{N}_{2n})+1)}{n}.
			\end{equation}
			Setting $\delta=1/n$,  $\beta_n=c_1\max\{\mathcal{B},n^{1/p}\}$ and combining (\ref{bound0}), (\ref{bound1}), (\ref{bound2}) and  (\ref{bound3}), we get
			\begin{equation} \label{bound5}
				\mathbb{E}\big\{\mathcal{R}(\hat{f}_\phi)-\mathcal{R}(f^*)\big\}\leq\frac{c_2\lambda_L\mathcal{B}\log\mathcal{N}_{2n}(\frac{1}{n},\Vert\cdot\Vert_\infty,\mathcal{F}_{\phi})}{n^{1-1/p}}+ 2\big\{\mathcal{R}(f^*_\phi)-\mathcal{R}(f^*)\big\},
			\end{equation}
			where $c_2>0$ is a constant not depending on $n,d,\mathcal{B}$ and $\lambda_L$. This shows
			(\ref{bound5a}).
			
			\subsubsection*{Step 3: Bounding the covering number}
			Lastly, we will give an upper bound on the  covering number by the VC dimension of $\mathcal{F}_\phi$ through its parameters. Denote ${\rm Pdim}(\mathcal{F}_\phi)$ by the pseudo dimension of $\mathcal{F}_\phi$, by Theorem 12.2 in \cite{anthony1999}, for $2n\geq {\rm Pdim}(\mathcal{F}_\phi)$
			$$\mathcal{N}_{2n}(\frac{1}{n},\Vert \cdot\Vert_\infty,\mathcal{F}_\phi)\leq\Big(\frac{2e\mathcal{B}n^2}{{\rm Pdim}(\mathcal{F}_\phi)}\Big)^{{\rm Pdim}(\mathcal{F}_\phi)}.$$
			Besides, according to Theorem 3 and Theorem 6 in \cite{bartlett2019nearly}, there exist universal constants $c$ and $C$ such that
			$$c\cdot\mathcal{S}\mathcal{D}\log(\mathcal{S}/\mathcal{D})\leq{\rm Pdim}(\mathcal{F}_\phi)\leq C\cdot\mathcal{S}\mathcal{D}\log(\mathcal{S}).$$
			Combining the upper bound of the covering number and pseudo dimension, together with (\ref{bound5}), we have
			\begin{equation} \label{bound6}
				\mathbb{E}\big\{\mathcal{R}(\hat{f}_\phi)-\mathcal{R}(f^*)\big\}\leq c_3\lambda_L\mathcal{B}\frac{\log(n)\mathcal{S}\mathcal{D}\log(\mathcal{S})}{n^{1-1/p}}+ 2\big\{\mathcal{R}(f^*_\phi)-\mathcal{R}(f^*)\big\},
			\end{equation}
			for some constant $c_3>0$  not depending on $n,d,\lambda_L,\mathcal{B},\mathcal{S}$ and $\mathcal{D}$, which leads to (\ref{oracle}).
			This completes the proof of Lemma \ref{lemma1}.
		\end{proof}

	\subsection{Proof of Lemma \ref{lemma2}}
		The proof is straightforward under Assumption \ref{A1}. For the risk minimizer $f^*$ and any $f\in\mathcal{F}_\phi$, we have
		\begin{align*}
			R(f) -R^\phi(f^*)&=\mathbb{E}\{L(f(X),Y)-L(f^*(X),Y)\}\leq \mathbb{E}\{\lambda_L\vert f(X)-f^*(X)\vert\},
		\end{align*}
		thus $$\inf_{f\in\mathcal{F}_\phi}\{	R(f) -R^\phi(f^*)\}\leq \lambda_L\inf_{f\in\mathcal{F}_{\phi}}\mathbb{E}\vert f(X)-f^*(X)\vert=:\lambda_L\inf_{f\in\mathcal{F}_{\phi}}\Vert f-f^*\Vert_{L^1(\nu)},$$
		where $\nu$ denotes the marginal probability measure of $X$ and $\mathcal{F}_\phi=\mathcal{F}_{\mathcal{D},\mathcal{W},\mathcal{U},\mathcal{S},\mathcal{B}}$ denotes the class of feedforward neural networks with parameters $\mathcal{D},\mathcal{W},\mathcal{U},\mathcal{S}$ and $\mathcal{B}$.

			\subsection{Proof of Lemma \ref{quadapprox}}
		Invoking the proof of Lemma \ref{lemma2}, for any $f\in\mathcal{F}_\phi$, we firstly have
		\begin{align*}
			\mathcal{R}(f) -\mathcal{R}(f^*)&\leq \lambda_L\mathbb{E}\{\vert f(X)-f^*(X)\vert\}.
		\end{align*}
		Then for function $f\in\mathcal{F}_\phi$ satisfying $\Vert f-f^*\Vert_{L^\infty(\nu)}>\varepsilon_{L,f^*}$,
		\begin{align*}
			\mathcal{R}(f) -\mathcal{R}(f^*)&\leq \lambda_L\mathbb{E}\{\vert f(X)-f^*(X)\vert\}\\
			&\leq \lambda_L\mathbb{E}\big\{\frac{\vert f(X)-f^*(X)\vert^2}{\varepsilon_{L,f^*}}\big\}\\
			&\leq \frac{\lambda_L}{\varepsilon_{L,f^*}}\Vert f(X)-f^*(X)\Vert^2_{L^2(\nu)}.
		\end{align*}
		Secondly, under Assumption \ref{quadratic},
		$$\mathcal{R}(f) -\mathcal{R}(f^*)\leq \Lambda_{L,f^*} \Vert f-f^* \Vert^2_{L^2(\nu)},$$
		for any $f$ satisfying $\Vert f-f^* \Vert_{L^\infty(\nu)}\leq \varepsilon_{L,f^*}$. Combining above results, for any constant $\lambda_{L,f^*}\geq \max\{\Lambda_{L,f^*},\lambda_L/\varepsilon_{L,f^*}\}$, we have
		$$\mathcal{R}(f) -\mathcal{R}(f^*)\leq \lambda_{L,f^*} \Vert f-f^* \Vert^2_{L^2(\nu)},$$
		for any $f\in\mathcal{F}_\phi$.

		\subsection{Proof of Theorem \ref{thm2}}
		\begin{proof}
			Let $K\in\mathbb{N}^+$ and $\delta\in(0,1/K)$, define a region $\Omega([0,1]^d,K,\delta)$ of $[0,1]^d$ as $$\Omega([0,1]^d,K,\delta)=\cup_{i=1}^d\{x=[x_1,x_2,\ldots,x_d]^T:x_i\in\cup_{k=1}^{K-1}(k/K-\delta,k/K)\}.$$
			By Theorem 2.1 of \cite{shen2019deep}, for any $M, N\in\mathbb{N}^+$, there exists a function $f^*_\phi\in\mathcal{F}_{\mathcal{D},\mathcal{W},\mathcal{S},\mathcal{B}}$ with depth $\mathcal{D}=12M+14$ and width $\mathcal{W}=\max\{4d\lfloor N^{{1}/{d}}\rfloor+3d,12N+8\}$ such that $\Vert f^*_\phi\Vert_\infty \leq
			\vert f^*(0)\vert+\omega_{f^*}(\sqrt{d})$ and
			$$\vert f^*_\phi(x)-f^*(x)\vert\leq18\sqrt{d} \omega_{f^*}(NM)^{-2/d},$$
			for any $x\in[0,1]^d\backslash \Omega([0,1]^d,K,\delta)$ where $K=\lfloor N^{1/d}\rfloor^2\lfloor M^{2/d}\rfloor$ and $\delta$ is an arbitrary number in $(0,{1}/{3K}]$. Note that the Lebesgue measure of $\Omega([0,1]^d,K,\delta)$ is no more than $dK\delta$ which can be arbitrarily small if $\delta$ is arbitrarily small. Since $\nu$ is absolutely continuous with respect to the Lebesgue measure,  we have
			$$\Vert f^*_\phi-f^*\Vert_{L^1(\nu)}\leq18\sqrt{d}\omega_{f^*}N^{-2/d}M^{-2/d},$$
			and
			$$\Vert f^*_\phi-f^*\Vert^2_{L^2(\nu)}\leq384d\big\{\omega_{f^*}(N^{-2/d}M^{-2/d})\}^2.$$
			By Lemma \ref{lemma1}, we finally have
			$$\mathbb{E}\{\mathcal{R}(\hat{f}_\phi)-\mathcal{R}(f^*)\}\leq C\lambda_L\mathcal{B}\frac{\mathcal{S}\mathcal{D}\log(\mathcal{S})}{n^{1-1/p}}+18\lambda_L\sqrt{d}\omega_{f^*}(N^{-2/d}M^{-2/d}),$$
			where $C$ does not depend on $n,d,N,M,\lambda_L,\mathcal{D},\mathcal{W},\mathcal{B}$ or $\mathcal{S}$. This completes the proof of Theorem \ref{thm2}. Additionally, if Assumption \ref{quadratic} holds, by Lemma \ref{quadapprox}, we finally have
					$$\mathbb{E}\{\mathcal{R}(\hat{f}_\phi)-\mathcal{R}(f^*)\}\leq C\lambda_L\mathcal{B}\frac{\mathcal{S}\mathcal{D}\log(\mathcal{S})}{n^{1-1/p}}+384\lambda_{L,f^*}d\big\{\omega_{f^*}(N^{-2/d}M^{-2/d})\big\}^2,$$
			where $\lambda_{L,f^*}=\max\{\Lambda_{L,f^*},\lambda_L/\varepsilon_{L,f^*}\}$ is defined in Lemma \ref{quadapprox}.
		\end{proof}

		\subsection{Proof of Theorem \ref{thm3}}
		\begin{proof}
			Based on Theorem 3.1 in \cite{baraniuk2009random}, there exists a linear projector $A\in\mathbb{R}^{d_\delta\times d}$ that maps a low-dimensional manifold in a high-dimensional space to a low-dimensional space nearly preserving the distance. To be exact, there exists a matrix $A\in\mathbb{R}^{d_\delta\times d}$ such that
			$AA^T=(d/d_\delta)I_{d_\delta}$, where $I_{d_\delta}$ is an identity matrix of size $d_\delta\times d_\delta$, and
			$$(1-\delta)\vert x_1-x_2\vert\leq\vert Ax_1-Ax_2\vert\leq(1+\delta)\vert x_1-x_2\vert,$$
			for any $x_1,x_2\in\mathcal{M}.$ And it is not hard to check that
			$$A(\mathcal{M}_\rho)\subseteq A([0,1]^d)\subseteq \left[-\sqrt{\frac{d}{d_\delta}},\sqrt{\frac{d}{d_\delta}}\right]^{d_\delta}.$$
			For any $z\in A(\mathcal{M_\rho})$, define $x_z=\mathcal{SL}(\{x\in\mathcal{M}_\rho: Ax=z\})$ where $\mathcal{SL}(\cdot)$ is a set function which returns a unique element of a set. Note that if $Ax=z$, then it is not necessary that $x=x_z$. For the high-dimensional function $f^*: \mathbb{R}^{d}\to\mathbb{R}^1$, we define its low-dimensional representation $\tilde{f}^*:\mathbb{R}^{d_\delta}\to\mathbb{R}^1$ by
			$$\tilde{f}^*(z)=f^*(x_z), \quad {\rm for\ any} \ z\in A(\mathcal{M}_\rho)\subseteq\mathbb{R}^{d_\delta}.$$
			For any $z_1,z_2\in A(\mathcal{M}_\rho)$, let $x_i=\mathcal{SL}(\{x\in\mathcal{M}_{\rho}, Ax=z_i\})$. By the definition of $\mathcal{M}_\rho$, there exist $\tilde{x}_1,\tilde{x}_2\in\mathcal{M}$ such that $\vert x_i-\tilde{x}_i\vert\leq\rho$ for $i=1,2$. Then,
			\begin{align*}
				\vert \tilde{f}^*(z_1)-\tilde{f}^*(z_2)\vert&=\vert {f}^*(x_1)-{f}^*(x_2)\vert\leq\omega_{f^*}(\vert x_1-x_2\vert)\leq \omega_{f^*}(\vert \tilde{x}_1-\tilde{x}_2\vert+2\rho)\\
				&\leq \omega_{f^*}(\frac{1}{1-\delta}\vert A\tilde{x}_1-A\tilde{x}_2\vert+2\rho)\\
				&\leq \omega_{f^*}(\frac{1}{1-\delta}\vert A{x}_1-A{x}_2\vert+\frac{2\rho}{1-\delta}\sqrt{\frac{d}{d_\delta}}+2\rho)\\
				&\leq \omega_{f^*}(\frac{1}{1-\delta}\vert z_1-z_2\vert+\frac{2\rho}{1-\delta}\sqrt{\frac{d}{d_\delta}}+2\rho).
			\end{align*}
			Then, by Lemma 4.1 in \cite{shen2019deep}, there exists a $\tilde{g}$ defined on $\mathbb{R}^{d_\delta}$ which has the same type modulus of continuity as $f^*$ such that
			$$\vert \tilde{g}(z)-\tilde{f}^*(z)\vert\leq\omega_{f^*}(\frac{2\rho}{1-\delta}\sqrt{\frac{d}{d_\delta}}+2\rho)$$
			for any $z\in A(\mathcal{M}_\rho)\subseteq\mathbb{R}^{d_\delta}$. With $E=[-\sqrt{{d}/{d_\delta}},\sqrt{{d}/{d_\delta}}]^{d_\delta}$, by Theorem 2.1 and 4.3 of \cite{shen2019deep}, for any $N,M\in\mathbb{N}^+$, there exists a function $\tilde{f_\phi}: \mathbb{R}^{d_\delta}\to\mathbb{R}^1$ implemented by a ReLU FNN with width $\mathcal{W}=\max\{4d_\delta \lfloor N^{{1}/{d_\delta}}\rfloor+3d_\delta,12N +8\}$ and depth $\mathcal{D}=12M+14$ such that
			$$\vert \tilde{f}_\phi(z)-\tilde{g}(z)\vert\leq18\sqrt{d_{\delta}} \omega_{f^*}\big((NM)^{-2/d_\delta}\big)$$
			for any $z\in E\backslash \Omega(E)$ where $\Omega(E)$ is a subset of $E$ with an arbitrarily small Lebesgue measure as well as $\Omega:=\{x\in\mathcal{M}_\rho: Ax\in\Omega(E)\}$ does. In addition, for any $x\in\mathcal{M}_\rho$, let $z=Ax$ and $x_z=\mathcal{SL}(\{x\in\mathcal{M}_\rho: Ax=z\})$. By the definition of $\mathcal{M}_\rho$, there exist $\bar{x},\bar{x}_z\in\mathcal{M}$ such that $\vert x-\bar{x}\vert\leq\rho$ and $\vert x_z-\bar{x}_z\vert\leq\rho$. Then we have
			\begin{align*}
				\vert x-x_z\vert&\leq\vert\bar{x}-\bar{x}_z\vert+2\rho\leq\vert A\bar{x}-A\bar{x}_z\vert/(1-\delta)+2\rho\\
				&\leq(\vert A\bar{x}-A{x}\vert+\vert A{x}-A{x}_z\vert+\vert A{x}_z-A\bar{x}_z\vert)/(1-\delta)+2\rho\\
				&\leq(\vert A\bar{x}-A{x}\vert+\vert A{x}_z-A\bar{x}_z\vert)/(1-\delta)+2\rho\\
				&\leq 2\rho\{1+\sqrt{d/d_\delta}/(1-\delta)\}.
			\end{align*}
			
			If we define $f^*_\phi=\tilde{f}_\phi\circ A$ which is $f^*_\phi(x)=\tilde{f}_\phi(Ax)$ for any $x\in[0,1]^d$, then $f^*_\phi\in\mathcal{F}_{\mathcal{D},\mathcal{W},\mathcal{U},\mathcal{S},\mathcal{B}}$ is also a ReLU FNN with the same parameter as $\tilde{f}_\phi$, and for any $x\in\mathcal{M}_\rho\backslash\Omega$ and $z=Ax$, we have
			\begin{align*}
				\vert f^*_\phi(x)-f^*(x)\vert&\leq\vert f^*_\phi(x)-f^*_\phi(x_z)\vert+\vert f^*_\phi(x_z)-f^*(x)\vert\\
				&\leq\omega_{f^*} (\vert x-x_z\vert)+\vert \tilde{f}_\phi(z)-\tilde{f}^*(z)\vert\\
				&\leq\omega_{f^*} (\vert x-x_z\vert)+\vert \tilde{f}_\phi(z)-\tilde{g}(z)\vert+\vert \tilde{g}(z)-\tilde{f}^*(z)\vert\\
				&\leq2\omega_{f^*}(\frac{2\rho}{1-\delta}\sqrt{\frac{d}{d_\delta}}+2\rho)+18\sqrt{d_{\delta}} \omega_{f^*}\big((NM)^{-2/d_\delta}\big)\\
				&\leq (2+18\sqrt{d_{\delta}}) \omega_{f^*}\big((C_2+1)(NM)^{-2/d_\delta}\big),
			\end{align*}
			where $C_2>0$ is a constant not depending on any parameter. The last inequality follows from $\rho\leq C_2(NM)^{-2/d_{\delta}}(1-\delta)/\{2(\sqrt{d/d_\delta}+1-\delta)\}$. Since  the probability measure $\nu$ of $X$ is absolutely continuous with respect to the Lebesgue measure,
			\begin{equation}\label{approxbound}
				\Vert f^*_\phi -f^*\Vert_{L^1(\nu)} \leq (2+18\sqrt{d_{\delta}}) \omega_{f^*}\big\{(C_2+1)(NM)^{-2/d_\delta}\big\},
			\end{equation}
			where $d_\delta=O(d_\mathcal{M}{\log(d/\delta)}/{\delta^2})$ is assumed
			to satisfy $d_\delta \ll d$. By Lemma \ref{lemma1}, we have
			$$\mathbb{E}\{\mathcal{R}(\hat{f}_\phi)-\mathcal{R}(f^*)\}\leq C_1\lambda_L\mathcal{B}\frac{\mathcal{S}\mathcal{D}\log(\mathcal{S})}{n^{1-1/p}}+(2+18\sqrt{d_{\delta}}) \omega_{f^*}\big((C_2+1)(NM)^{-2/d_\delta}\big),$$
			where $C_1,C_2>0$ are constants  not depending on $n,d,\mathcal{B},\mathcal{S},\mathcal{D},\lambda_L,N$ or $M$.
			 If Assumption \ref{quadratic} also holds, then
			 \begin{equation*}
			 	\Vert f^*_\phi -f^*\Vert^2_{L^2(\nu)} \leq \Big[(2+18\sqrt{d_{\delta}}) \omega_{f^*}\big\{(C_2+1)(NM)^{-2/d_\delta}\big\}\Big],
			 \end{equation*}
		 	and
			\begin{align*}
				\mathbb{E} \{\mathcal{R}(\hat{f}_\phi)-\mathcal{R}(f^*)\}\leq& C_1\lambda_L\mathcal{B}\frac{\mathcal{S}\mathcal{D}\log(\mathcal{S})\log(n)}{n^{1-1/p}}\\
				&\qquad+\lambda_{L,f^*}\Big[(2+18\sqrt{d_{\delta}}) \omega_{f^*}\big\{(C_2+1)(NM)^{-2/d_\delta}\big\}\Big]^2,
			\end{align*}
			where $\lambda_{L,f^*}=\max\{\Lambda_{L,f^*},\lambda_L/\varepsilon_{L,f^*}\}>0$ is a constant defined in Lemma \ref{quadapprox}.
			Furthermore, if Assumption \ref{calib} also holds, then
			\begin{align*}
				\mathbb{E} \{\Delta^2(\hat{f}_\phi,f^*)\}\leq C_{L,f^*}\Bigg( C_1\lambda_L\mathcal{B}&\frac{\mathcal{S}\mathcal{D}\log(\mathcal{S})\log(n)}{n^{1-1/p}}\\
				&+\lambda_{L,f^*}\Big[(2+18\sqrt{d_{\delta}}) \omega_{f^*}\big\{(C_2+1)(NM)^{-2/d_\delta}\big\}\Big]^2\Bigg),
			\end{align*}
			where $C_{L,f^*}>0$ is a constant defined in Assumption \ref{calib}. This completes the proof of Theorem \ref{thm3}.
		\end{proof}

		\subsection{Proof of Corollary \ref{c-all}}
		We prove the Corollary for the case of deep with fixed width networks. The proofs of other cases can be carried out similarly.
		Recall  that $f^*$ is H\"older continuous of order $\alpha$ and constant $\theta$.
		Under the conditions of Theorem \ref{thm2},  for any $N, M\in\mathbb{N}^+$, the function class of ReLU multi-layer perceptrons $\mathcal{F}_{\phi}=\mathcal{F}_{\mathcal{D},\mathcal{W},\mathcal{U},\mathcal{S},\mathcal{B}}$ with depth $\mathcal{D}=12M+14$ and width $\mathcal{W}=\max\{4d\lfloor N^{1/d}\rfloor+3d,12N+8\}$, we have the excess risk of the ERM $\hat{f}_\phi$ satisfies
		\begin{equation}\label{eqc1}
		\mathbb{E} \{\mathcal{R}(\hat{f}_\phi)-\mathcal{R}(f^*)\}\leq C\lambda_L\mathcal{B}\frac{\mathcal{S}\mathcal{D}\log(\mathcal{S})\log(n)}{n^{1-1/p}}
		+18\lambda_L\theta\sqrt{d}N^{-2\alpha/d}M^{-2\alpha/d},
		\end{equation}
		for $2n \ge \text{Pdim}(\mathcal{F}_\phi)$, where $C>0$ is a constant which does not depend on $n,d,\mathcal{B},\mathcal{S},\mathcal{D},\mathcal{W},\lambda_L,N$ or $M$. Here $\lambda_L$ is the Lipschitz constant of the robust loss function $L$ under Assumption \ref{A1}.
		
		For deep with fixed width networks, given any $N\in\mathbb{N}^+$, the network width is fixed and
		$$\mathcal{W}= \max\{4d\lfloor N^{1/d}\rfloor+3d,12N+8\},$$
		where $\mathcal{W}$ does not depend on $d$ if $\mathcal{W}= 12N+8$ when $d\lfloor N^{1/d}\rfloor+3d/4\leq3N+2$. Recall that for any multilayer neural network in $\mathcal{F}_\phi$, its parameters naturally satisfy
		\begin{equation*}
			\max\{\mathcal{W},\mathcal{D}\}\leq
			\mathcal{S}\leq
			\mathcal{W}(d+1)+(\mathcal{W}^2+\mathcal{W})(\mathcal{D}-1)+\mathcal{W}+1\leq 2\mathcal{W}^2\mathcal{D}.
		\end{equation*}
		Then by plugging $\mathcal{S}\leq2\mathcal{W}^2\mathcal{D}$ and $\mathcal{D}=12M+14$ into (\ref{eqc1}), we have
		\begin{equation*}
			\mathbb{E} \{\mathcal{R}(\hat{f}_\phi)-\mathcal{R}(f^*)\}\leq C\lambda_L\mathcal{B}
			(\log n) \mathcal{W}^2(12M+14)^2\log\{2\mathcal{W}^2(12M+14)\}\frac{1}{n^{1-1/p}}
			+18\lambda_L\theta\sqrt{d}(NM)^{-2\alpha/d}.
		\end{equation*}
		Note that the first term on the right hand side is increasing in $M$ while the second term is decreasing in $M$. To achieve the optimal rate with respect to $n$, we need a balanced choice of $M$ such  that
		$$M^2\log(M)/n\approx M^{-2\alpha/d},$$
		in terms of their order. This leads to the choice of $M=\lfloor n^{(1-1/p)d/2(d+\alpha)}(\log n)^{-1}\rfloor$ and the network depth and size, where
			\begin{eqnarray*}
			\mathcal{D}&=&12\lfloor n^{(1-1/p)d/2(d+\alpha)}(\log n)^{-1}\rfloor+14, \\
			\mathcal{W}&=& \max\{4d\lfloor N^{1/d}\rfloor+3d,12N+8\}, \\
			\mathcal{S}&= &O(n^{(1-1/p)d/2(d+\alpha)}(\log n)^{-1}),
		\end{eqnarray*}
	
		%the excess risk of the empirical risk minimizer
		Combining the above inequalities, we have
		$$\mathbb{E} \{\mathcal{R}(\hat{f}_\phi)-\mathcal{R}(f^*)\}\leq \lambda_L\left\{C\mathcal{B}
		+18\theta\sqrt{d}\right\}n^{-(1-1/p)\alpha/(d+\alpha)},$$
		%	for large enough $n$,	
		for $2n \ge \text{Pdim}(\mathcal{F}_\phi)$, where $C>0$ is a constant which does not depend on $n,\lambda,\alpha$ or $\mathcal{B}$, and $C$ does not depend on $d$ and $N$ if $d\lfloor N^{1/d}\rfloor+3d/4\leq3N+2$, otherwise $C=O(d^2\lfloor N^{2/d}\rfloor)$.
		Additionally, if Assumption \ref{calib} holds, then
		$$\mathbb{E} \{\Delta^2(\hat{f}_\phi,f^*)\}\leq \lambda_LC_{L,f^*}\left(C\mathcal{B}
		+18\theta\sqrt{d}\right)n^{-(1-1/p)\alpha/(d+\alpha)},$$
		where $C_{L,f^*}>0$ is a constant defined in Assumption \ref{calib}.
		
		\section{Figures}
In this section, we show the training process for the models ``Blocks'', ``Bumps'', ``Heavisine'' and ``Dopller'' \citep{dj1994},  using the stochastic gradient algorithm Adam  \textit{Adam}  \citep{kingma2014adam} with default learning rate $0.01$ and default $\beta=(0.9,0.99)$ (coefficients used for computing running averages of gradient and its square). We set the batch size to be $n/4$ in all cases where $n$ is the size of the training data, and train the network for at least $400$ epochs until the training loss converges or becomes sufficient steady.
	
		\begin{figure}[H]
			\centering
			\begin{subfigure}{\textwidth}
				\includegraphics[width=1\textwidth]{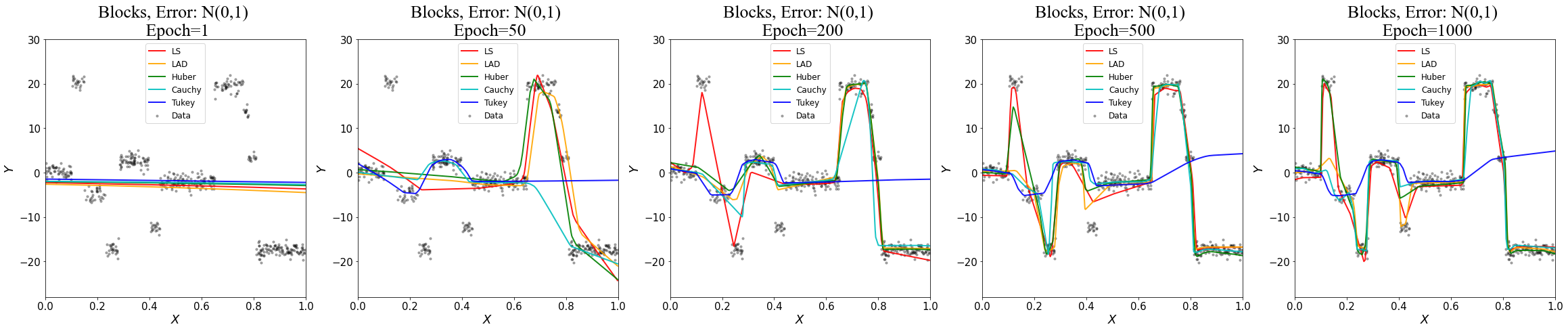}
			\end{subfigure}
			
			\begin{subfigure}{1\textwidth}
				\includegraphics[width=1\textwidth]{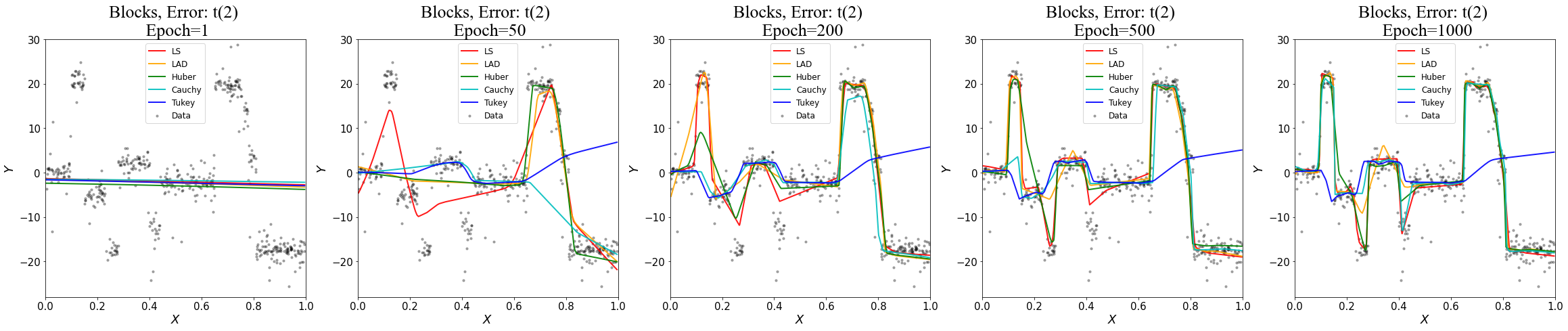}
			\end{subfigure}
			
			\begin{subfigure}{\textwidth}
				\includegraphics[width=\textwidth]{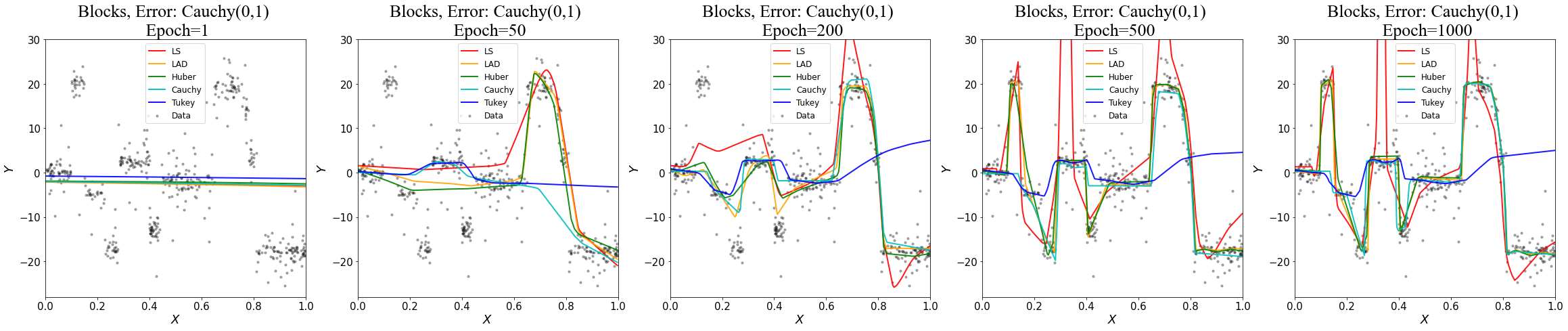}
			\end{subfigure}
		
			\begin{subfigure}{\textwidth}
			\includegraphics[width=\textwidth]{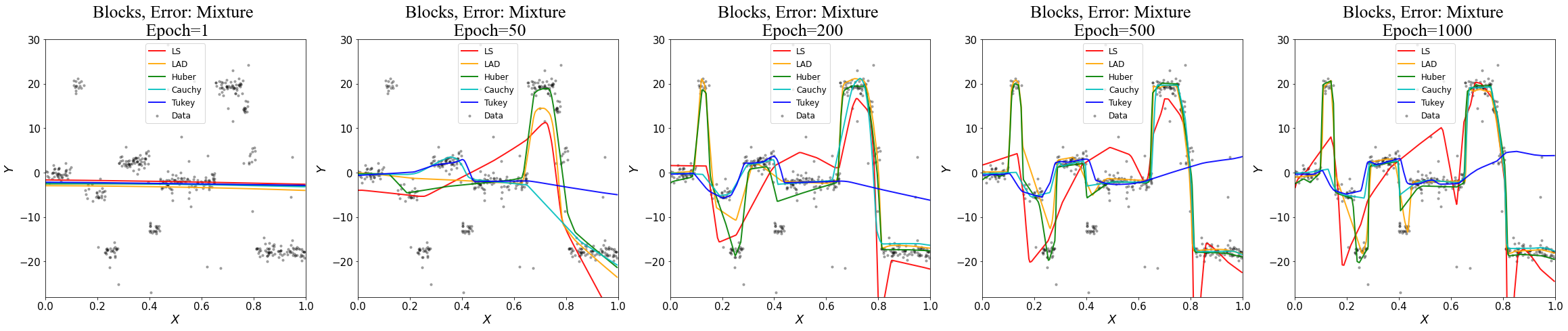}
			\end{subfigure}
			
			\caption{The training process of ``Blocks'' model. The training data is depicted as black dots and colored curves stand for predictions of the trained estimators under different loss functions at different time steps (epochs). From the top to the bottom, each row corresponds a case with a certain type of error: $N(0,1),t(2),Cauchy(0,1)$ and normal mixture. From the left to right, each column corresponds a certain point of the training process: epoch at $1,50,200,500$ and $1000$.}
			\label{fig:blcoks}
		\end{figure}

		\begin{figure}[H]
			\centering
			\begin{subfigure}{\textwidth}
				\includegraphics[width=1\textwidth]{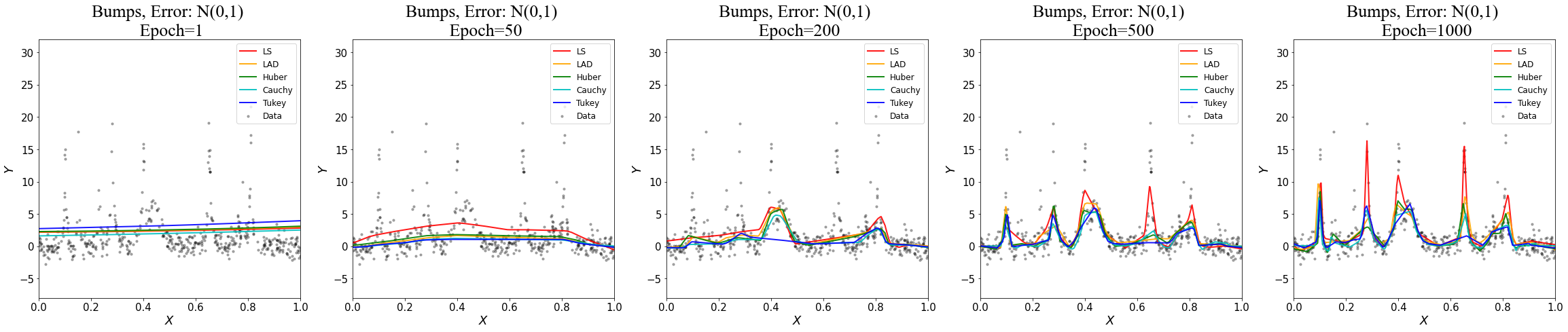}
			\end{subfigure}
			
			\begin{subfigure}{1\textwidth}
				\includegraphics[width=1\textwidth]{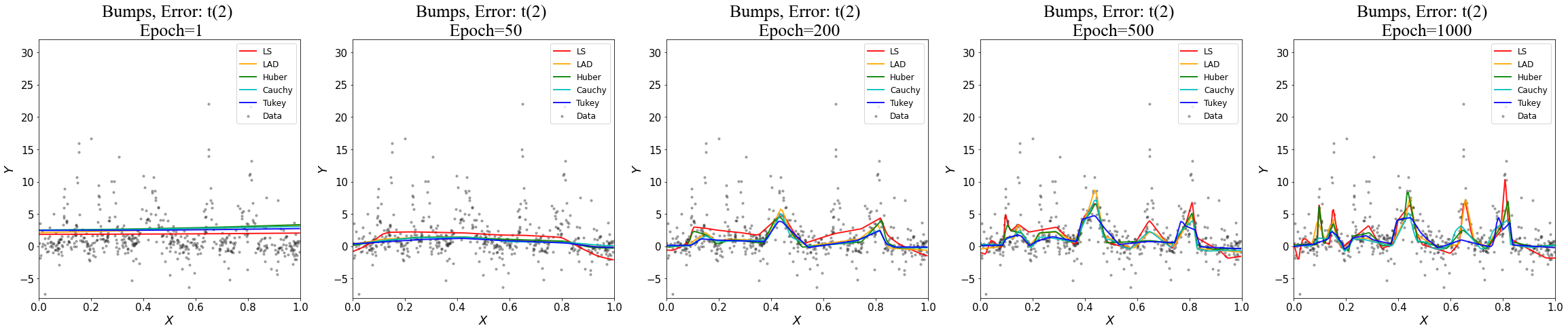}
			\end{subfigure}
			
			\begin{subfigure}{\textwidth}
				\includegraphics[width=\textwidth]{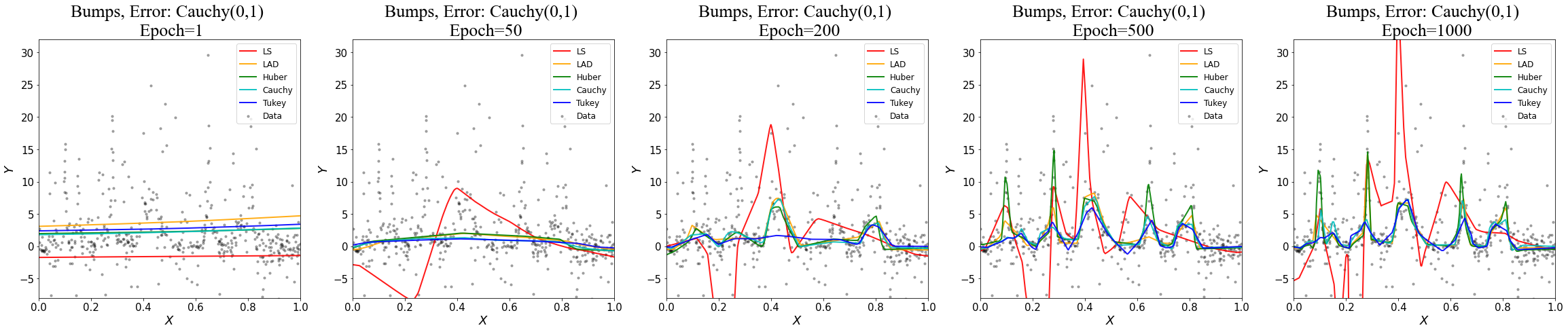}
			\end{subfigure}
			
			\begin{subfigure}{\textwidth}
				\includegraphics[width=\textwidth]{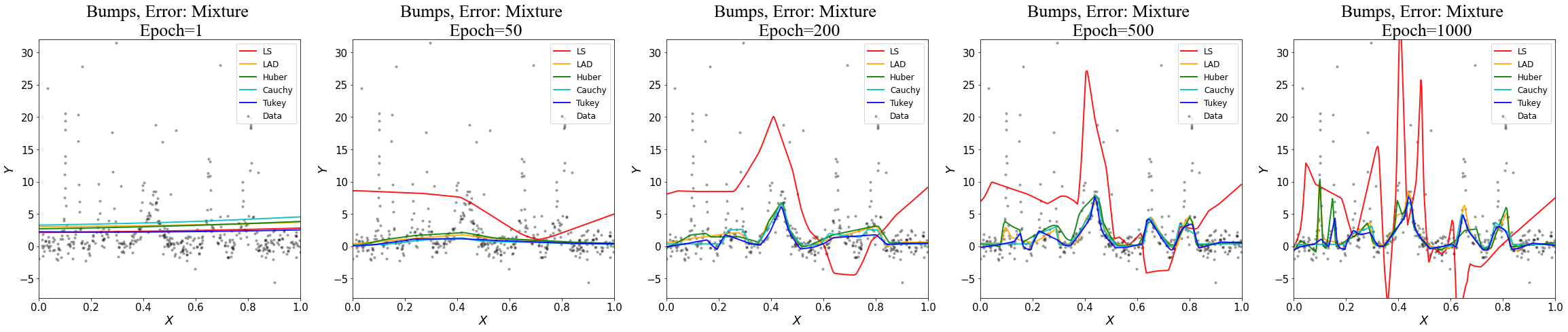}
			\end{subfigure}
		
			\caption{The training process of ``Bumps'' model. The training data is depicted as black dots and colored curves stand for predictions of the trained estimators under different loss functions at different time steps (epochs). From the top to the bottom, each row corresponds a case with a certain type of error: $N(0,1),t(2),Cauchy(0,1)$ and normal mixture. From the left to right, each column corresponds a certain point of the training process: epoch at $1,50,200,500$ and $1000$.}
			\label{fig:bumps}
		\end{figure}

		\begin{figure}[H]
			\centering
			\begin{subfigure}{\textwidth}
				\includegraphics[width=1\textwidth]{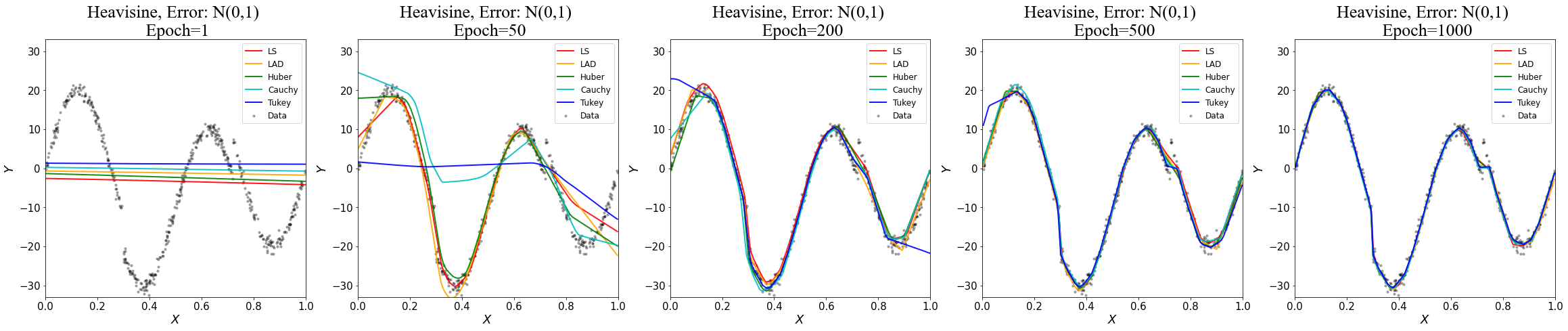}
			\end{subfigure}
			
			\begin{subfigure}{1\textwidth}
				\includegraphics[width=1\textwidth]{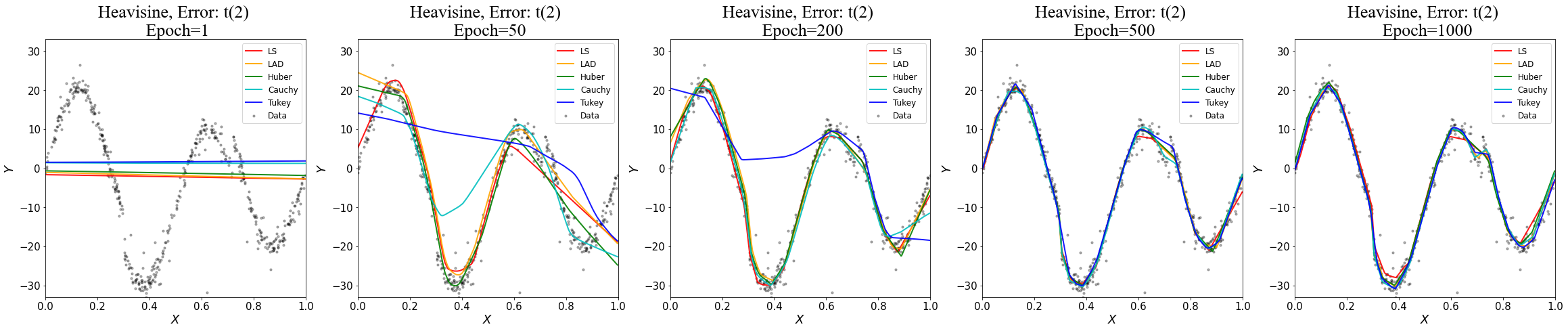}
			\end{subfigure}
			
			\begin{subfigure}{\textwidth}
				\includegraphics[width=\textwidth]{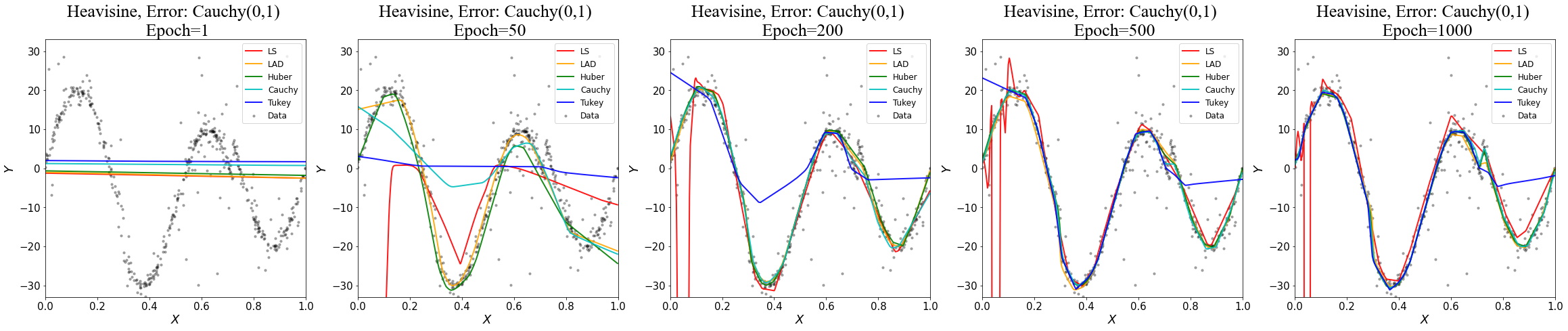}
			\end{subfigure}
			
			\begin{subfigure}{\textwidth}
				\includegraphics[width=\textwidth]{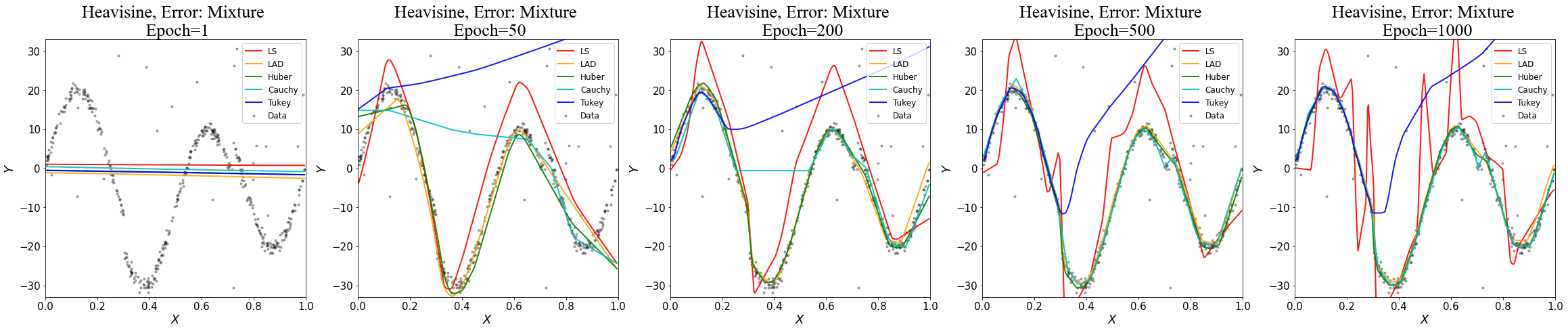}
			\end{subfigure}
		
			\caption{The training process of ``Heavisine'' model. The training data is depicted as black dots and colored curves stand for predictions of the trained estimators under different loss functions at different time steps (epochs). From the top to the bottom, each row corresponds a case with a certain type of error: $N(0,1),t(2),Cauchy(0,1)$ and normal mixture. From the left to right, each column corresponds a certain point of the training process: epoch at $1,50,200,500$ and $1000$.}
			\label{fig:heavisine}
		\end{figure}

		\begin{figure}[H]
			\centering
			\begin{subfigure}{\textwidth}
				\includegraphics[width=1\textwidth]{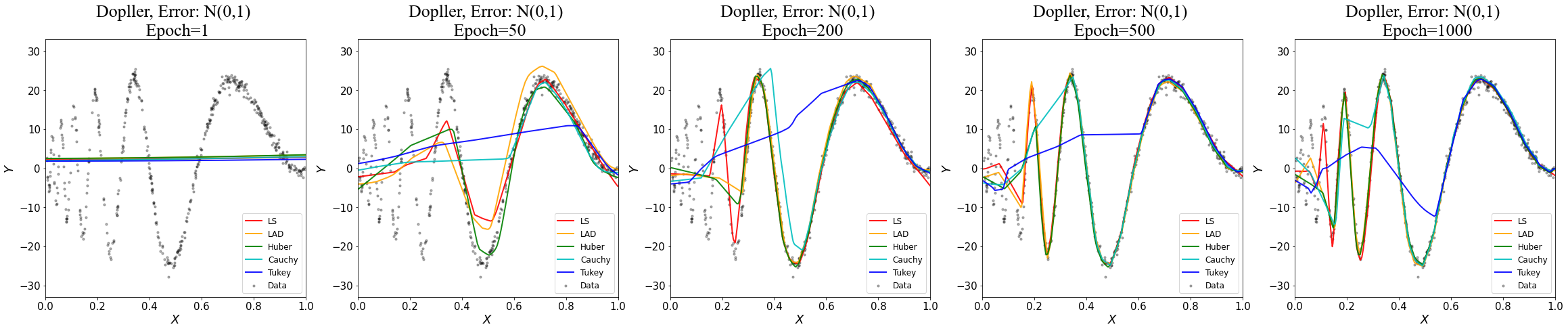}
			\end{subfigure}
			
			\begin{subfigure}{1\textwidth}
				\includegraphics[width=1\textwidth]{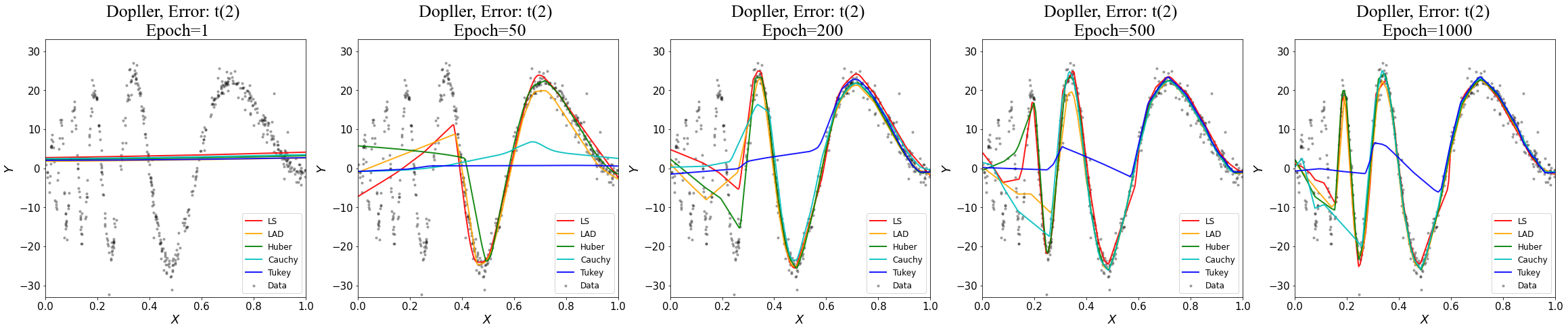}
			\end{subfigure}
			
			\begin{subfigure}{\textwidth}
				\includegraphics[width=\textwidth]{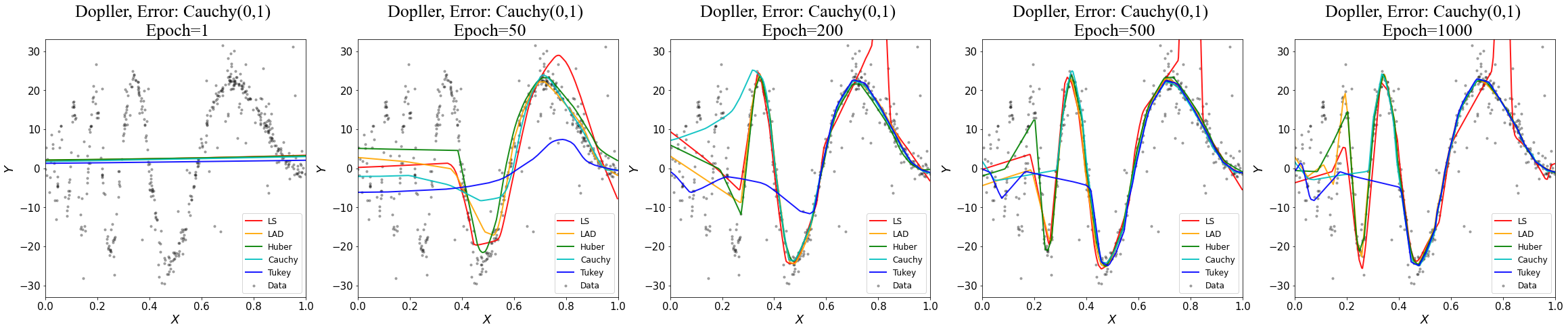}
			\end{subfigure}
		
			\begin{subfigure}{\textwidth}
			\includegraphics[width=\textwidth]{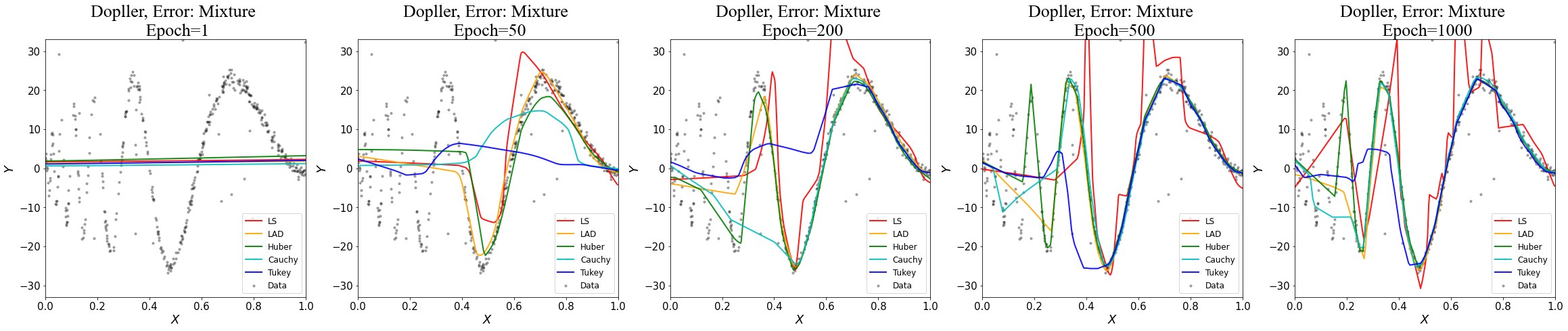}
			\end{subfigure}	
				
			\caption{The training process of ``Dopller'' model. The training data is depicted as black dots and colored curves stand for predictions of the trained estimators under different loss functions at different time steps (epochs). From the top to the bottom, each row corresponds a case with a certain type of error: $N(0,1),t(2),Cauchy(0,1)$ and normal mixture. From the left to right, each column corresponds a certain point of the training process: epoch at $1,50,200,500$ and $1000$.}
			\label{fig:doppler}
		\end{figure}

	\end{appendix}

\end{document}